\documentclass[a4paper,10pt]{amsart}

\usepackage[centertags]{amsmath}
\usepackage{amsfonts}
\usepackage{amssymb}
\usepackage{amsthm}
\usepackage{newlfont}
\usepackage{euscript}

\newcommand{\1}{1\!\!{\mathrm I}}
\renewcommand{\Re}{{\Bbb R}}
\newcommand{\eps}{\varepsilon}
\newcommand{\kap}{\varkappa}
\newcommand{\vO}{\varOmega}
\newcommand{\bu}{\mathbf{u}}
\newcommand{\bp}{\mathbf{p}}
\newcommand{\ax}{\Re^+}
\newcommand{\prt}{\partial}
\newcommand{\se}{\mathsf{E}}
\newcommand{\cs}{\mathsf{C}}
\newcommand{\Ff}{{\EuScript F}}
\newcommand{\Ef}{{\EuScript E}}
\newcommand{\Bf}{{\EuScript B}}
\newcommand{\Tf}{{\EuScript T}}
\newcommand{\Cf}{{\EuScript C}}

\newcommand{\Df}{{\EuScript D}}
\newcommand{\Nf}{{\EuScript N}}
\newcommand{\Gf}{{\EuScript G}}
\newcommand{\Qf}{{\EuScript Q}}
\newcommand{\Lf}{{\EuScript L}}

\newcommand{\Sf}{{\EuScript S}}
\newcommand{\pf}{\Pi_{fin}}
\newcommand{\Kb}{\mathbf{K}}
\newcommand{\ZZ}{{\Bbb Z}}
\newcommand{\NN}{{\Bbb N}}
\newcommand{\QQ}{{\Bbb Q}}
\newcommand{\TT}{{\Bbb T}}

\newcommand{\DD}{{\Bbb D}}

\newcommand{\demo}{\emph{Proof.} }
\newcommand{\Cd}{C_\bullet\,}
\newcommand{\Span}{\mathrm{span}\,}
\newcommand{\bro}{\hbox{{\boldmath $\rho$}}}
\newcommand{\sbro}{\hbox{{\boldmath $\rho$}}}
\newcommand{\bkap}{\hbox{{\boldmath $\vartheta$}}}
\newcommand{\bnu}{\hbox{{\boldmath $\nu$}}}
\newcommand{\bpi}{\hbox{{\boldmath $\Pi$}}}
\newcommand{\sbpi}{\hbox{{\small\boldmath $\Pi$}}}
\newcommand{\bg}{\hbox{{\boldmath $\Gamma$}}}
\newcommand{\ba}{\hbox{{\boldmath $\alpha$}}}
\newcommand{\sba}{\hbox{{\boldmath $\alpha$}}}
\newcommand{\sbg}{\hbox{{\small\boldmath $\Gamma$}}}
\newcommand{\pa}{\prt_\alpha}
\newcommand{\pba}{\prt_{\sba}}
\newcommand{\be}{\begin{equation}}
\newcommand{\ee}{\end{equation}}

\begin{document}

\numberwithin{equation}{section}
 \theoremstyle{plain}
\newtheorem{thm}{Theorem}[section]
\newtheorem{cor}{Corollary}[section]
\newtheorem{lem}{Lemma}[section]
\newtheorem{prop}{Proposition}[section]
\theoremstyle{definition}
\newtheorem{dfn}{Definition}[section]
\newtheorem{ex}{Example}[section]

\title[Stochastic calculus of variations for general L\'evy processes]
{Stochastic calculus of variations for general L\'evy processes
and its applications to jump-type SDE's  with non-degenerated
drift}
\author{Alexey M.Kulik}%
\address{Institute of Mathematics,
Ukrai\-ni\-an National Academy of Sciences, 3, Tereshchenkivska Str., Kyiv 01601, Ukraine}
 \abstract{We consider an SDE in $\Re^m$
of the type $dX(t)=a(X(t))dt+d U_t$ with a L\'evy process $U$ and
study the problem for the distribution of a solution to be regular
in various senses. We do not impose any specific conditions on the
L\'evy measure of the noise, and this is the main difference
between our method and the known methods by J.Bismut or J.Picard.
The main tool in our approach is the stochastic calculus of
variations for a L\'evy process, based on the time-stretching
transformations of the trajectories.

Three problems are solved in this framework. First, we prove that
if the drift coefficient $a$ is non-degenerated in an appropriate
sense, then the law of the solution to the Cauchy problem for the
initial equation is absolutely continuous, as soon as the L\'evy
measure of the noise satisfies one of the rather weak {intensity
conditions}, for instance the so-called \emph{wide cone
condition}. Secondly, we provide the sufficient conditions for the
density of the distribution of the solution to the Cauchy problem
to be smooth in the terms of the family of the so-called
\emph{order indices} of the L\'evy measure of the noise (the drift
again is supposed to be non-degenerated). At last, we show that an
invariant distribution to the initial equation, if exists,
possesses a $C^\infty$-density provided the drift is
non-degenerated and the L\'evy measure of the noise satisfies the
wide cone condition.}
\endabstract
\email{kulik@imath.kiev.ua}%
\subjclass[2000]{Primary 60H07; Secondary 60G51}%
\keywords{L\'evy process, admissible time-stretching
transformations, differential grid, stochastic calculus of
variations}
\maketitle

\centerline{\textsc{Introduction}}
 \vskip 10pt
\setcounter{section}{0}

 In this paper, we consider an SDE in $\Re^m$ of the type
\be\label{01}
 dX(t)=a(X(t))dt+d U_t,
\ee
 where $a\in C^1(\Re^m,\Re^m)$ satisfies the linear growth
condition and $U_\cdot$ is a L\'evy process in $\Re^m$. We study
the properties of the distribution of both the solution
$X(x,\cdot)$ to the Cauchy problem associated with (\ref{01}) and
a stationary solution $X^*(\cdot)$ to (\ref{01}), supposing latter
to exist.  The question under discussion is the following one: do
the distributions $P_{x,t}(dy)\equiv P(X(x,t)\in dy),P^*(dy)\equiv
P(X^*(t)\in dy) $ of these solutions have densities $p_{x,t},p^*$
w.r.t. the Lebesgue measure $\lambda^m$ in $\Re^m$? Do these
densities possess any additional regularity property, for
instance,  belong to the class $C^\infty$? This question is a
natural analog for the classical hypoellipticity
 problem for partial differential equations, and it can be
 reformulated in analytic terms in the following way. Let
  $L$ be the {\it L\'evy-type} pseudo-differential operator
$$
 Lf(x)= (\nabla f(x),a(x))_{\Re^m}+\int\limits_{\|u\|_{\Re^m}>1}\Bigl[f(x+u)-f(x)
 \Bigr]\Pi(du)+\int\limits_{\|u\|_{\Re^m}\leq 1}\Bigl[f(x+u)-f(x)-(\nabla f(x),u)_{\Re^m}
 \Bigr]\Pi(du)
$$
 associated with (\ref{01}), where $\Pi$ is the L\'evy measure for
 $U$. Then $P_{x,t}(dy)$ is the
  fundamental solution to the Cauchy problem for the  operator
  $\prt_t-L$ and $P^*(dy)$ is the invariant measure for the
  operator $L$.

 The hypoellipticity problem for equations
of the type (\ref{01}) and the more general equations \be\label{02}
d X(t)=a(X(t))dt+\int_{\Re^m} c(X(t-),u)\tilde \nu(dt,du)
 \ee
 with a compensated Poisson point measure $\tilde \nu$ was
 studied by numerous authors.

 First of all, let us mention the analytic approach,  see \cite{kochubey}
  and survey in \cite{kolokoltsov}. This approach uses
 some version of the {parametrix method}, and
 the typical conditions demanded here contain
 the assumptions on a smoothness and a growth rate of  the
 probability density of the initial process $U$ (roughly speaking,
 the noise should be close to the one generated by a stable process).

 There  also exist two groups of probabilistic results inspired
 by the  Malliavin's approach to the hypoellipticity problem
 in the diffusion (i.e., parabolic) setting.
 The first group is based on the method, in which a Malliavin-type calculus on the space of
the trajectories of L\'evy processes is introduced via the
transformations of trajectories that change values of their jumps.
This approach was proposed by J.Bismut (\cite{bismut}). In this
method the L\'evy measure was initially supposed to have some
(regular) density w.r.t Lebesgue measure. This is a natural
condition sufficient for the transformations, changing values of
the jumps, to be admissible. There exists a lot of works in this
direction, weakening both the non-degeneracy conditions on
coefficients and regularity claims on the L\'evy measure, cf.
\cite{bict_grav_Jac}, \cite{leandre},\cite{kom_takeuchi}. There
also exists a closely related approach based on a version of
Yu.A.Davydov's {\it stratification method}, cf. \cite{Dav_Lif},
\cite{Dav_Lif_Smor}. One can say that this group of results is
based on \emph{a spatial regularity} of the noise, which through
either Malliavin-type calculus or stratification method guarantees
the regularity of the distribution of the functional under
investigation.

Another group of results is based on the approach developed by
J.Picard, see \cite{picard} and
\cite{ishikawa},\cite{ishikawa_Kunita}. Here the  perturbations of
the point measure by adding a point into it are used. Since the single
perturbation of such a kind generates not a derivative but a
difference operator, one should use an ensemble of such
perturbations. Therefore \emph{a frequency regularity} is
needed, i.e. limitations on the asymptotic behavior of the L\'evy
measure at the origin should be imposed.

Our aim is to study  the hypoellipticity problem for equation
(\ref{01}) in a situation where the conditions imposed on the
L\'evy measure of the noise are as weak as possible. In
particular, the noise is not supposed to possess  neither spatial
nor frequency regularities.

Three problems are solved in this paper.
 The first one is
concerned with the absolute continuity of the law of the solution
to (\ref{01}) with non-degenerated drift. We give a general
sufficient condition for the absolute continuity without any
restrictions on $U$.
 The same problem was solved in
\cite{Me_TViMc},\cite{Me_jumps_UMZH} for the equation of the type
(\ref{02}) with some moment restriction on the jump part, and in
\cite{Nou_sim} for the one-dimensional SDE of the type (\ref{01}).

The second problem is to provide the conditions on the L\'evy
measure of the noise, which would be sufficient and close to the
necessary ones for the smoothness of the density of the law of
$X(x,t)$. This problem is unsolved even in the case $a=0,
c(x,u)=u$; for the L\'evy process $U$, the criterion for the
distribution of $U_t$ to possess a $C^\infty$-density is not
known.  We show that if the drift coefficient in equation
(\ref{01}) is non-degenerated in an appropriate sense, then for
the law of $X(x,t)$ such a criterion  can be given in the terms of
properly defined \emph{order indices} $\bro_r, r\in \NN$ of the
L\'{e}vy measure $\Pi$.

The claim on the drift $a$ to be non-degenerated is least
restrictive while the problem of the investigation of the properties
of the invariant distribution to (\ref{01}) is considered. Such a
claim is very natural since the invariant distribution have to
exist, and appears to be sufficient for an invariant distribution
to possess the $C^\infty$-density under very mild conditions on
the jump noise.

Our approach is motivated by a natural idea that, without any
conditions on the L\'evy measure of $U$, there always exist
admissible transformations of $U$ changing the moments of jumps,
and one can construct some kind of stochastic calculus of
variations based on these transformations. This idea is not very
new, it was mentioned in the introduction to \cite{picard}. We
also believe that it was one of the motivations for the
construction of an integration-by-parts framework for the pure
Poisson process in \cite{Carl_Pard} and \cite{TsoiEl}. However the
detailed version of the calculus of variation, based on the time
changing transformations, which would give opportunity to study
$m$-dimensional SDE's, was not available till the recent papers of
the author \cite{Me_TViMc},\cite{Me_jumps_UMZH} (the preliminary
version of such a calculus was proposed by the author in
\cite{Me_1996}; the similar approach was proposed in
\cite{Nou_sim} with an application to a one-dimensional SDE of the
type (\ref{01})).

The structure of the paper is the following. In Section 1 we
formulate the main results of the paper, in Section 2 we  make a
detailed discussion of these results and give some sufficient
conditions and corollaries. In Section 3 the stochastic calculus
for L\'evy processes, based on the time-stretching
transformations, is introduced. The proofs of the statements about
the existence of the density for $P_{x,t}(dy)$, smoothness of this
density, and smoothness of the density for $P^*(dy)$ are given in
Sections 4, 5 and 6, respectively.

\section{Main results}

\subsection{Auxiliary definitions and notation.} Before
formulating the main results of the paper, let us introduce a
notation. Denote, by $S_m=\{v\in\Re^m|\|v\|_{\Re^m}=1\}$, a unit
sphere in $\Re^m$. For $v\in S_m, \varrho\in(0,1)$, denote by
$V(v, \varrho)\equiv \{y\in\Re^m||(y,v)|_{\Re^m}\geq
\varrho\|y\|_{\Re^m}\}$ the two-sided cone with the axis $\langle
v\rangle\equiv\{tv, t\in\Re\}$.

\begin{dfn}\label{d14a} For $r\in\NN$, we define
$$
\rho_r(\varrho,\eps)= \Bigl[\varepsilon^r\ln{1\over
\varepsilon}\Bigr]^{-1}\cdot\inf_{v\in S_m}\int_{V(v, \varrho)}
(|(u, v)_{\Re^m}|\wedge \varepsilon)^r\Pi(du),\quad \eps>0,\quad
\bro_r=\lim_{\varrho\to 0+}\mathop{\lim\inf}\limits_{\eps\to
0+}\rho_r(\varrho,\eps)\in[0,+\infty].
$$
\end{dfn}

We call $\bro_r$ \emph{the upper order index of power $r$} for the
L\'evy measure $\Pi$. The main role in our considerations plays
the index $\bro_2$; we denote this index by  $\bro$.

\begin{dfn}\label{d14b}Define
$$
\vartheta(\eps)= \Bigl[\varepsilon^2\ln{1\over
\varepsilon}\Bigr]^{-1}\cdot\sup_{v\in S_m}\int_{\Re^m} (|(u,
v)_{\Re^m}|\wedge \varepsilon)^2\Pi(du),\quad \eps>0,\quad
\bkap=\mathop{\lim\inf}\limits_{\eps\to
0+}\vartheta(\eps)\in[0,+\infty].
$$
\end{dfn}

We call $\bkap$ \emph{the lower order index} for the L\'evy
measure $\Pi$. In the one-dimensional case, the definition of the
order indices is most simple, since $S_1=\{-1,+1\}$ and
$V(v,\varrho)=\Re$ for $v=\pm 1,\varrho\in(0,1)$. In the case
$m=1$, we have
$$\rho_r(\varrho,\eps)=\rho_r(\eps)=
\Bigl[\varepsilon^r\ln{1\over
\varepsilon}\Bigr]^{-1}\cdot\int_{\Re} (|u|\wedge
\varepsilon)^r\Pi(du),\quad
\vartheta(\eps)=\Bigl[\varepsilon^2\ln{1\over
\varepsilon}\Bigr]^{-1}\cdot\int_{\Re} (|u|\wedge
\varepsilon)^2\Pi(du),
$$
and  $\bkap=\bro$.

\begin{dfn}\label{d15} The function $a$ belongs to the class
$\Kb_r, r\in \NN,$ if, for every $\varrho\in(0,1)$, there exists
$D=D(a,r,\varrho)>0$ such that, for every $x\in\Re^m, \,v \in
S_m$, there exists $w=w(x,v)\in S_m$ with
\be\label{13b}\left|(a(x+y)-a(x), v)_{\Re^m}\right|\geq
D|(y,w)_{\Re^m}|^r,\quad y\in V(w,\varrho), \|y\|_{\Re^m}\in
(-D,D). \ee

The function   $a$ belongs to the class $\Kb_{r,loc}^O$ ($r\in
\NN,$ $O$ is some open subset of $\Re^m$) if, for every $x\in O,
\varrho\in(0,1)$, there exists $D=D(a,r,\varrho,x)>0$ such that,
for every $v \in S_m$, there exists   $w=w(x,v)\in S_m$ with
(\ref{13b}) being true. The function $a$ belongs to the classes
$\Kb_\infty$ or $\Kb_{\infty,loc}^O$, if $\exists r\in\NN: a\in
\Kb_r$ or $a\in \Kb_{r,loc}^O$, respectively.
\end{dfn}

\begin{ex}\label{e16} a) The function $a\in C^1(\Re^m,\Re^m)$ belongs to the class
$\Kb_{1,loc}^O$ if, for every $x\in O$, $\det \nabla a(x)\not=0$.

b) The function $a\in C^1(\Re^m,\Re^m)$ belongs to the class
$\Kb_1$ if $\sup_{x\in\Re^m}\left\|[\nabla
a(x)]^{-1}\right\|_{\Re^{m\times m}}<+\infty$ and $\nabla a$ is
uniformly continuous.

c) The function $a\in C^r(\Re, \Re)$ belongs to the class $\Kb_r$
if, for some $R,c>0$, the inequality $|a'(x)|\geq c$ holds for all
$x$ with $|x|>R$, and, for every $x$, one of the derivatives
$a'(x),a''(x),\dots, a^{(r)}(x)$ differs from $0$.
\end{ex}

\begin{dfn}\label{d310} The measure $\Pi$ satisfies {\it the wide cone
condition} if, for every $v\in S_m$, there exists
$\varrho=\varrho(v)\in (0,1)$ such that
$\Pi(V(v,\varrho))=+\infty$.
\end{dfn}

{\it Remarks. 1.} For $m=1$, the measure $\Pi$ satisfies the wide
cone condition iff $\Pi(\Re)=+\infty$.

{\it 2.}  In Definition \ref{d310}, the value of the parameter
$\rho$ can be chosen to be independent of $v$; this follows from
the compactness of $S_m$.

Denote, by $CB^k(\Re^m)$, the set of the real-valued functions $f$
on $\Re^m$ such that $f$ has $k$ Sobolev derivatives and its
$k$-th derivative is a bounded function on $\Re^m,
CB^0(\Re^m)\equiv L_\infty(\Re^m)$. Denote also, by
$C_b^\infty(\Re^m)$, the set of the real-valued infinitely
differentiable functions on $\Re^m$ that are bounded together with
every their derivative. It is clear that $CB^{k}(\Re^m)\subset
C^{k-1}(\Re^m)$ and $C_b^\infty(\Re^m)=\bigcap_{k=1}^\infty
CB^k(\Re^m)$.

\subsection{Absolute continuity of the law of $X(x,t)$.}

 In this subsection, the coefficient $a$ is supposed to belong to $C^1(\Re^m,\Re^m)$ and to
 satisfy the linear growth condition.

\begin{thm}\label{t11} Suppose that for a given $x_*\in\Re^m$
 there exists $\eps_*>0$ such that
 for
 arbitrary $v\in \Re^m\backslash\{0\}, x\in \bar B(x_*,\eps_*)\equiv
 \{y|\|y-x_*\|\leq \eps_*\}$
\be\label{11}
 \Pi\Bigl(u : (a(x+u)-a(x),  v)_{\Re^m}\not=0
\Bigr)= +\infty. \ee
 Then, for every $t>0$,
$$
P\circ[X(x_*,t)]^{-1}\ll\lambda^m.
$$
\end{thm}

This statement is analogous to that  of Theorem 3.2
\cite{Me_TViMc}, but the moment restriction analogous to condition
(\ref{13}) below, that was used in \cite{Me_TViMc}, is removed
here.

The statement of Theorem \ref{t11} can be generalized in the
following way. Consider the sequence of equations of the type
\be\label{12} X_n(x,t)=x+\int_0^ta_n(X_n(x,s))\,ds+U_t^n+V_t^n,
\quad t\in\ax, \ee
 where $V^n$ are non-random functions from the
Skorokhod's space $\DD(\Re^+,\Re^m)$, and the L\'evy processes
$U^n$ are given by stochastic integrals
$$
U_t^n=U_0+\int_0^t\int_{\|u\|>1}c_n(u)\nu(ds,du)+\int_0^t\int_{\|u\|\leq
1}c_n(u)\tilde \nu(ds,du),\quad t\in\ax, n\in\NN.
$$

\begin{thm}\label{t12} Suppose that the following conditions hold
true:

1) the coefficients $a_n, n\geq 1$ belong to $C^1(\Re^m,\Re^m)$
and satisfy the uniform linear growth condition;

2) $a_n\to a, \nabla a_n\to \nabla a, n\to+\infty, $  uniformly on
every compact set;

3) the functions $\|c_n\|$ are dominated by a function
$\mathbf{c}$ with $\int_{\Re^m} \left[\1_{\|u\|\leq 1}
\mathbf{c}^2(u)+\1_{\|u\|>1}\mathbf{c}(u)\right]\Pi(du)<+\infty$;

4) $c_n(u)\to u,n\to+\infty$ for $\Pi$-almost all $u\in\Re^m$;

5) $V^n\to V, n\to +\infty$ in $\DD(\ax,\Re^m)$;

6) $x_n\to x_*, t_n\to t>0, n\to +\infty$ and the function $V$ is
continuous at the point $t$.

Suppose also that the function $a$, the measure $\Pi$ and the
point $x_*$ satisfy the condition of Theorem \ref{t11}.

Then the laws of $X_n(x_n,t_n)$  converge in variation to the law
of the solution $X(x_*,t)$ to the equation
$$
X(x_*,t)=x_*+\int_0^ta(X(x_*,s))\,ds+U_t+V_t, \quad t\in\ax.
$$
\end{thm}
 As a corollary, we obtain the following uniform version of
 Theorem \ref{t11}.

\begin{cor}\label{p13} Suppose that the conditions of Theorem \ref{t12} hold
true. Suppose also that, for every $n\in \NN$, the function $a_n$,
the measure $\Pi_n(du)=c_n(u)\Pi(du)$, and the point $x_n$ satisfy
the condition of Theorem \ref{t11}, and $t_n>0$. Then the
family of the distributions of $X_n(x_n,t_n), n\geq 1$ is
uniformly absolutely continuous.
\end{cor}

Let us  also give a partial form of the Corollary \ref{p13},  that
is important by itself.

\begin{cor} Suppose that the condition of Theorem
\ref{t11} holds true for every $x_*\in\Re^m$. Then the map
$$
\Re^m\times (0,+\infty)\ni(x,t)\mapsto p_{x,t}\in L_1(\Re^m)
$$
is continuous, and therefore the process $X$ is strongly Feller.
\end{cor}

\subsection{Smoothness of the density $p_{x,t}$.}

In this paper, while solving the problem of the smoothness of the
density (both of the law of $X(x,t)$ and of the law of $X^*(t)$),
we restrict ourselves by the L\'evy processes satisfying the
following moment condition: \be\label{13} \int_{\|u\|_{\Re^m}\leq
1}\|u\|_{\Re^m}\Pi(du)<+\infty. \ee This supposition is crucial
for the specific form of the calculus of variations developed
below. We believe that this limitation can be removed, and the
results given below also holds true for the L\'evy processes
without any additional moment conditions. But such an expansion
should involve some more general version of the calculus of
variations, based on a "more singular" integration-by-parts
formula. This is a subject for the further investigation.

The coefficient $a$ is supposed to be infinitely differentiable
and to have all the derivatives bounded. We also suppose that
 \be\label{13a}
 \int_{\{\|u\|>1\}}\|u\|^p\Pi(du)<+\infty\hbox{ for every } p<+\infty.
 \ee
These conditions imply, in particular, that \be\label{13c}
E\sup_{s\leq t}\|X(x,s)-x\|^p<+\infty,\quad p<+\infty. \ee
Conditions on the coefficient $a$ and condition (\ref{13a}) are
technical ones and, unlike condition (\ref{13}), can be replaced
by more weak analogs in the formulation of the most of the results
given below. In order to make the exposition transparent and
reasonably short, we omit these considerations.

The main regularity  result is given by the following theorem.
Denote $\mathbf{c}(k,m)={2e\over e-1}(km+m^2+2m-2)$, $k\geq 0,
m\in\NN$.

\begin{thm}\label{t17}Let $a\in \Kb_{r}$ and $\bro_{2r}\in
(0,+\infty]$ for some $r\in \NN$. Then, for every $x\in\Re^m$ and
$t\in\ax$ with $t{\sbro_{2r}\over 2r}>\mathbf{c}(k,m)$, the
density $p_{x,t}$ belongs to the class $CB^k(\Re^m)$. In
particular, if $a\in \Kb_{r}$ and $\bro_{2r}=+\infty$ for some
$r\in \NN$, then $p_{x,t}\in C^\infty_b(\Re^m)$ for every
$t\in\ax$.
\end{thm}

The following theorem shows that the conditions given before are
rather precise. Denote, by $\Theta$, the set of $(x,t)$ such that
$P(X(x,t)\in dy)=p_{x,t}(y)dy$. We do not claim $\Theta$ to
coincide with $\Re^m\times (0,+\infty)$ and give the properties of
$p_{x,t}$ for $(x,t)\in \Theta$.

\begin{thm}\label{t18} $\mathbf{a}${\bf .}  The density $p_{x,t}$ does not
belong to $L_{r,loc}(\Re^m)$ for $t\bkap<m(1-{1\over r})$, $r>1$.

$\mathbf{b}${\bf .}  The density $p_{x,t}$ does not belong to
$C(\Re^m)$ for $t\bkap<m$.

If the condition (\ref{13}) fails, then the following analogues of
$\mathbf{a}$,$\mathbf{b}$ hold true:

$\mathbf{a1}${\bf.}  the density $p_{x,t}$ does not belong to
$L_{r}(\Re^m)$ for $t\bkap<m(1-{1\over r})$;

$\mathbf{b1}${\bf.}  the density $p_{x,t}$ does not belong to
$CB^0(\Re^m)$  for $t\bkap<m$.
\end{thm}

\subsection{Smoothness of the invariant distribution.}
Like in the previous subsection, the coefficient $a$ is supposed
to be infinitely differentiable and to have all the derivatives
bounded. The jump noise is claimed to satisfy the moment
conditions (\ref{13}), (\ref{13a}). Consider the invariant
distribution $P^*$ of (\ref{01}) or, equivalently, the
distribution of $X^*(t)$, where $X^*(\cdot)$ is a stationary
process  satisfying (\ref{01}). We suppose the invariant
distribution to exist and to have all the moments (we do not claim
this distribution to be unique).

{\it Remark.} The most simple sufficient condition here is the
claim for the drift coefficient $a$ to be "dissipative at the
infinity": \be\label{16} \exists R\in\ax, \gamma>0:\quad
(a(x),x)_{\Re^m}\leq -\gamma \|x\|_{\Re^m}^2, \quad
\|x\|_{\Re^m}\geq R.\ee Condition (\ref{16}), together with
(\ref{13a}),  guarantees both
 that $P^*$ exists and that $P^*$ has all the moments.

\begin{thm}\label{t114} Let $\Pi$ satisfy the wide cone
condition and $a\in \Kb_\infty$.

Then $P^*(dy)=p^*(y)dy$ with $p^*\in C^\infty_b(\Re^m)$.
\end{thm}

\section{Sufficient conditions, examples and  discussion}

In this section, we would like to demonstrate by a detailed
discussion the general results formulated in Theorems \ref{t11} --
\ref{t114}.

\subsection{Absolute continuity of the law of $X(x,t)$.}
Let us formulate several sufficient conditions for the condition
(\ref{11}) to hold true.  We are interested in the conditions on
the drift $a$, such that, under minimal assumptions on the jump
noise, the solution to (\ref{01}) has the absolutely continuous
distribution. Obviously, the necessary assumption here is that
$\Pi(\Re^m)=+\infty$, because otherwise the distribution of $X(t)$
has an atom.

The first condition is given in the case  $m=1$. Everywhere below
$x_*$ is used for the initial value of the solution. Denote
$N(a,y)=\{x\in \Re| a(x)=y\}$.

\begin{prop}\label{p36} Suppose that $\Pi(\Re)=+\infty$ and
there exists some $\delta_*>0$ such that
$$
\forall y\in\Re\quad \#\Bigl[N(a,y)\cap
(x_*-\delta_*,x_*+\delta_*)\Bigr]<+\infty.
$$

Then (\ref{11}) holds true, and therefore, for every $t>0$,
$$
P\circ[X(x_*,t)]^{-1}\ll\lambda^1.
$$
\end{prop}

In \cite{Nou_sim}, in the case $m=1$ only,  the law of $X(t)$ was
proved to be absolutely continuous under condition that $a(\cdot)$
is strictly monotonous at some neighborhood of $x_*$. One can see
that this condition is somewhat more restrictive than the one of
Proposition \ref{p36}. The proof of Proposition \ref{p36}, as well
as the proofs of Propositions \ref{p310}, \ref{p37} below, is
given in the subsection 4.3.

 The second sufficient condition is formulated for multidimensional case.

\begin{prop}\label{p310} Let the measure  $\Pi$ satisfy the wide cone
condition. Suppose that there exists a neighborhood $O$ of the
initial point $x_*$ such that $a\in K_{\infty,loc}^O\equiv
\bigcap_r K_{r,loc}^O$.

Then (\ref{11}) holds true, and therefore, for every $t>0$,
$$
P\circ[X(x_*,t)]^{-1}\ll\lambda^m.
$$
\end{prop}

 One can give some more precise versions of the sufficient condition in the
multidimensional case, if the structure of the drift coefficient
is specified in more details.

Define a \emph{proper smooth surface} $S\subset\Re^m$ as any set
of the type $S=\{x|\phi(x)\in L\}$, where $L$ is a proper linear
subspace of $\Re^m$ and $\phi\in C^1(\Re^m,\Re^m)$ is such that
$\det \nabla \phi(0)\not=0$ and $\phi^{-1}(\{0\})=\{0\}$.

\begin{prop}\label{p37} Suppose that one of the following group of conditions holds true:

{$\mathbf{a.}$} $a\in C^1(\Re^m, \Re^m)$, $\det \nabla
a(x_*)\not=0$ and \be\label{310} \Pi(\Re^m\backslash
S)=+\infty\quad  \hbox{for every proper smooth surface $S$;} \ee

{$\mathbf{b.}$} $a(x)=Ax, A\in \Lf(\Re^m,\Re^m)$ is non-degenerate
and \be\label{311} \Pi(\Re^m\backslash L)=+\infty\quad \hbox{for
every proper linear subspace $L\subset \Re^m$.}\ee

Then (\ref{11}) holds true, and therefore, for every $t>0$,
$$
P\circ[X(x_*,t)]^{-1}\ll\lambda^m.
$$
\end{prop}

Condition (\ref{310}) is less restrictive than the wide cone
condition introduced in Definition \ref{d310}. It  holds true, for
instance, if $\Pi(\Re^m\backslash Y)=+\infty$ for every set
$Y\subset \Re^m$, whose Hausdorff dimension does not exceed $m-1$.

Condition (\ref{311}) is close to the necessary one, this is
illustrated by the following simple example. Let (\ref{311}) fail
for some $L$, and let $L$ be invariant for $A$. Then, for $x_*\in
L$ and any $t\geq 0$, $P(X(x_*,t)\in L)>0.$ Therefore, the law of
$X(x_*,t)$ is not absolutely continuous.

Condition (\ref{311}) was introduced by M.Yamazato in the paper
\cite{yamazato}, where the problem of the absolute continuity of
the distribution of the L\'evy process was studied. This condition
obviously is  necessary for the law of $U_t$ to possess a density.
In \cite{yamazato}, some sufficient conditions were also given.
Statement 4 of the main theorem in \cite{yamazato} guarantees the
absolute continuity of the law of $U_t$ under the following  three
assumptions:

(a) condition (\ref{311}) is valid;

(b) $\Pi(L)=0$ for every linear subspace $L\subset\Re^m$ with
dim\,$L\leq m-2$;

(c)  the conditional distribution of the radial part of some
\emph{generalized polar coordinate} is absolutely continuous.

We would like to note that assumption (c) is some kind of a
"spatial regularity" assumption (in the sense we have used in
Introduction) and is crucial in the framework of \cite{yamazato}.
Without such an assumption, condition (\ref{311}) is not strong
enough to guarantee $U_t$ to possess a density, this is
illustrated by the following example.

\begin{ex}\label{e39} Let $m=2, \Pi=\sum_{k\geq 1}\delta_{z_k}$,
where $z_k=({1\over k!}, {1\over (k!)^2}), k\geq 1.$ Every point
$z_k$ belongs to the parabola $\{z=(x,y)|y=x^2\}$. Since every
line intersects this parabola at not more than  two points,
condition (\ref{311}) together with assumption (b) given before
hold true. On the other hand, for any $t>0$, it is easy to
calculate the Fourier transform  of the first coordinate $U_t^1$
of $U_t=(U_t^1,U_t^2)$ and show that
$$
\lim_{N\to +\infty} E\exp\{i 2\pi N! U_t^1\}=1.
$$
 This means that the law of $U_t^1$ is singular, and
consequently the law of $U_t$ is singular too.
\end{ex}

Due to  Proposition \ref{p37}, (\ref{311}) is the exact condition
for the \emph{linear} multidimensional equation (\ref{01})  to
possess the same regularization feature with the one given in
Introduction. We have seen that the process $U_t$ may satisfy this
condition and fail to have an absolutely continuous distribution.
However, adding a non-degenerated linear drift, we obtain  the
solution to (\ref{01}) (i.e., an Ornstein-Uhlenbeck process with
the jump noise) with the absolutely continuous distribution. At
this time, we cannot answer the question whether (\ref{311})  is
strong enough to handle the non-linear case, i.e. whether
statement {\bf a} of Proposition \ref{p37} is valid with
(\ref{310}) replaced by (\ref{311}).

\subsection{Smoothness of the density $p_{x,t}$.}

 Theorems \ref{t17},\ref{t18} allows one to
completely describe
 the regularity properties of the distribution
density of the solution to (\ref{01}) in the case $m=1$. These
properties are determined by the value of the order index $\bro$
(remind that for $m=1$ the upper order index $\bro$ coincides with
the lower order index $\bkap$), the  only possible cases here are
$\bro=+\infty,\bro=0,\bro\in(0,+\infty).$

The case of $\bro=+\infty$ is "diffusion-like", which means that
if $a\in \Kb_1$ then the density $p_{x,t}$ instantly (i.e., for
every positive $t$) becomes infinitely differentiable. The
opposite case $\bro=0$ means that the intensity of the noise is
too low to produce the regular density and for every $x\in\Re,
t\in\ax, p>1$ the density $p_{x,t}$, if exists, does not belong to
$L_{p, loc}(\Re)$.

If we compare equation (\ref{01}) with the diffusion equations, an
essentially new feature occurs in the intermediate case
$\bro\in(0,+\infty)$. On the one hand, if $a\in \Kb_1$, then we
see from Theorem \ref{t17} that there exists a sequence
$\{\mathbf{a}_k={2e(k+1)\over \sbro(e-1)}, k\geq 0\}$ such that
$p_{x,t}\in CB^{k}(\Re)$ as soon as $t>\alpha_k$. On the other
hand, $p_{x,t}\not\in CB^0(\Re)$ for $t$ small enough. We believe
that such a feature was not known before and introduce for it the
term {\it gradual hypoellipticity}.

Thus, if $m=1$ and $a\in \Kb_1$, then the only possibilities for
the law of $P_{x,t}$ are
\begin{itemize} \item $P_{x,t}$ does not have a density of the class
$\bigcup_{p>1}L_{p,loc}$ for any $t>0$ ($\bro=0$);

\item the density of $P_{x,t}$ becomes $C^k$-differentiable after
some non-trivial period of time ($\bro\in(0,+\infty)$);

\item the density of $P_{x,t}$  instantly becomes infinitely
differentiable $(\bro=+\infty)$.
\end{itemize}

 In some cases the gradual hypoellipticity feature  can be described in more
details.

\begin{prop}\label{p111} Let $m=1$ and $\bro_1<+\infty$. Let
$\Pi$ be one-sided, i.e.
$\Pi((-\infty,0))\cdot(\Pi(0,+\infty))=0$. Then the density
$p_{x,t}$ does not belong to $CB^k(\Re)$ for $t\bro_1<k+1$.
\end{prop}

For the proof of Proposition  \ref{p111} see subsection 5.1. If
the conditions of this Proposition hold true, $\bro>0$ and $a\in
\Kb_1$, then the rate of smoothness of the density is increasing
gradually: there exist two progressions $\{\mathbf{a}_k=\alpha
k+\beta\}$ and $\{\mathbf{b}_k=\gamma k+\delta\}$
($\alpha,\gamma>0$) such that $p_{x,t}\not\in CB^{k}$ while
$t<\mathbf{b}_k$, but $p_{x,t}\in CB^{k}(\Re)$ as soon as
$t>\mathbf{a}_k$.

\begin{ex}\label{e112}
 Let $\Pi=\sum_{n\geq 1}\delta_{\gamma^{-n}},
\gamma>1$, then $\bkap=\bro=\bro_1={1\over \ln \gamma}$, and the
conditions of Proposition  \ref{p111} hold true.
\end{ex}

 The gradual hypoellipticity feature can also occur  in the multidimensional
case. If $m>1$, $\bkap>0$, $\bro_{2r}<+\infty$ and $a\in \Kb_{r}$
for some $r\in \NN$, then, on the one hand, for every $k\in \NN$
$p_{x,t}\in CB^k(\Re^m)$ while $t$ is large enough, but, on the
other hand, for every $p>1$ $p_{x,t}\not \in L_{p,loc}(\Re^m)$
while $t$ is small enough.

 Let us discuss one more question related to Theorems \ref{t17},
\ref{t18}. In Theorem \ref{t18}, no specific conditions on $a$ are
imposed.
 In particular, we can take $a\equiv 0$  and  establish the properties of the
distribution of the initial  L\'evy process $U$. It is easy to see
that any condition involving the order indices cannot provide the
distribution of $U_t$ to be singular: if $\Pi(du)=\pi(u)\,du$ and
$\Pi(\Re)=+\infty$, then the distribution of $U_t$ for every $t>0$
has a density. On the contrary, due to Theorem \ref{t18}, the
condition on $\bkap$ appears to be the proper type of a necessary
condition for the distribution of $U_t$ to have \emph{a regular}
density. Take for simplicity $m=1$ and consider the property
$$
\hbox{for every $t>0$, the distribution of $U_t$ has the density
from the class } C^\infty_b(\Re). \leqno
\hbox{{\boldmath$UC^\infty_b$}}:
$$
Due to Theorem \ref{t18}, the condition $\bro=+\infty$ is
necessary for {\boldmath$UC^\infty_b$} to hold true. On the other
hand, it is known (see \cite{kallenberg},\cite{sato82}) that if
\be\label{14} \lim_{\eps\to 0+}\Bigl[\varepsilon^2\ln{1\over
\varepsilon}\Bigr]^{-1}\int_{\{|u|\leq \eps\}}
u^2\Pi(du)=+\infty,\ee then {\boldmath$UC^\infty_b$} holds true.
The conditions $\bro=+\infty$ and (\ref{14}) are in fact very
similar, since we can rewrite the first one to the form
$$
\lim_{\eps\to 0+}\left\{\Bigl[\varepsilon^2\ln{1\over
\varepsilon}\Bigr]^{-1}\int_{\{|u|\leq \eps\}}
u^2\Pi(du)+\Bigl[\ln{1\over
\varepsilon}\Bigr]^{-1}\Pi(|u|>\eps)\right\}
  =+\infty.
$$
However, the following example shows that there exists a
non-trivial gap between these two conditions.

\begin{ex}\label{e113} Let $\Pi=\sum_{n\geq 1}n\delta_{{1\over
n!}}$. Then, for every $r\in \NN$,
$$
\bro_r\geq \lim\inf_{\eps\to 0+}\left\{\Bigl[\ln{1\over
\varepsilon}\Bigr]^{-1} \Pi(|u|>\eps)\right\}\geq
\lim\inf_{N\to+\infty}{1\over \ln N!} \sum_{n\leq N-1} n \geq
\lim\inf_{N\to+\infty}{N(N-1)\over 2 N\ln N}=+\infty.
$$
This means that if the  coefficient $a$ belongs to $K_r$ for some
$r\in \NN,$ then the solution  of (\ref{01}) possesses the
$C^\infty$-density. On the other hand, for any $t>0$, one has
$$
\lim_{N\to+\infty} \Bigl|E\exp\{i 2\pi N!
U_t\}\Bigl|=\lim_{N\to+\infty}\prod_{n>N}\Bigl|\exp\{tn(e^{i2\pi
N!\over n!}-1-{i2\pi N!\over n!})\}\Bigr|=1,
$$
thus the law of $U_t$ for every $t>0$ is singular. This provides
the example of the situation where $\bro=+\infty$, but the
distribution of $U_t$ for every $t$ is essentially singular in a
sense that \be\label{15}
\lim\sup_{|z|\to+\infty}|\phi_{U_t}(z)|=1,\ee where $\phi_{U_t}$
is used for the Fourier transform of $U_t$. Moreover, this
provides the following new and interesting feature. We  say that
the L\'evy noise in Example \ref{e113} possesses some {\it hidden
hypoelipticity} (another new term) in the following sense. The law
of $U_t$ for every $t\in \ax$  is  singular due to (\ref{15}).
But, for any drift coefficient $a\in \Kb_\infty$ (that is  a
rather general non-degeneracy condition on $a$), the law of the
solution to (\ref{01}) possesses the $C^\infty$-density.
\end{ex}

\subsection{General overview} Let us  summarize  the answers on
three questions formulated at the beginning of the Introduction.
Let us formulate in a compact form some of the previous results.
We omit additional technical conditions in the formulation.

\begin{thm}\label{p112} $\mathrm {I}.$  If
$a\in\Kb_{\infty,loc}^{\Re^m}$ and $\Pi$ satisfies the wide cone
condition,  then, for every $t>0,x\in\Re^m$,
$P_{x,t}\ll\lambda^m.$

$\mathrm {II}.$ If $a\in\Kb_{\infty}$ and $\Pi$ satisfies the wide
cone condition, then $P^*(dy)=p^*(y)dy$ with $p^*\in
C^\infty_b(\Re^m)$.

$\mathrm {III.a}.$ If $a\in\Kb_{r}$ and $\bro_{2r}=+\infty$, then
$P_{x,t}(dy)=p_{x,t}(y)dy$ with $p_{x,t}\in C^\infty_b(\Re^m)$ for
every $t>0$.

$\mathrm{\phantom{III}b}.$ If $a\in\Kb_{r}$ and
$\bro_{2r}\in(0,+\infty)$, then $P_{x,t}(dy)=p_{x,t}(y)dy$ with
$p_{x,t}\in CB^k(\Re^m)$ for every $t>\mathbf{a}_k.$

$\mathrm {\phantom{III}c}.$ If $\bkap=0$, then $p_{x,t}$, if
exists, does not belong to $L_{p,loc}$ for any $t>0, p>1$.
\end{thm}

Let us note that, surprisingly, the sufficient conditions for an
invariant distribution to possess {\it smooth\,} density (the part
II. for  Theorem \ref{p112}) look like much more similar to the
sufficient conditions for $P_{x,t}$ to possess {\it some} density
(the part I.) than the conditions for $P_{x,t}$ to possess {\it
smooth} density (the part III.).

We would like to finish Section 2 with  one more remark. It is
known that the property for the distribution of the L\'evy
 process to be absolutely continuous is time-dependent: one can construct a process $U_t$
 in such a way that the law of $U_t$ is singular for $t<t_*$ and absolutely continuous for
 $t>t_*$ for some $t_*>0$ (see \cite{rubin},\cite{tucker} and   more
 recent paper \cite{sato94}). The  results given before show that such
 a feature is still valid for the solutions of equations of
 the type (\ref{01}) with   non-degenerated drift coefficient, but in a different form.
 On the one hand, the part I. of  Theorem \ref{p112} shows that the law
 of $X(x,t)$ is absolutely continuous for every $t>0$ as soon as $a\in\Kb_{\infty,loc}^{\Re^m}$
 and $\Pi$ satisfies the wide cone condition. Thus
 the type of the distribution of $X(x,t)$,
 unlike the one of the distribution of $U_t$, is not
 time-dependent. The proper form of such a dependence is the
 "gradual hypoellipticity" feature. Recall that
 such a feature occurs when
 $\bkap >0, \bro_{2r}<+\infty$ and  $a\in \Kb_{r}$ for some $r\in \NN$.

Another form of such a dependence is given by parts II., III. of
Theorem \ref{p112}, that show that  the regularity  properties of
the distribution density of the stationary solution essentially
differ from those of the solution to the Cauchy problem. The
stationary solution can be informally considered as the solution
to the Cauchy problem with the initial point $-\infty$. Thus one
should conclude that while any finite time interval in the case
$\bkap=0$ is "not long enough" for a non-degenerated drift to
generate a smooth density,
 the infinite time interval is "long
enough", provided that $a$ is weakly non-degenerated
($a\in\Kb_\infty$) and $\Pi$ satisfies the wide cone condition.
These considerations show that the hypoellipticity properties of
the solution to (\ref{01}), in general, are essentially
time-dependent.

\section{Time-stretching transformations and
associated stochastic calculus for a  L\'evy process}

\subsection{Basic constructions and definitions.}
In this subsection we  introduce the stochastic calculus on a
space of trajectories of the general L\'evy process, that is the
basic tool in our approach. This calculus is based on the
time-stretching transformations of the jump noise and associated
differential structure. Differential constructions of  a similar
kind have been known for some time, say, the integration-by-parts
framework for a pure Poisson process was introduced independently
in \cite{Carl_Pard} and \cite{TsoiEl}, some analytic properties of
the corresponding differential structure on a configuration space
(over $\ax$ or a Riemannian manifold) were described in a cycle of
the papers by N.Privault, cf.
\cite{Privault_SDE},\cite{Privault99},\cite{Privault2001}. Our
construction (introduced initially in \cite{Me_1996}) is slightly
different and is applicable in the general situation where a
spatial variable of the noise is non-trivial. The more detailed
exposition, as well as some related notions, such as the joint
stochastic derivative and the extended stochastic integral w.r.t.
the compensated Poisson point measure, can be found  in
\cite{Me_TSP2001}.

Let us introduce the notation. By $\nu$ and $\tilde \nu$, we denote
the point measure and the compensated point measure, involved in
the L\'evy---Khinchin representation for the process $U$:
$$
U_t=U_0+\int_0^t\int_{\|u\|>1}u\nu(ds,du)+\int_0^t\int_{\|u\|\leq
1}u\tilde \nu(ds,du),
$$
$\nu$ is a Poisson point measure on
$\ax\times(\Re^m\backslash\{0\})$ with the intensity measure
$dt\Pi(du)$, $\tilde \nu(dt,du)=\nu(dt,du)-dt\Pi(du)$.
  We use the
standard terminology from the theory of Poisson point measures without any additional discussion.
The term "(locally finite) configuration" for a realization of the
point measure is frequently used. We suppose that the basic
probability space $(\Omega,{\Ff},P)$ satisfies condition
${\Ff}=\sigma(\nu)$, i.e. every random variable is a functional of
$\nu$ (or $U$). This means that in fact one can treat $\Omega$ as
the configuration space over $\ax\times(\Re^m\backslash\{0\})$
with a respective $\sigma$-algebra. Also the notion of the point
process $p(\cdot)$ associated with the process $U$ (and the
measure $\nu$) is used in the exposition. The domain $\Df$ of this
process is equal to the (random) set of $t\in \ax$ such that
$U_t\not=U_{t-}$, and $p(t)=U_t-U_{t-}$ for $t\in\Df$.

The notation $\nabla_x$ for the gradient w.r.t. the space variable
$x$ is frequently used. If the function depends only on $x$, then
the subscript $x$ is omitted. If it does not
cause misunderstanding, we omit the subscript and write, for
instance, $\|x\|$ instead of $\|x\|_{\Re^m}$.

Denote $H=L_2(\Re^+), H_0=L_\infty(\Re^+)\cap L_2(\Re^+),
Jh(\cdot)=\int_0^\cdot h(s)\,ds,h\in H.$ For a fixed $h\in H_0$, we
define  the family $\{T_h^t,t\in\Re\}$ of transformations of the
axis $\Re^+$ by putting $T^t_hx, x\in\ax$ equal to the value at
the point $s=t$ of the solution of the Cauchy problem
\be\label{21} z'_{x,h}(s)=Jh(z_{x,h}(s)),\quad s\in \Re, \qquad
z_{x,h}(0)=x.\ee
 Since (\ref{21}) is the Cauchy problem for the
time-homogeneous ODE, one has that $T^{s+t}_h=T^s_h\circ T^t_h$,
and in particular $T_h^{-t}$ is the inverse transformation to
$T_h^t$. Multiplying $h$ by some $a>0$, we multiply, in fact, the
symbol of the equation by $a$. Now, taking the time change $\tilde
s={s\over a}$, we see that $T_h^a=T_{ah}^1,a>0$, which together
with the previous considerations gives that $T_h^t=T_{th}^1, h\in
H_0, t\in\Re$.

Denote $T_h\equiv T_h^1$, we have just proved that $T_{sh}\circ
T_{th}=T_{(s+t)h}$. This means that $\Tf_h\equiv \{T_{th}, t\in
\Re\}$ is a one-dimensional group of transformations of the time
axis $\ax$. It follows from the construction that ${d\over dt}
|_{t=0}T_{th}x=Jh(x), x\in\ax.$

{\it Remark.} We call $T_h$ the \emph{time stretching
transformation} because, for $h\in C(\ax)\cap H_0$, it can be
constructed in a more illustrative way: take the sequence of
partitions $\{S^n\}$ of $\ax$ with $|S_n|\to 0, n\to +\infty$. For
every $n$, we make the following transformation of the axis: while
preserving an initial order of the segments, every segment of the
partition should be stretched by $e^{h(\theta)}$ times, where
$\theta$ is some inner point of the segment (if $h(\theta)<0$ then
the segment is in fact contracted). After passing to the limit
(the formal proof is omitted here in order to shorten the
exposition) we obtain the transformation $T_h.$ Thus one can say
that $T_h$ performs the stretching of every infinitesimal segment
$dx$ by $e^{h(x)}$ times.

 Denote $\Pi_{fin}=
 \{\Gamma\in \Bf(\Re^d),\Pi(\Gamma)<+\infty\}$ and
define, for $h\in H_0, \Gamma\in \pf$, a transformation $T_h^\Gamma$
of the random measure $\nu$ by
$$
[T_h^\Gamma \nu]([0,t]\times \Delta)=\nu([0,T_{h}t]\times
(\Delta\cap\Gamma))+ \nu([0,t]\times (\Delta\backslash\Gamma))
,\quad t\in\Re^+,\Delta\in \pf.
$$
 An easy calculation gives that $r_h(t)\equiv {d\over dt}(T_ht)=
 \int_0^1 h(T_{sh}t)\,ds,t\in\Re^+$. We put
 $$
  p_h^\Gamma=\exp\left\{\int_{\Re^+}
 r_h(t)\nu(dt,\Gamma)-\lim_{t\to+\infty}[T_ht-t]\Pi(\Gamma)\right\}.
 $$
Since $T_h^\Gamma \nu$ is again a random Poisson point measure,
its intensity measure can be expressed through $r_h(\cdot),\Pi$
explicitly. Thus the following statement is a corollary of the
classical absolute continuity result for L\'evy processes, see
\cite{Skor_nez_prir}, Chapter 9.

\begin{lem}\label{l21} The transformation $T_h^\Gamma$ is admissible for
the distribution of $\nu$ with the density $p_h^\Gamma$, i.e., for
every $\{t_1,\dots,t_n\}\subset\ax,
\{\Delta_1,\dots,\Delta_n\}\subset \pf$ and the Borel function
$\phi:\Re^n\to \Re$,
$$
E\phi( [T_h^\Gamma \nu]([0,t_1]\times \Delta_1),\dots,[T_h^\Gamma
\nu]([0,t_n]\times \Delta_n))=Ep_h^\Gamma\phi(\nu ([0,t_1]\times
\Delta_1),\dots,\nu([0,t_n]\times \Delta_n)).
$$
\end{lem}

 The statement of the lemma and the fact that  ${\Ff}$ is
generated by $\nu$ imply that the transformation $T_h^\Gamma$
generates the corresponding transformation of the random
variables, we denote it also by $T_h^\Gamma$.

 The image of a configuration of the point measure $\nu$ under
 $T_h^\Gamma$
 can be described in a following way: every  point
$(\tau,x)$ with $x\not\in \Gamma$ remains unchanged; for every
point $(\tau,x)\in N$ with $x\in \Gamma$, its ``moment of the jump"
$\tau$ is transformed to $T_{-h}\tau$; neither any point of
the configuration is eliminated nor any new point is added to
the configuration. In a sequel, we suppose that the probability space
$\Omega$ coincides with the space of locally finite configurations
on ${\Re^+}\times \Re^d$ and denote, by the same symbol
$T_h^\Gamma$, the bijective transformation of this space described
above.

Let $\Cf$ be the set of functionals $f\in \cap_pL_p(\Omega,P)$
satisfying the following condition: for every $\Gamma\in \pf$,
there exists the random element $\nabla^\Gamma_H f\in\cap_p
L_p(\Omega,P,H)$ such that, for every $h\in H_0$, \be
(\nabla_H^\Gamma f, h)_H=\lim_{\eps\to 0} {1\over \eps}[T_{\eps
h}^\Gamma\circ f-f] \ee with convergence in every $L_p,
p<+\infty$.

\begin{ex}\label{e22}  Let
$\Delta\in\pf, f=\tau_n^\Delta\equiv \inf\{t|\nu([0,t]\times
\Delta)=n\}$. Then $f\in \Cf$ and
$$
[\nabla_H^\Gamma f](\cdot)=-\1_{[0,\tau_n^\Delta]}(\cdot)
\1_{p(\tau_n^\Delta)\in\Gamma}.
$$

 We denote
$$ (\rho^\Gamma,h)= -\int_0^\infty
 h(t)\,\tilde\nu(dt,\Gamma)
$$
and note that $L_p-\lim_{\eps\to 0}{p^\Gamma_{\eps h}-1\over
\eps}=-(\rho^\Gamma,h), \quad p\in (1,+\infty). $
\end{ex}

\begin{lem}\label{l23} For every $\Gamma\in \pf$, the pair
$(\nabla_H^\Gamma,\Cf)$ satisfies  the following conditions:

 1) For every $f_1,\ldots,f_n\in\Cf$ and $F\in C^1_b(\Re^n)$,
$$
F(f_1,\ldots,f_n)\in\Cf\quad \hbox{ and }\quad \nabla_H
F(f_1,\ldots,f_n)=\sum_{k=1}^nF'_k(f_1,\ldots,f_n)\nabla_H f_k
$$
(chain rule).

 2) The map $\rho^\Gamma:h\mapsto (\rho^\Gamma,h)$ is a weak random element
 in $H$ with weak moments of all orders, and
 $$
 E(\nabla_H^\Gamma f,h)_H=-Ef(\rho^\Gamma,h),\quad h\in H,f\in\Cf
$$
(integration-by-parts formula).

3) There exists a countable set $\Cf_0\subset \Cf$ such that
$\sigma(\Cf_0)=\Ff$.
\end{lem}

Conditions 1),2) follow from the definition of the class $\Cf$ and Lemma
\ref{l21}; condition 3) holds true due to Example \ref{e22}.

For a given $h\in H, \Gamma\in \pf, p>1$, consider the map
$$
\nabla_h^\Gamma:\Cf \ni f\mapsto (\nabla_H^\Gamma f, h)_H\in
L_p(\Omega, \Ff,P)
$$
as a densely defined unbounded operator.  Lemma \ref{l23} provides
that its adjoint operator is well defined on $\Cf\subset
L_q(\Omega, \Ff,P), {1\over p}+{1\over q}=1$, by the equality
$$
[\nabla_h^\Gamma]^*g=-(\rho^\Gamma,h) g- \nabla_h^\Gamma g.
$$
Since $\Cf$ is dense in $L_q(\Omega, \Ff,P)$, this means that
$\nabla_h^\Gamma$ is closable in the $L_p$ sense.

\begin{dfn}\label{d24} The closure $D^\Gamma_{h,p}$ of
$\nabla_H^\Gamma$ in the $L_p$ sense is  called the \emph{stochastic
derivative} in the direction $(h,\Gamma)$ of order $p$. The
 $\Gamma$-stochastic derivative $D^\Gamma_{p}$ of
order $p$ is defined for $f\in\cap_{h\in H}
Dom(D^\Gamma_{h,p})$ such that there exists $g\in
L_p(\Omega,P,H)$ with
$$
(g,h)_H=D_{h,p}^\Gamma f,\quad h\in H,
$$
by the equality $D_p^\Gamma f=g$. If $p=2$, then $p$ is omitted in
the notation.
\end{dfn}

Now a differential structure on the initial  space of trajectories
is constructed, and it is natural to try to develop some calculus
which would provide statements of the type "if for a functional
$f$ the family $\{D^\Gamma f, \Gamma\in\pf\}$ is non-degenerate in
some sense, then the law of $f$ is regular." The stratification
method or the Malliavin-type calculus of variations is supposed to
be a natural tool here. However, the differential structure
developed before has some new specific properties that does not
allow us to apply these tools immediately. The most important
feature is illustrated by the following example.

\begin{ex}\label{e25} Let $f=\tau_n^\Gamma$, $h,g\in C_b(\ax)\cap
L_2(\ax)$ be such that $h(t)\int_0^tg(s)\,ds\not=g(t)\int_0^t
g(s)\, ds, t>0$, then
$$
D^\Gamma_h D^\Gamma_g f=h(\tau^\Gamma_n)\int_0^{\tau_n^\Gamma}
g(s)\,ds\not= g(\tau^\Gamma_n)\int_0^{\tau_n^\Gamma}
h(s)\,ds=D^\Gamma_g D^\Gamma_h f
$$
almost surely. In particular, this means that the family of
transformations $\{T^\Gamma_h, h\in H_0\}$ is not commutative and
therefore cannot be considered as an infinite-dimensional additive
group of transformations. Roughly speaking, the differential
structure described by $\Gamma$-stochastic derivative is
\emph{non-flat}.
\end{ex}

One possible way to overcome this difficulty and to introduce an
analog of the stratification method in the framework described
before was developed in \cite{Me_TViMc}. There, some transformation
(corresponding to the transformation of the L\'evy process into
the associated point process), that changes the non-flat gradient
$D^\Gamma$ to some linear-type gradient over a space $\Re^\infty$,
was used. The relation between these two gradients is close to the
one between  the "damped" and "intrinsic" gradients on the
configuration space over the Riemannian manifold (see
\cite{Privault2001}).

 The analysis based on the change of the space and the
gradient allows one to apply the stratification method and obtain
efficient conditions for the absolute continuity of the
distribution of a solution to (\ref{01}) or (\ref{02}). However,
this analysis appears to be rather complicated. Below we introduce
another approach based on the new notion of a \emph{differential
grid}. This approach not only simplifies the way the
stratification method can be applied,  but also allows us to
develop the efficient stochastic calculus of variations and
consider the question of the smoothness of the density.

\subsection{Differential grids and associated Sobolev
classes.}

\begin{dfn}\label{d26}
 A family $\Gf=\{[a_i,b_i)\subset \ax, h_i\in
H_0, \Gamma_{i}\in\pf, i\in\NN\}$ is called  \emph{a differential
grid} (or simply \emph{a grid}) if

(i) for every $i\not= j$,  $\Bigl([a_i,b_i)\times \Gamma_i\Bigr)
\cap \Bigl( [a_j,b_j)\times \Gamma_j\Bigr) =\emptyset$;

(ii) for every $i\in\NN$, $Jh_i>0$ inside $(a_i,b_i)$ and $Jh_i=0$
outside $(a_i,b_i)$.
\end{dfn}

Any grid $\Gf$ generates a partition of some part of the phase space
$\ax\times (\Re^m\backslash \{0\})$ of the random measure $\nu$
into the cells $\{\Gf_{i}=[a_i,b_i)\times \Gamma_{j}\}$. We call the
grid $\Gf$ \emph{finite}, if $\Gf_i=\emptyset$ for all indices
$i\in\NN$ except some finite number of indices.

Denote $T_t^{i}=T_{th_{i}}^{\Gamma_{i}}$. For any $i\in
\NN, t,\tilde t\in \Re $, the transformations $T_t^{i}$,$T_{\tilde t}^{
i}$ commute because so do the time axis transformations $T_{t
h_i}$,$T_{\tilde t h_i}$. It follows from the construction of the
transformations $T_h^\Gamma$ that, for a given $i\in \NN, t\in \Re$,
$$
T_t^{i}\tau_n^{\Gamma_{i}}=T_{t h_i}\tau_n^{\Gamma_{i}}\quad
\begin{cases}
=\tau_n^{\Gamma_{i}},& \tau_n^{\Gamma_{i}}\not\in [a_i,b_i)\\
\in [a_i,b_i),& \tau_n^{\Gamma_{i}}\in [a_i,b_i)
\end{cases}\quad \hbox{for every } n
$$
(see Example \ref{e22} for the notation $\tau_n^\Gamma$). In other
words, $T_t^{i}$ does not change points of configuration outside
the cell $\Gf_{i}$ and keeps the points from this cell in it.
Therefore, for every $i,\tilde i \in \NN, t,\tilde t\in \Re $, the
transformations $T_t^{i}$,$T_{\tilde t}^{\tilde i}$ commute, which
implies the following proposition. Denote, by $\ell_0\equiv
\ell_0(\NN)$, the set of all sequences $l=\{l_{i}, i\in \NN\}$ such
that $\#\{i| l_{i}\not=0\}<+\infty$.

\begin{prop}\label{p27} For a given grid $\Gf$ and $l\in \ell_0$,
define the transformation
$$
T^{\Gf}_{l}=T^{1}_{l_{1}}\circ T^{2}_{l_{2}}\circ\dots .
$$
This definition is correct since the  transformation
$T^{i}_{l_{i}}$ differs from the identical one only for a finite
number of indices $i$. Then $\Tf^{\Gf}=\{T^\Gf_l, l\in\ell_0\}$ is
the group of admissible transformations of $\Omega$ which is
additive in the sense that $T^\Gf_{l_1+l_2}=T^\Gf_{l_1}\circ
T^\Gf_{l_2}, l_{1,2}\in\ell_0.$
\end{prop}

It can be said that, by fixing the grid $\Gf$, we choose, from the
whole variety of admissible transformations $\{T_h^\Gamma, h\in
H_0,\Gamma\in \pf\}$, the additive family that is more convenient
to deal with. Let us introduce the gradients and Sobolev classes
associated with such families.

Denote, by $\ell_2$,
 the Hilbert space of the sequences
$$
l=\{l_{i},i\in\NN\}:\quad \|l\|_{\ell_2}\equiv\left[\sum_{i\in
\NN} l_{i}^2\right]^{1\over 2}<+\infty,\quad (l,\tilde
l)_{_{\ell_2}}\equiv \sum_{i\in \NN}\sigma_{i} l_{i}\tilde l_i.
$$
Define $\mathbf{l}^{i}\in\ell_2,i\in\NN$  by
$\mathbf{l}^{i}_{i}=1, \mathbf{l}^{i}_{j}=0, i\not=j.$

\begin{dfn}\label{d28} The random element $f\in L_p(\Omega,P,E)$
$(p\in(1,+\infty)$), taking values in a separable Hilbert space
$E$, belongs to the domain of the stochastic derivative
$D^{\Gf}_p$ if

1) for every $i\in \NN, e\in E$, $(f,e)_E\in
Dom(D^{\Gamma_i}_{h_i,p})$;

2) there exists $g\in L_p(\Omega, P,\ell_2\otimes E)$ such that
$D^{\Gamma_i}_{h_i,p}(f,e)_E=(g,{\mathbf{l}^{i}})_{\ell_2\otimes
E}$ for  $e\in E, i\in\NN$.\\
The element $g$ is denoted by $D^{\Gf}_p f$. If $p=2$, then $p$ is
omitted in the notation.
\end{dfn}

The class of all elements $f\in L_p(\Omega,P,E)$ stochastically
differentiable in the  sense of Definition \ref{d28} is denoted by
$W_p^1(\Gf,E)$. This class  is a Banach space w.r.t. the norm
$$
\|f\|_{p,1}^{\Gf,E}\equiv \left\{E \|f\|^p_E + E \Bigl\|D^{\Gf}_p
f\Bigr\|_{\ell_2\otimes E}^p\right\}^{1\over p}
$$
since the operator $D^{\Gf,E}_p$ is closed in $L_p$.

Similarly, define the Sobolev class
$W_p^d(\Gf,E)$ for $d\geq 1, p>1$ as the domain of the operator $[D^{\Gf}_p]^kf$, it
is  a Banach space w.r.t. the norm
$$
\|f\|_{p,d}^{\Gf,E}\equiv \left\{E \|f\|_E^p + E\sum_{k=1}^d
\Bigl\|[D^{\Gf}_p]^kf\Bigr\|_{[\ell_2]^{\otimes k}\otimes
E}^p\right\}^{1\over p}.
$$

At last, define $I^{\Gf}_p$ as the adjoint operator to
$D^{\Gf}_p$.   This operator is called the \emph{stochastic
integral}, which is natural, in particular, due to the following
example (see also \cite{Me_TSP2001}, Theorems 1.1 and 1.2).

\begin{ex}\label{e29} It follows from Lemma \ref{l23}  that a
non-random element $\mathbf{l}^i\in\ell_2$ belongs to the domain
of every $I_p^{\Gf}$, and
$$I_p^{\Gf}(\mathbf{l}^i)=-\rho_{h_i}^{\Gamma_i}=\int_{(a_i,b_i)\times \Gamma_i} h_i(s)\tilde \nu(ds,du).
$$
\end{ex}

The following properties of $D^{\Gf}_p$,$I^{\Gf}_p$ are due to the
chain rule (Lemma \ref{l23}, statement 1). The proof is analogous
to the proof of the same properties of the stochastic derivative
and integral w.r.t. the Wiener process and is omitted.

\begin{lem}\label{l210} 1) Let $f_j\in W_p^1(\Gf,E_j),j=1,\dots,n$,
$F:E_1\times\dots\times E_n\to E$ be Frechet differentiable,
continuous, and bounded together with its derivative. Then
$F(f_1,\dots,f_n)\in W_p^1(\Gf,E)$ and
$$
D_p^{\Gf}F(f_1,\dots,f_n)=\sum_{j=1}^nF'_j(f_1,\dots,f_n)\cdot
D_p^{\Gf}f_j.
$$
2) Let $g\in Dom(I^{\Gf}_{p_1}), f\in W_{p_2}(\Gf,\Re)$,
$p_1>p_2$. Then $fg\in Dom(I^{\Gf}_{p}),$ where $p={p_1q_2\over
p_1q_2-p_1-q_2}, q_2={p_2\over p_2-1},$ and
$$
I_p^{\Gf}(fg)=f\cdot I_{p_2}^{\Gf}(g)-(D^{\Gf}_{p_1}f,g)_{\ell_2}.
$$
\end{lem}

\subsection{Existence of the
density via the stratification method.}

In this subsection, we give two sufficient conditions for the existence
of the density for a functional on $(\Omega,\Ff,P)$.  The first
condition is formulated in terms of the Sobolev-type
stochastic derivative introduced in the previous subsection.

\begin{thm}\label{t212}
Consider the $\Re^m$-valued  random vector
$f=(f_1,\dots,f_m)$ which belongs for some grid $\Gf$ to
$W_2^1(\Gf,\Re^m)$. Denote, by
$\Sigma^{f,\Gf}=(\Sigma^{f,\Gf}_{k,r})_{k,r=1}^m$, the Malliavin
matrix for $f$,
$$\Sigma^{f,\Gf}_{k,r}\equiv (D^{\Gf} f_k,D^{\Gf}
f_r)_{\ell_2},\quad k,r=1,\dots,m,$$ and put $
\Nf(f,\Gf)=\{\omega|\Sigma^{f,\Gf}(\omega)$ is non-degenerate$\}$.
Then
$$
P\Bigl|_{\Nf(f,\Gf)}\circ f^{-1}\ll \lambda^m.
$$
\end{thm}

{\it The proof } is made in the framework of the stratification
method (see \cite{Dav_Lif_Smor}, Chapter 2 for the basic
constructions of this method) and contains several standard steps.
First, let us choose
  a countable set $\ell_*\subset \ell_0$ dense in
$\ell_2$. For any $\bar l=(l^1,\dots,l^m)\in [\ell_*]^m$, we denote
$$
\Nf(f,\bar l)=\{\omega| \hbox{ the matrix }\left((D^{\Gf}
f_k,l^r)_{\ell_2}\right)_{k,r=1}^m \hbox{ is non-degenerate}\}.
$$
Then $\Nf(f,\Gf)=\cup_{\bar l\in [\ell_*]^m }\Nf(f,\bar l)$ and
thus, in order to prove the statement of the theorem, it is
sufficient to prove that, for every fixed $\bar l\in [\ell_0]^m$,
\be\label{22}
P\Bigl|_{\Nf(f,\bar l)}\circ f^{-1}\ll \lambda^m.
\ee
The set $\bar l$ generates the commutative  group of admissible
transformations of $(\Omega,\Ff,P)$, indexed by $\Re^m$:
$$
T_t\equiv T_{t_1l^1}^{\Gf}\circ \dots\circ T_{t_ml^m}^{\Gf},\quad
t=(t_1,\dots,t_m).
$$
In order to prove  (\ref{22}), we proceed in the following way. Consider
the stratification of $(\Omega,\Ff,P)$ on the orbits of the group
$\{T_t,t\in\Re^m\}$, which can be considered in our case after a proper
parametrization as $\Re^m$ or some proper linear
subspaces of $\Re^m$. The group $\{T_t\}$ generates a measurable
parametrization of $(\Omega,\Ff,P)$ (the detailed exposition will
be given further), and thus $P$ can be decomposed into a regular
family of conditional distributions such that every conditional
distribution is supported by some orbit. Denote, by $\rho_{\bar
l}\equiv
(\sum_{i\in\NN}l^1_i\rho_{h_i}^{\Gamma^i},\dots,\sum_{i\in\NN}l^m_i\rho_{h_i}^{\Gamma^i})$,
the logarithmic derivative of $P$ w.r.t. $\{T_t\}$. Then, for
almost all orbits $\gamma$, the conditional distribution
$P_\gamma$, supported by the orbit $\gamma$, possess the
logarithmic derivative $\rho_{\bar l,\gamma}$, that is equal to
the restriction of $\rho_{\bar l}$ on the orbit $\gamma$. Since
$\rho_{\bar l}$ has an exponential moment, $\rho_{\bar l,\gamma}$ has such a moment too for almost all
$\gamma$. This
implies (see \cite{Boga}, Proposition 4.3.1) that, for almost all
$\gamma$, $P_\gamma$ possesses a positive continuous density.

On the almost every orbit  $\gamma$, the function $f_\gamma$ is
equal to the restriction of $f$ on $\gamma$ and  belongs to the
Sobolev class $\cap_p W^1_p(P_\gamma)$. This fact is more or less
standard and we do not give the proof here. In a linear framework,
this subject was discussed in details in \cite{Me_Mark_Uniq}. The
non-linear case of a commutative admissible group $\{T_t\}$ is quite
analogous. We refer the interested reader to \cite{Me_Mark_Uniq}
and references therein.

Taking into account this analytic background, we can apply the change-of-variables formula on
the almost every orbit $\gamma$
and obtain
 the absolute continuity of the image of the measure
$P_\gamma$ under the map $f_\gamma$. After all, (\ref{22}) is obtained
by the Fubini theorem. We omit this part of the exposition,
referring the reader to \cite{Dav_Lif_Smor}, Chapter 2, or
\cite{Pilip_1996}.

Now let us verify that our specific group $\{T_t\}$ generates a
measurable parametrization of $(\Omega,\Ff,P)$, i.e. there exists
a measurable map $\Phi:\Omega\to \Re^m\times \tilde \Omega$ such
that $\tilde \Omega$ is a Borel measurable space and the image of
every orbit of the group $\{T_t\}$ under $\Phi$ has the form
$L\times\{\tilde\omega\}$, where $L$ is a linear subspace of
$\Re^m$. This condition was supposed to hold true under the
considerations made before.

In order to shorten the notation, we restrict ourselves to the case
where
$$
l_{i}^r=\begin{cases}1,&i=r\\
0,&\hbox{otherwise}\end{cases},\quad r=1,\dots,m,
$$
the general case is quite analogous. For $i=1,\dots, m$, we denote
$\Df_i=\{\tau\in \Df\cap (a_i,b_i)| p(\tau)\in\Gamma_{i}\}$,
$c_i={b_i-a_i\over 2}$. Let $\omega\in \Omega$ be fixed. We recall
that $\omega$ is interpreted as a (locally finite) configuration.
Set $I(\omega)=\{i|\Df_i(\omega)\not=\emptyset\}$ and, for $i\in
I(\omega)$, we define $\tau_i(\omega)=\inf \Df_i(\omega).$ Note that,
due to condition (ii) of Definition \ref{d26} for every
$i=1,\dots,m$ and $x\in(a_i,b_i)$, the transformation
$$
\Re\ni z\mapsto T_{zh_i}x
$$
is strictly monotonous and its image is equal to $(a_i,b_i)$.
Therefore, for every $i\in I(\omega)$ there exists the unique
$z_i(\omega)\in\Re$  such that $T_{z_i}^{i}\tau_i=c_i$. Denote
$z(\omega)=(z_1(\omega),\dots,z_m(\omega))\in\Re^m$, where
 $z_i(\omega)=0$ for $i\not\in I(\omega)$. Denote by $\tilde
 \Omega$ the set of all configurations satisfying the following
 additional condition: for every cell $\Gf_{i}, i=1,\dots,m$, either the
 configuration is  empty in this cell, or the moment of the first
 jump in this cell is equal to $c_i$. Now put, for every $\omega\in\Omega\,$,
 $\varpi(\omega)= [T_{z(\omega)}\omega]\in\tilde \Omega$. Then
 the map
$$
\Phi:\omega\mapsto (z(\omega), \varpi(\omega))
$$
provides the needed parametrization. The theorem is proved.

Another version of the previous result can be  given  in the terms
of the \emph{almost sure} stochastic derivative. Although we will
not use the framework of {almost sure} stochastic derivatives
while studying equation (\ref{01}), it can be very useful  while
studying the distributions of some other classes of functionals.
Thus we formulate briefly the main points of this framework.

\begin{dfn}\label{d213} For a given grid $\Gf$, the functional
$f$ is called to be almost surely (a.s.) differentiable w.r.t.
$\Gf$,  if there exists a random element $\tilde D^{\Gf}f$ with
values in $\ell_2$ such that, for every
 $l\in \ell_0$, $$
{1\over t}\left[f\circ T^\Gf_{l}-f\right]\to (\tilde D^{\Gf}
f,l)_{\ell_2}, \quad t\to 0\quad \hbox{almost surely.}
$$
 The element $D^{\Gf} f$ is called the almost
sure (a.s.) derivative of $f$ w.r.t. $\Gf$.
\end{dfn}

\begin{thm}\label{t214} Consider the random vector
$f=(f_1,\dots,f_m)$ such that, for some grid $\Gf$, every functional
$f_r, r=1,\dots,m$  is a.s. differentiable w.r.t. $\Gf$. Denote
$\tilde \Sigma^{f,\Gf}=(\tilde \Sigma^{f,\Gf}_{k,r})_{k,r=1}^m$,
$$\tilde \Sigma^{f,\Gf}_{k,r}\equiv (\tilde D^{\Gf} f_k,\tilde D^{\Gf}
f_r)_{\ell_2},\quad k,r=1,\dots,m,$$ and put $ \tilde
\Nf(f,\Gf)=\{\omega|\tilde \Sigma^{f,\Gf}(\omega)$ is
non-degenerate$\}$. Then
$$
P\Bigl|_{\tilde \Nf(f,\Gf)}\circ f^{-1}\ll \lambda^m.
$$
\end{thm}
\demo Due to the arguments given in the proof of the previous
theorem, it is sufficient to prove the same statement in a
finite-dimensional case, i.e. when $\Omega$ is $\Re^m$ and
$\Tf^{\Gf}$ is the canonical  group of linear shifts in $\Re^m$.
In this situation the needed statement holds true due to the
standard change-of-variables formula and the following lemma.

\begin{lem}\label{l215} Let, for some $m,n\in \NN$, the function
$F:\Re^m\to\Re^n, G:\Re^m\to \Re^m\times \Re^n$ be such that, for
every $a\in\Re^m$ for $\lambda^m$-almost all $x\in \Re^m$,
$$
{1\over t}\|F(x+ta)-F(x)-t(G(x),a)_{\Re^m}\|_{\Re^n}\to 0,\quad
t\to 0.
$$
Then, for every $\eps>0$, there exists $F_\eps\in C^1(\Re^m,\Re^n)$
such that
$$ \lambda^m(\{x|
F(x)\not=F_\eps(x)\}\cup\{x| G(x)\not=\nabla F_\eps(x)\})<\eps.
$$
\end{lem}

This result is a straightforward consequence of the Lebesgue
theorem about the points  of density for a measurable set and the
following statements.

\begin{prop}\label{p216} I. (\cite{Federer}, Theorem 3.1.4). Let
the function $f:\Re^m\to \Re^n$ be approximatively
differentiable at every point of a set $A\subset \Re^m$  along all
the vectors from the basis. Then, for $\lambda^m$-almost all points
$a\in A$, the function $f$ has the approximative derivative at $a$.

II. (\cite{Federer}, Theorem 3.1.16). Let $A\subset \Re^m, f:A\to
\Re^n$ and
\be\label{23}
\mathop{ap\lim\sup}\limits_{x\to
a}{\|f(x)-f(a)\|_{\Re^n}\over \|x-a\|_{\Re^m}}<+\infty
\ee
 for $\lambda^m$-almost all $a\in A$. Then, for every
$\eps>0$, there exists $g\in C^1(\Re^m, \Re^n)$ such that
$$ \lambda^m(\{x|f(x)\not=g(x)\}<\eps.
$$
\end{prop}
We are not going to discuss definitions of the approximative limit
and derivative here, referring the reader to \cite{Federer}. Let
us only mention that the usual differentiability along some
direction implies the approximative differentiability along this
direction, and if the approximative derivative exists, then
(\ref{23}) holds true. Theorem \ref{t214} is proved.

The following theorem gives the convergence in variation of the
distribution of random vectors in terms of their derivatives, and
will be used in the proof of Theorem \ref{t12}.

\begin{thm}\label{t217}
For some given  grid $\Gf$ and $p>m$,
consider the sequence of $\Re^m$-valued random vectors
$\{f^n\}\subset W_p^1(\Gf,\Re^m)$ such that
$$
f_n\to f\hbox{ in }W_{p,1}^{\Gf, \Re^m},\quad n\to +\infty.
$$
 Then, for every $A\subset{\Nf(f,\Gf)}$,
$$
P\Bigl|_{A}\circ f^{-1}_n\to P\Bigl|_{A}\circ f^{-1},\quad n\to
+\infty
$$
in variation.
\end{thm}

The statement of the theorem follows, via
  the stratification
arguments analogous to those given in the proof of Theorem \ref{t212},
  from the
finite-dimensional criterion for the convergence in variation of
the sequence of induced measures, given in
\cite{Pilipenko_convergence_by_variation} (see
\cite{Pilipenko_convergence_by_variation}, Theorem 2.1 and
Corollary 2.7).

Let us mention that the analog of Theorem \ref{t217} can be also
given in the terms of the almost sure derivatives, but an
additional uniform condition on the sequence $\{f_n\}$ should be
imposed in this case. We do not discuss this subject here,
referring the interested reader to \cite{Me_conv_var}.

\section{Absolute continuity of the distribution of
a solution to an SDE with jumps}

\subsection{Differential properties of the solution to an SDE
with jumps}

 We are going to apply the general results about
the existence of the density obtained in the previous section to the
specific class of functionals: solutions to SDE's with jumps.
The first step, that is necessary here, is to verify whether such
solutions are either stochastically or a.s. differentiable. In
this subsection, we give the answer to  this question.

Consider the Cauchy problem for equation (\ref{01}) of the type
\be\label{31} X(x,t)=x+\int_0^t a(X(x,s))\,ds+U_t-U_0,\quad
t\in\ax. \ee We suppose that $a$ belongs to $C^1(\Re^m,\Re^m)$. We
also impose the linear growth condition on $a$:
$$
\exists K: \|a(x)\|^2\leq K(1+\|x\|^2).
$$
These conditions  provide that equation (\ref{31})
 has the unique strong solution. Moreover, these solutions
 considered for different $x,t$ form a \emph{stochastic flow of
 diffeomorphisms}.

Denote $\Delta(x,u)=a(x+u)-a(x), x\in\Re^m, u\in\Re^d.$

\begin{thm}\label{t31} \textbf{I.} For every
$x\in\Re^m,t\in\ax,\Gamma\in \pf,h\in H_0$, every component of the
vector $X(x,t)$ is a.s. differentiable w.r.t.
$\{T_{rh}^\Gamma,r\in\Re\}$, i.e. there exist a.s. limits
$$
Y_k(x,t)=\lim_{\eps\to 0} {1\over \eps} [T_{\eps h}^\Gamma
X_k(x,t)-X_k(x,t)],\quad k=1,\dots,m.
$$
The process $Y(x,\cdot)$ satisfies the equation
\be\label{32}
 Y(x,t)=\int_0^t\int_\Gamma \Delta(X(x,s-),u
 )Jh(s)\,\nu(ds,du)+
\int_0^t [\nabla a](X(x,s)) Y(x,s)\,ds ,\quad t\geq 0.
\ee

\textbf{II.} The solution $X(x,t)$ is stochastically
differentiable with the derivative given by (\ref{32}).
\end{thm}
\emph{Remark.} In a sequel, we use only statement II. Statement I
provides here the main part of the proof and is emphasized only
for the convenience of the reader.

{\it Remark.} The statement close to statement I was proved in
\cite{Me_TViMc}.  The statement close to statement II
 was proved in \cite{Nou_sim} for $m=1$. We cannot use
 straightforwardly the result from  \cite{Nou_sim} since the proof
 there contains some specifically one-dimensional features such
 as an exponential formula for the derivative of the flow
 corresponding
 to the solution of the ODE (Lemma 1 \cite{Nou_sim}).

\emph{Proof of statement }I. It is sufficient to consider only
the case where $a,\nabla a$ are bounded. The general case follows
from this one due to the standard localization arguments.

Denote $\Df^\Gamma=\{\tau\in\Df:p(\tau)\in\Gamma\}$,
$\Omega_k=\{\Df\cap\{0,t\}=\emptyset,\#(\Df^\Gamma\cap (0,t))=k,
\},k\geq 0.$ Since $\Gamma\in\pf$, $\Omega=\cup_k\Omega_k$ almost
surely and it is enough to verify that the needed statement  holds
true a.s. on every $\Omega_k$. The case $k=0$ is trivial.

Denote $\nu^*(t,A)=\nu(t,A\backslash
\Gamma),U^*_t=\int_0^t\int_{\Re^d\backslash \Gamma}u\tilde
\nu(ds,du)$. For a given $t>0, \tau\in(0,t),p\in\Re^d,x\in\Re^m$,
consider the process $X_\cdot^{\tau}$ on $[0,t]$ such that
$$
X_{t}^{\tau}=\begin{cases}x+\int_{0}^t
a(X^{\tau}_s)\,ds+U_t^*,&t<\tau\\
x+\int_{0}^t b(X^{\tau}_s)\,ds+p+ U_t^*,&t\geq \tau\end{cases}.
$$
 Note that the point process $\{p(T),
T\in\Df^\Gamma\}$ is independent of $\nu^*$, and the distribution
of the variable $\tau_1^\Gamma\equiv \min \Df^\Gamma$, while this
variable is restricted to $\Omega_1,$  is absolutely continuous.
Then statement I on $\Omega_1$ follows immediately from
Example \ref{e22} and the following lemma.

\begin{lem}\label{l32} With probability 1 for
$\lambda^1$-almost all $\tau\in (0,t)$,
$$
{d\over d\eps}\Bigl|_{\eps=0}X_{t}^{\tau+\eps}=-
\Delta(X_{\tau-}^{\tau},p)\Ef^{*}_{t},
$$
where $\Ef^{*}$ is the stochastic exponent defined by the
equation
$$
\Ef^{*}_{r}=I_{\Re^m}+\int_\tau^r \nabla
b(X^{\tau}(s))\Ef^{*}_{s}\,ds,\quad r\geq \tau.
$$
\end{lem}

\demo $X_{t}^{\tau}$ is the value at the point $t$ of the
solution to the equation
\be\label{33}
d\tilde X_t=a(\tilde X_t)\,dt+dU_t^*,
\ee
with the starting point $\tau$ and the initial value
$$
X_{\tau}^\tau=X_{\tau-}^{\tau}+p.
$$
 Suppose that
 $\eps<0$.  Then $X_{s}^\tau=X_{s}^{\tau+\eps}, s<\tau+\eps.$ Thus
  $X_t^{\tau+\eps}$ is also the value of the solution to the same
 equation with the same starting and terminal points and with  the initial value being equal to
$$
X_{(\tau+\eps)-}^{\tau}+p+\int_{\tau+\eps}^\tau
a(X_{s}^{\tau+\eps})\,ds+ [U^*_\tau-U^*_{\tau+\eps}].
$$
Thus the difference $\Phi(\tau,\eps)$ between the initial values
for $X_{t}^{\tau+\eps},X_{t}^{\tau}$ is equal to
 $\int_{\tau+\eps}^\tau[a(X_s^{\tau+\eps})-
a(X_s^{\tau})]\,ds$.

The process  $\{U_t^*\}$ has  c\`{a}dl\`{a}g trajectories, and
therefore almost surely the set of discontinuities for its
trajectories is at most countable.  Therefore almost surely there
exists the set $\TT=\TT(\omega)\subset\ax$ of the full Lebesgue
measure such that
$$
\delta(t,\gamma)\equiv \sup_{|s-t|\leq\gamma}[\|U_s^*-U_t^*\|]\to
0,\quad \gamma\to 0,\quad t\in\TT.
$$
Then, for $s\in(\tau+\eps,\eps)$,
$$
\|X_{s}^\tau-X_{\tau-}^\tau\|+\|X_{s}^{\tau+\eps}-X_{\tau-}^\tau
-p\|\leq \Cd\{|\eps|+\delta(\tau,|\eps|)\}.
$$
Here and below, we denote, by $\Cd$, any constant such that it can be
calculated explicitly, but its  exact form is not needed in a
further exposition. Thus, for $\tau\in\TT$,
$$
\|\Phi(\tau,\eps)+\eps[a(X_{\tau-}^\tau+p)-a(X_{\tau-}^\tau)]\|\leq
\Cd|\eps| \{|\eps|+\delta(\tau,|\eps|)\},
$$
which implies the needed statement.

The case $\eps>0$ is analogous, let us discuss it briefly. Again,
take $\tau\in\TT$ and represent $X_{t}^{\tau}$ as the solution
to (\ref{33})  with the initial value $X_{\tau-}^\tau+p$.
$X_{t}^{\tau+\eps}$ is also the solution to (\ref{33}) but with the
other starting point $\tau+\eps$. The estimates analogous to ones
made before show that, up to the $o(|\eps|)$ terms,
$$
X_{\tau+\eps}^{\tau+\eps}-X_{\tau+\eps}^\tau=\eps \Bigl\{
-a(X_{\tau-}^\tau+p) +a(X_{\tau-}^\tau)\Bigr\},
$$
which implies the statement of the lemma. The lemma is proved.

Now let $k>1$ be fixed. Consider the countable family $\Qf_k$ of
the partitions $Q=\{0=q_0<q_1\dots<q_k=t\}$ with
$q_1,\dots,q_{k-1}\in\QQ$ and denote
$$
\Omega_Q=\{\Df\cap\{q_i, i=0,k\}=\emptyset,
\Df^\Gamma\cap(q_{i-1},q_i)=1,i=1,\dots,k\},\quad Q\in\Qf_k.
$$
We have $\Omega_k=\cup_{Q\in\Qf_k}\Omega_Q$. Therefore it is
enough to verify the statement of Theorem \ref{t31} on $\Omega_Q$
for a given $Q$. The distributions of the variables
$\tau_j^\Gamma, j=1,\dots,k$ (see Example \ref{e22} for the
notation $\tau_j^\Gamma$), while these variables are restricted to
$\Omega_k,$ are absolutely continuous. Then statement I on
$\Omega_Q$ follows immediately from Example \ref{e22}, the
standard theorem  about differentiation of the solution to
equation (\ref{31})  w.r.t. the initial value, and the statements
analogous of one of Lemma \ref{l32} and written on the intervals
$[0,q_1],[q_1,q_2],\dots,[q_{k-1},t]$. Statement I is proved.

\emph{Proof of statement }II. Again, suppose first that
$a,\nabla a$ are bounded. In the framework of Lemma \ref{l32}, one has
the estimate
\be\label{34}
\|X_t^{\tau+\eps}-X_t^{\tau}\|\leq \Cd |\eps|
\ee
valid point-wise. Indeed, both $X_t^{\tau+\eps}$ and $X_t^{\tau}$
are the solutions to (\ref{33}) with the same initial point ($\tau$ for
$\eps<0$ and $\tau+\eps$ for $\eps>0$) and different initial
values. The difference between the initial values are estimated by
$$
\left\|\int_{\tau+\eps}^\tau[a(X_s^{\tau+\eps})-
a(X_s^{\tau})]\,ds\right\|\leq -2\|a\|_\infty \eps
$$
for $\eps<0$ and by
$$
\left\|\int^{\tau+\eps}_\tau[a(X_s^{\tau+\eps})-
a(X_s^{\tau})]\,ds\right\|\leq 2\|a\|_\infty \eps
$$
for $\eps>0$. Thus, inequality (\ref{34}) follows from  the Gronwall lemma. Using
the described before technique, involving partitions $Q\in\Qf_k$,
and applying the Gronwall lemma once again, we obtain that almost
surely on the set $\Omega_k$
$$
\|T_{\eps h}^\Gamma X(x,t)-X(x,t)\|\leq k\Cd |\eps|.
$$
This means that the family $\{{1\over \eps}[T_{\eps h}^\Gamma
X(x,t)-X(x,t)]\}$ we already have proved to converge to the
solution to (\ref{32}) almost surely as $\eps\to 0$ is dominated
by the variable
$$
\Cd \cdot \nu(t,\Gamma)\in\cap_p L_p(\Omega,\Ff,P).
$$
Therefore the convergence holds true also in the $L_p$ sense for any
$p$, and $X(x,t)$ is stochastically differentiable with the
derivative given by (\ref{32}).

The last thing we need to do is to remove the claim on $a$
to be bounded. Consider a sequence $\{a_n\}\subset
C^1_b(\Re^m,\Re^m)$ such that $a_n(x)=a(x)$ for $\|x\|\leq n$. We
have just proved that the solution $X_n(x,t)$ to an equation of
the type (\ref{31}) with $a$ replaced by $a_n$ is stochastically
differentiable and its derivative $Y_n(x,t)$ is given by an
equation of the type (\ref{32}) with $a$ replaced by $a_n$. The
sequence $\{a_n\}$ can be chosen in such a way  that it satisfies
the linear growth condition uniformly w.r.t. $n$. Under such a
choice,
$$
X_n(x,t)\to X(x,t),\quad Y_n(x,t)\to Y(x,t),\quad n\to +\infty
$$
in every $L_p(\Omega, P,\Re^m)$. Since the stochastic derivative is a
closed operator, this implies the needed statement for $X(x,t)$.
The theorem is proved.

\subsection{The proofs of Theorems \ref{t11}, \ref{t12}.}

 The proof
of Theorem \ref{t11} is an essentially simplified version of the
proof of the analogous statement in \cite{Me_TViMc}. It is based
on the other version of the absolute continuity result, with the
conditions formulated in the terms of the point process
$\{p(\tau),\tau\in\Df\}$.  Below the initial value $x_*$ is fixed,
and we omit it in the notation writing $X(s)\equiv X(x_*,s)$.

Denote, by $\{\Ef_r\}$, the stochastic exponent,  i.e. the
 $m\times m$-matrix-valued process satisfying
the equation
$$
\Ef_r=I_{\Re^m}+\int_0^r \nabla a(X(s))\Ef_{s}\,ds,\quad r\in\ax.
$$
This process has continuous trajectories. The matrix $\Ef_r$ is
a.s. invertible for every $r$, and, moreover, almost surely
$$
\sup_{r\leq t} \|\Ef_r^{-1}\|_{\Re^m\times\Re^m}<+\infty.
$$
We do not discuss this fact in details, since the technique is
quite standard here (see, for instance \cite{protter}, Chapter 5,
\S 10).

 \begin{lem}\label{l33} Denote by $S_t$
a linear span of the set of vectors $\{\Ef_{\tau}^{-1}\cdot
\Delta(X(\tau-),p(\tau)),\tau\in\Df\cap (0,t)\}$ and put
$\Omega_t=\{\omega|\,\hbox{{\rm dim}}\,S_t(\omega)=m\}.$ Then
$$
P|_{\Omega_t}\circ [X(t)]^{-1}\ll\lambda^m.
$$
\end{lem}
\demo Denote, by $S^n_t$, a linear span of the set of vectors
$\{\Ef_{\tau}^{-1}\cdot \Delta(X(\tau-),p(\tau)),\tau\in\Df\cap
(0,t),\|p(\tau)\|_{\Re^d}\geq {1\over n}\}$ and put
$\Omega_t^n=\{\omega|\,\hbox{{\rm dim}}\,S_t^n(\omega)=m\}.$ It is
clear that $\Omega_t=\bigcup_{n\geq 1}\Omega_t^n$, and thus it is
enough to prove that  $P|_{\Omega^n_t}\circ
[X(t)]^{-1}\ll\lambda^m$ for a given $n$.

Let $n$ be fixed. Consider the family of differential grids
$\{\Gf^N,N\in\NN\}$ of the form
$$
\Gamma_i^N=\Gamma^n\equiv \{u| \|u\|\geq {1\over n}\}, \quad
a_i^N=b_{i-1}^N={i-1\over N}, \quad h_i^N(s)=h({s-a_i^N\over
b_i^N-a_i^N}),  \quad s\in(a_i,b_i), i\in\NN,
$$
where $h\in H_0$ is some function such that $Jh>0$ inside $(0,1)$
and $Jh=0$ outside $(0,1)$.

Our aim is to show that almost surely \be\label{35}
\Omega_t^n\subset \bigcup_{N}\,\{\omega|\tilde
\Sigma^{X(t),\Gf^N}(\omega) \hbox{ is non-degenerate }\}. \ee Here
$\Sigma^{X(t),\Gf^N}$ is the Malliavin matrix for the random
vector $X(t)$ (see Theorem \ref{t212}). Theorem \ref{t212}
together with (\ref{35}) immediately imply the needed statement.

Denote $\Df^n\equiv\Df^{\Gamma^n}$, $A_N^{n,t}=\Bigl\{\omega| \Df
\cap\{{i-1\over N},i\in \NN\}=\emptyset,  \#\{\tau\in \Df^n\cap
(a_i,b_i)\}\subset\{0,1\}, i=1,\dots, [Nt+1]\Bigr\}$. Since
$\Gamma^n\in\pf$, one has that  almost surely
$$
\Omega_t^n\subset\bigcup_N\, [\Omega_t^n\cap A_N^{n,t}].
$$
Thus in order to prove (\ref{35}), it is sufficient to show that,
for every $N$, the matrix $\Sigma^{X(t),\Gf^N}$ is non-degenerate
on the set $\Omega_t^n\cap A_N^{n,t}$.

A change of the point measure outside $[0,t]$ does not change
$X(t)$, thus
$$
(D^{\Gf^N} X(t),l)_{\ell_2(\Gf^N)}=0\hbox{ for any }l: \quad
l_{i}=0, \quad i\leq [Nt+1].
$$
This means that the matrix $\Sigma^{X(t),\Gf^N}$ is  the Grammian
for the finite family of the vectors in $\Re^m$
$$
Y^1\equiv(D^{\Gf^N}
X(t),l^1)_{\ell_2(\Gf^N)},\dots,Y^{[Nt+1]}\equiv(D^{\Gf^N}
X(t),l^{[Nt+1]})_{\ell_2(\Gf^N)},\quad l^r=(l^r_{i}),
\quad l_{i}^r=\begin{cases} 1,& i=r,\\
0,& \hbox{otherwise}.\end{cases}
$$
Therefore $\Sigma^{X(t),\Gf^N}$ is non-degenerate iff the family
$\{Y^r, r=1,\dots, [Nt+1]\}$ is of the maximal rank.

The family  $\{Y^r\}$ on the set $A_N^{n,t}$ can be given
explicitly. First of all, let us write the solution to equation
(\ref{32}) in the following form: \be\label{36}
Y_h^\Gamma(t)=\Ef_t\int_0^t\int_{\Gamma} Jh(s)\Ef_{s}^{-1}
\Delta(X(s-),u)\,\nu (ds,du),\quad t\geq 0, \quad h\in
H_0,\Gamma\in\pf. \ee Taking in (\ref{36}) $\Gamma=\Gamma^n$ and
$h=h_r^N, r=1,\dots,[Nt+1],$ we obtain that, on the set
$A_N^{n,t}$ $Y^r=c_r\Ef_t\tilde Y^r,$ $r=1,\dots,[Nt+1]$, where
$$
(c_r,\tilde Y^r)=\begin{cases} (
h_r^N(\tau_r),\Ef_{\tau_r}^{-1}\cdot
\Delta(X(\tau_r-),p(\tau_r))),& \{\tau\in
\Df^n\cap (a_r,b_r)\}=\{\tau_r\}, \\
(1,0),& \{\tau\in \Df^n\cap (a_r,b_r)\}=\emptyset.
\end{cases}
$$
The matrix $\Ef_t$ is non-degenerate, the constants $\{c_r\}$ are
positive on $A_N^{n,t}$. This means that $\{Y^r\}$ has the maximum
rank iff the same holds true for $\{\tilde Y^r\}$. But the family
$\{\tilde Y^r\}$ contains all the vectors
$$\{\Ef_{\tau}^{-1}\cdot
\Delta(X(\tau-),p(\tau)),\tau\in\Df^n\cap (0,t)\}
$$
and therefore has the maximal rank on $\Omega^n_t$. This means
that $\{Y^r\}$ has the maximal rank on $\Omega_t^n\cap A_N^{n,t}$
and (\ref{35}), together with the statement of the lemma, holds true.
The lemma is proved.

\begin{lem}\label{l34} Under the condition of Theorem \ref{t11},
\be\label{37} \gamma_n\equiv\inf\limits_{x\in \bar B(x_*,\eps_*),
v\not=0}\Pi\Bigl(u:\|u\|\geq {1\over
n},(\Delta(x,u),v)_{\Re^d}\not=0\Bigr)\to +\infty,\quad
n\to+\infty.\ee
\end{lem}

This statement follows immediately from the Dini theorem applied
to the monotone sequence of lower semi-continuous functions
$$
\phi_n: \bar B(x_*,\eps_*)\times \{\|v\|=1\}\ni(x,v)\mapsto
\Pi\Bigl(u:\|u\|\geq {1\over
n},(\Delta(x,u),v)_{\Re^d}\not=0\Bigr).
$$

 \emph{Proof of the Theorem \ref{t11}}  Denote by $\Sf$ the set of all proper
subspaces of $\Re^m$. This set can be parametrized in such a way
that it becomes a Polish space, and, for every of the random vectors
$\xi_1,\dots,\xi_k$, the map $\omega\mapsto
\Span(\xi_1(\omega),\dots,\xi_k(\omega))$ defines the random
element in $\Sf$.

For every $n\geq 1$, consider the set
$\Df^n=\{\tau_1^n,\tau_2^n,\dots\}$. For a given $S^*\in \Sf$,
$\delta>0$,  let us consider the event
$$
A_\delta^n=\{S_\delta^n\not \subset S^*\}=\{\exists i:
\tau_i^n\leq \delta,
\Ef^{-1}_{\tau_i^n}\Delta(X(\tau_i^n-),p(\tau_i^n))\not\in S^*\}
$$
(see the beginning of the proof of Lemma \ref{l33} for the
notation $S_t^n$). One has that $\Omega\backslash
A_\delta^n\subset B_\delta\cup C_\delta^n,$ where $
B_\delta=\{\exists s\in[0,\delta]: X(s-)\not\in \bar
B(x_*,\eps_*)\},$
 $$C_\delta^n=\bigcap_i\Bigl[\{\tau_i^n>\delta\}\cup \{X(\tau_i^n-)\in \bar
B(x_*,\eps_*),
\Ef^{-1}_{\tau_i^n}\Delta(X(\tau_i^n-),p(\tau_i^n))\in
S^*,\tau_i^n\leq \delta\}\Bigr].
$$

The distribution of the value $p(\tau_i^n)$ is equal to
$\lambda_n^{-1} \Pi|_{\Gamma^n}$, where $\Gamma^n=\{u| \|u\|\geq
{1\over n}\}, \lambda^n=\Pi(\Gamma^n)$. Moreover, this value is
independent with the $\sigma$-algebra $\Ff_{\tau_i^n-}$, and, in
particular, with the variables $X(\tau_i^n-), \Ef_{\tau_i^n}$.
This provides the estimate
\be\label{38}
P[\{\tau_i^n>\delta\}\cup \{X(\tau_i^n-)\in \bar B(x_*,\eps_*),
\Ef^{-1}_{\tau_i^n}\Delta(X(\tau_i^n-),p(\tau_i^n))\in
S^*,\tau_i^n\leq \delta\}\Bigr|\Ff_{\tau_i^n-}\Bigr]\leq
\1_{\{\tau_i^n>\delta\}}+(1-{\gamma_n\over
\lambda_n})\1_{\{\tau_i^n\leq\delta\}}.
\ee
It follows from (\ref{38}) that
$$
P(C_\delta^n)\leq E\left(1-{\gamma_n\over
\lambda_n}\right)^{\nu([0,\delta]\times\Gamma^n)}=\exp\{-\delta\gamma_n\}\to
0,n\to+\infty.
$$
Since $A_\delta^n\subset\{S_\delta\not \subset S^*\}$, this means
that almost surely
\be\label{39}
\{S_\delta \subset S^*\}\subset B_\delta.
\ee
Now we take $\delta<{t\over m}$ and  iterate (\ref{39}) on the time
intervals $[0,\delta],[\delta,2\delta],\dots,
[(m-1)\delta,m\delta]$ with $S^*_1=\{0\}, S^*_2=S_\delta,\dots,
S_m^*=S_{(m-1)\delta}$ (we can do this due to the Markov property
of $X$). We obtain that
$$
\{\dim S_t<m\}\subset \bigcup_{k=1}^m\{\dim S_{(k-1)\delta}=\dim
S_{k\delta}<m\}\subset B_{m\delta}.
$$
Since $P(B_{m\delta})\to 0, \delta\to 0+$, this provides  that
$P\{\dim S_t<m\}=0$, which together with Lemma \ref{l33} gives the
needed statement. The theorem is proved.

\emph{Proof of Theorem \ref{t12}.} Due to statement II of
Theorem \ref{t31}, the solutions $X_n(x_n,t_n)$ to (\ref{12})
are stochastically differentiable and their derivatives are given
by SDEs of the form (\ref{32}). The usual localization
arguments allows us to restrict the consideration to the case where
$\{a_n\}$ are uniformly bounded together with their derivatives
and $\Pi$ is supported by some bounded set. Then, applying Theorem
4, \cite{Gikh_Skor}, Chapter 4.2,  we obtain that, for any $p>1$,
$X_n(x_n,t_n)$ converge to $X(x_*,t)$ in the $L_p$ sense, together
with their stochastic derivatives given by (\ref{32}). This means
that, for every {finite} differential grid $\Gf$ and any $p>1$,
$$
X_n(x_n,t_n)\to X(x_*,t)\quad\hbox{in }W_p^1(\Gf,\Re^m),\quad n\to
+\infty.
$$
Thus the statement of Theorem \ref{t12} follows from Theorem
\ref{t217}.

\subsection{The proofs of Propositions \ref{p36} -- \ref{p37}}

{\it Proof of the Proposition \ref{p36}.} Take
$\eps_*={\delta_*\over 2}$. Then, for every $x\in \bar
B(x_*,\eps_*)$,
$$\{u|\Delta(x,u)=0\}=\{u|a(x+u)=a(x)\}\subset \{|u|>\Cd
\delta_*\}\cup\{u|x+u\in N(a,a(x))\cap
(x_*-\delta_*,x_*+\delta_*)\}=\Delta_1\cup\Delta_2.
$$
Here we used that $a$ is Lipschitz. The set $\Delta_2$ is finite
and therefore $\Pi(\Delta_2)<+\infty$. The set $\Delta_1$ is
separated from $0$ and therefore $\Pi(\Delta_1)<+\infty$. Since
$\Pi(\Re)=+\infty$, this means that
$\Pi(\Delta(x,u)\not=0)=+\infty$. Proposition is proved.

The proof of the Proposition \ref{p310} is almost trivial: for a
given $x\in B(x_*,\eps_*),v\in S_m$ one should take $w=w(x,v)$,
given by the Definition \ref{d15}, and for this $v$ choose
$\varrho\in(0,1)$ such that $\Pi(V(w,\varrho))=+\infty$ (this is
possible since $\Pi$ satisfies the wide cone condition). Then, for
every $D>0$, $\Pi(u\in V(w,\varrho), \|u\|\leq D)=+\infty$, and
(\ref{11}) follows from (\ref{13b}).

{\it Proof of the Proposition \ref{p37}.} Consider the set
$\Phi_{x_*,\eps_*}$ of the functions $\phi_x:\Re^m\ni u\mapsto
a(x+u)-a(x)\in \Re^m, x\in B(x_*,\eps_*)$. It is easy to see that
if for every linear subspace $L_v\equiv \{y|(y,v)=0\}, v\in S_m,$
$$
\Pi(u|u\not\in \phi^{-1}(L_v))=+\infty,\quad \phi\in
\Phi_{x_*,\eps_*},
$$
then (\ref{11}) holds true. In the case {\bf b},
$\Phi_{x_*,\eps_*}$ contains the unique function $\phi(u)=Au$.
Since $A$ is non-degenerate, $\phi^{-1} (L_v)$ is a proper linear
subspace of $\Re^m$ for every $v\in S_m$, and (\ref{311}) provides
(\ref{11}). In the case {\bf a}, $\phi_x\in C^1(\Re^m, \Re^m)$,
and for $\eps_*$ small enough $\det \nabla \phi_x(0)=\det \nabla
a(x)\not=0, x\in B(x_*,\eps_*)$. Then $\phi^{-1}_x (L_v)$ is a
proper smooth subspace of $\Re^m$ for every $v\in S_m$, and
(\ref{310}) provides (\ref{11}). Proposition is proved.

\section{Smoothness of the density of the solution
to the Cauchy problem}

\subsection{ The irregularity properties of the
density.}

We start our exposition with the easier part: the proof of Theorem
\ref{t18} and Proposition \ref{p111}, that give the irregularity
properties of $p_{x,t}$.

Recall that the function $a$ is supposed to be globally Lipschitz
and the jump noise is supposed to satisfy the moment condition
(\ref{13}).

{\it Proof of Theorem \ref{t18}: the case $m=1$.} For $\eps\in
(0,1)$, denote $M^\eps=\int_{\eps<|u|\leq 1} u\,\Pi(du)$ and
consider a decomposition of the process $U_t$ of the form
$$
U_t=U_0+R_t^\eps+V_t^\eps-tM^\eps,\quad
R_t^\eps=\int_0^t\int_{\|u\|\leq \eps} u \tilde \nu(ds,du),\quad
V_t^\eps=\int_0^t\int_{\|u\|>\eps} u \nu(ds,du).
$$
$R^\eps_t$ is a martingale, and its quadratic variation is equal
to
$$
[R^\eps]_t=\sum_{s\leq t}(R_s-R_{s-})^2=\int_0^t \int_{\|u\|\leq
\eps} u^2  \nu(ds,du).
$$
We have that
$$
E[R^\eps]_t^2=E\left[\int_0^t \int_{\|u\|\leq \eps} u^2 \tilde
\nu(ds,du)\right]^2+\left[\int_0^t \int_{\|u\|\leq \eps} u^2
\Pi(du)\, ds \right]^2=
$$
$$
=t \int_{\|u\|\leq \eps} u^4\Pi(du)+t^2\left[\int_{\|u\|\leq \eps}
u^2\Pi(du)\right]^2,
$$
and therefore, for $\eps$ small enough,
$$
E[R^\eps]_t^2\leq \Cd\left[\eps^{2} \ln{1\over\eps}\right]
\rho(\eps).
$$
 Applying the Chebyshev and
Burkholder inequalities, we obtain that, for every given $\alpha>0$,
$$
P(\sup_{s\leq t} |R_s^\eps|\geq \eps^{1-\alpha})\leq {\Cd E
[R^\eps]_t^2\over \eps^{2-2\alpha}}\leq \Cd\left[\eps^{2\alpha}
\ln{1\over\eps}\right] \rho(\eps).
$$
Next, for every $\eps\in(0,1)$
$$
P(V_s^\eps=0,s\in[0,t])=\exp[-t\Pi(|u|>\eps)]\geq \exp[-\eps^{-2}
t\int_{\Re} u^2\wedge \varepsilon^2 \Pi(du)]=\exp[t\rho(\eps)\ln
\eps]=\eps^{t\rho(\eps)}.
$$
Denote $A^\eps_\alpha=\{|R_s^{\eps_n}|\leq \eps_n^{1-\alpha},
V^\eps_s=0, s\in[0,t]\}$. Since $R^\eps$, $V^\eps$ are
independent, we have
$$
P(A^\eps_\alpha)\geq
\eps^{t\rho(\eps)}\left[1-\Cd\left[\eps^{2\alpha}
\ln{1\over\eps}\right] \rho(\eps)\right].
$$
Considering a sequence $\eps_n\to 0+$ such that $\rho(\eps_n)\to
\bro, n\to+\infty$, we obtain that, for $n$ big enough,
$$
P(A^{\eps_n}_\alpha)\geq {1\over 2}\eps_n^{t(\sbro+\alpha)}.
$$
Denote, by $X^n(x,t)$, the solution to the ODE \be\label{41}
X^n(x,t)=x+\int_0^ta(X^n(x,s))\,ds-tM^{\eps_n}. \ee By the
construction of the set $A_\alpha^{\eps_n}$, we have that on this
set
$$
|X(x, s)-X^n(x,s)|\leq L\int_0^s |X(x, r)-X^n(x,r)| dr
+\eps^{1-\alpha},\quad s\in [0,t],
$$
where $L$ denotes the Lipschitz constant for $a$. Then, by the
Gronwall lemma, $|X(x,t)-X^n(x,t)|\leq e^{Lt} \eps_n^{1-\alpha}$
on the set $A_\alpha^{\eps_n}$.  Thus there exist two sequences
$y_n=X^n(x,t)-e^{Lt}\eps_n^{1-\alpha},z_n=X^n(x,t)+e^{Lt}\eps_n^{1-\alpha}$
such that, for $n$ big enough, \be\label{42} P(y_n\leq X(x,t)\leq
z_n)\geq \Cd (z_n-y_n)^{(\bro+\alpha)\cdot {t\over 1-\alpha}}. \ee

Now we can complete the proof. For $y<z$ \be\label{43}\int_y^z
f(v)\, dv\leq \|f\|_{L_\infty}(z-y), \quad  \int_y^z f(v)\, dv\leq
\|f\|_{L_r}\left[ \int_y^z 1^{r\over r-1}\,dv\right]^{r-1\over
r}=\|f\|_{L_r} (z-y)^{r-1\over r}, r\in [1,+\infty). \ee Let
$t\bro<1-{1\over r}$. Then there exists $\alpha>0$ such that
$(\bro+\alpha)\cdot {t\over 1-\alpha}\ < 1-{1\over r}$ and
(\ref{42}) together with (\ref{43}) indicates that $p_{x,t}\not\in
L_r(\Re)$. This proves the statement {\bf a1}.
 Analogously, if $\bro t<1$, then there exists
$\alpha>0$ such that $(\bro+\alpha)\cdot {t\over 1-\alpha}\ < 1$
and (\ref{42}), (\ref {43}) indicate that $p_{x,t}$ is not
bounded, i.e. $p_{x,t}\not\in CB^0(\Re)$. This proves the
statement {\bf b1}. Under condition (\ref{13}) there exists
 $\lim_{n\to+\infty} M^{\eps_n}=M^0$ and  the sequences $\{y_n\},\{z_n\}$ are
bounded, that implies statements {\bf a},{\bf b}.

{\it Proof of Theorem \ref{t18}: the case $m>1$.} Consider
a decomposition of the process $U=(u^1,\dots,U^m)$ of the form
$$
U_t^i=U_0^i+R_t^{\eps,i}+V_t^{\eps,i}-tM^{\eps,i},\quad i=1,\dots,
m,
$$
where $R_t^{\eps,i}=\int_0^t\int_{|u^i|\leq \eps} u^i \tilde
\nu(ds,du), V_t^{\eps,i}=\int_0^t\int_{|u^i|>\eps} u^i
\nu(ds,du),M^{\eps,i}=\int_{\eps<|u^i|, \|u\|\leq 1}
u^i\,\Pi(du)$. Then, analogously to the proof of the case $m=1$, one
can verify that
$$
P(\sup_{s\leq t} \|R_s^\eps\|\geq \eps^{1-\alpha})\leq
\Cd\left[\eps^{2\alpha} \ln{1\over\eps}\right] \vartheta(\eps).
$$
On the other hand,
$$
P(V_s^\eps=0,s\in[0,t])=\exp\left[-t\Pi\left(\bigcup_{i=1}^m\{|u^i|>\eps\}\right)\right]\geq
\exp\left[-t\sum_{i=1}^m\Pi\left(\{|u^i|>\eps\}\right)\right]\geq
\eps^{mt\vartheta(\eps)}.
$$
Then, just as in the case $m=1$,  for every $\alpha\in(0,1)$ there
exist sequences $\{y_n^i\}\subset\Re^m, i=1,\dots, m$ and
$\{\delta_n\}\Re^+$ such that $\delta_n\to 0$ and \be\label{44}
P(X_n^i\in [y_n^i, y_n^i+\delta_n], i=1,\dots, m)\geq \Cd
\delta_n^{(\bkap+\alpha)\cdot {mt\over 1-\alpha}}. \ee The
arguments analogous to those used in the proof of the case $m=1$
show that (\ref{44}) implies statements {\bf a},{\bf b},{\bf
a1},{\bf b1} of Theorem \ref{t18}. The theorem is proved.

{\it Proof of Proposition \ref{p111}.} If $\bro_1=+\infty$, then
the statement is trivial. Thus we consider only the case
$\bro_1<+\infty$. Without losing generality, we can suppose that
$\Pi((-\infty,0))=0$.

Consider a sequence $\{\eps_n\}$ such that $\rho_1(\eps_n)\to
\bro_1$.  Since $\rho_1(\eps)\geq \rho_2(\eps)$, for the sequences
$\{y_n\},\{z_n\}$ given in the proof of Theorem \ref{t18} (the
case $m=1$), the following estimate holds true: \be\label{450}
P(y_n\leq X(x,t)\leq z_n)\geq \Cd (z_n-y_n)^{(\bro_1+\alpha)\cdot
{t\over 1-\alpha}}. \ee Denote, by $X^*(x,t)$, the solution to an
ODE of the type (\ref{41}) with $M^{\eps_n}$ replaced by $M^0$. It
follows from the comparison theorem that the law of $X(x,t)$ is
supported by $[X^*(x,t),-\infty)$ and the density $p_{x,t}$ is
equal to zero on $(-\infty,X^*(x,t)).$ On the other hand,
$M^0-M^{\eps_n}\leq \Bigl[\varepsilon_n \ln{1\over
\varepsilon_n}\Bigr]\rho_1(\eps_n)=o(\eps_n^{1-\alpha}),n\to+\infty$
and therefore, for $n$ big enough,
$(y_n,z_n)\cap(-\infty,X^*(x,t))\not=\emptyset.$ Therefore one can
show iteratively that if $p_{x,t}\in CB^k$, then
$$
\left|{d^{k-1}p_{x,t} \over dy^{k-1}}(y)\right|\leq \Cd
(z_n-y_n),\left|{d^{k-2}p_{x,t} \over dy^{k-2}}(y)\right|\leq \Cd
(z_n-y_n)^2,\dots,|p_{x,t}(y)|\leq \Cd (z_n-y_n)^k, \quad
y\in(y_n,z_n),
$$
and $P(y_n\leq X(x,t)\leq z_n)\leq \Cd(z_n-y_n)^{k+1}$. Comparing
this estimate with (\ref{450}) and taking $\alpha$ sufficiently small, we
obtain the needed statement. The proposition is proved.

\subsection{Smoothness of the density.}

The crucial difficulty in the proof of the smoothness of the
density is that the stochastic derivative $Y_h^\Gamma$ of the
variable $X(x,t)$, given by Theorem \ref{t31}, is not
stochastically differentiable w.r.t. $\{T_{rh}^\Gamma\}$. This
formally does not allow one to apply the standard Malliavin-type
regularity results. Moreover, the detailed analysis shows that
this difficulty is not only formal and the integration-by-parts
formula for the functionals of $X(x,t)$ (formula (\ref{421})
below) actually contains some additional "singular" terms. Below
we introduce the calculus of variations based on such
integration-by-parts formula and obtain the sufficient conditions
for the density of the law of the solution to (\ref{01}) to be
smooth.

Let us introduce some necessary constructions. We would like to
have an opportunity to divide any "portion of the jump mass" $\Pi$
into an arbitrary number of parts. Such an opportunity is
guaranteed by the following construction: we suppose that the
point measure $\nu$, correspondent to the process $U$, is in fact
a projection of another point measure $\bnu$ with a more wide
phase space and the specially constructed L\'evy measure $\bpi$.
To be precise, we suppose that the probability space is generated
by a Poisson random point measure $\bnu$ on $\ax\times\Re^{m+1}$
with the intensity measure $\lambda^1\times \bpi$, $\bpi\equiv
\Pi\times(\lambda^1|_{[0,1]})$, and $\nu$ is expressed through
$\bnu$ by
$$
\nu([0,t]\times \Gamma)=\bnu([0,t]\times \Gamma\times [0,1]),\quad
t\in\ax, \Gamma\in\pf.
$$
It is easy to see that such supposition does not restrict
generality, since for a given $\nu$ we can construct $\bnu$,
making an appropriate extension of the initial probability space.

For the "extended" random point measure $\bnu$, we will use the
terminology and constructions from Section 3. Further we denote
$\bu=(u,y)\in\Re^{m+1}$, the subsets of $\Re^{m+1}$ are denoted by
bold symbols, such as $\bg$. We also denote, by
$\bp(\cdot)=(p(\cdot),q(\cdot))$, the point process corresponding
to $\bnu$.

Given the measure $\Pi$, let us construct  the monotonously
decreasing sequence $\{\eps_n, n\in\ZZ\}$ in the following way:
$$\eps_0=1, \quad {\eps_{n+1}\over \eps_{n}}=1- (|n|+2)^{-1}, n\in \ZZ.
$$
By the construction, the sequence $\{\eps_n\}$ has  the following
properties:
$$\eps_n\downarrow 0, n\to +\infty, \quad \eps_n \to\infty,
n\to -\infty, \quad {\eps_{n+1}\over \eps_n}\to 1, n\to\infty,
\quad \sup_n{\eps_n\over \eps_{n+1}}\leq 2.$$

Denote $I_n\equiv \{u|\|u\|\in[\eps_{n+1},\eps_n)\}$. Let $t\in
\ax, \gamma\in(0,{1\over 2}),B>0$ be fixed, define the numbers
$K_n\in \NN, n\in \ZZ$ by
$$
K_n=\left[\max\left(B, 2t\Pi(I_n), {3\over
\gamma}\cdot2^{|n|-2}t^2\Pi(I_n)\right)\right]+2,
$$
where $[x]\equiv \max\{k\in \ZZ, k\leq x\}$. By the construction,
$$
K_n> B, \quad {t\over K_n}\Pi(I_n)<{1\over 2}, \quad {t^2\over
K_n}\Pi^2(I_n)<{2\gamma\over 3}\cdot 2^{-|n|}.
$$

We consider all the sets of the type  $I_n\times [{k-1\over
K_n},{k\over K_n})\subset \Re^{m+1}, k=1,\dots, K_n, n\in \ZZ,$
and enumerate them in an arbitrary way by the parameter $i\in
\NN$. The $i$-th set from this family will be denoted by
$\bg^\gamma_i$. Now, we can consider the the grid $\Gf^\gamma$ for
the random point measure $\bnu$ in the following way.

1) Every time interval $[a_i^\gamma,b_i^\gamma)$ is equal to $[0,t)$.

2) The family of sets $\{\bg^\gamma_i\}$ is the one constructed
before.

3) For every $i$, the function $h_i^\gamma$ has the form
$(\eps_n^{-1}\wedge 1) h$, where $n=n(i)$ is such that
$\bg^\gamma_i=I_n\times [{k-1\over K_n},{k\over K_n})$ for some
$k$. The function $h\in C^\infty(\Re)$ is such that $Jh=0$ outside
$(0,t)$, $Jh>0$ inside $(0,t)$ and $Jh=1$ on $(\beta,t-\beta)$,
where the constant  $\beta\in(0,{1\over2})$ will be determined
later on.

Denote $\Xi^\gamma=\bigcap_{i}\{\#\{\tau\in
\Df\cap[0,1)|\bp(\tau)\in\bg_i^\gamma\}\leq 1\}$. All the
variables $\#\{\tau\in \Df\cap[0,1)|\bp(\tau)\in\bg_i^\gamma\}$
are independent Poissonian variables with the intensities
$\lambda_i\equiv {t\over K_{n(i)}}\Pi(I_{n(i)})$.  For any
Poissonian variable $\xi$ with the intensity $\lambda$, the
inequality $P(\xi>1)\leq {\lambda^2\over 2}$ holds true. Thus
\be\label{45} P(\Xi^\gamma)
\geq\prod_{n\in\ZZ}\prod_{k=1}^{K_n}\left(1-{t^2\Pi^2(I_{n(i)})\over
2K_{n(i)}^2}\right)\geq 1-\sum_{n\in\ZZ}\sum_{k=1}^{K_n}
{t^2\Pi^2(I_{n(i)})\over 2K_{n(i)}^2}=
1-\sum_{n\in\ZZ}{t^2\Pi^2(I_{n(i)})\over 2K_{n(i)}}>
1-\sum_{n\in\ZZ}{\gamma\over 3}\cdot 2^{-|n|}=1-\gamma. \ee

Our trick is to replace the initial probability $P$ by
$$P^\gamma(\cdot)=P(\cdot|\Xi^\gamma)={P(\cdot\cap\Xi^\gamma)\over
P(\Xi^\gamma)}.
$$
 We will study firstly the distribution of $X(x,t)$
w.r.t. $P^\gamma$ and then tend $\gamma$ to $0$. The key point
here is the following analog of the classical Fourier lemma (see
\cite{malliavin} or Lemma 8.1 \cite{Iked_Wat}). Below we denote,
by $E^\gamma$, the expectation w.r.t. $P^\gamma$.

\begin{lem}\label{l_mal} Suppose that, for some $k\geq 0$, there exists
constants $\cs_1,\dots, \cs_{k+m}\in\Re^+$ such that, for every
$\gamma\in(0,{1\over 2}), F\in C_b^\infty(\Re^m), n\leq k+m,
\alpha_1,\dots,\alpha_n\in\{1,\dots,m\}$, \be\label{e_mal}
\left|E^\gamma\Bigl[{\prt\over \prt x_{\alpha_1}}\dots {\prt\over
\prt x_{\alpha_n}}F\Bigr](X(t))\right|\leq \cs_n\sup_x|F(x)|. \ee
Then $P(X(t)\in dx)=p(x)dx$ with $p\in CB^k(\Re^m)$.
\end{lem}

\demo The Fourier lemma provides that
 $P^\gamma(X(t)\in dx)=p^\gamma(x)dx$  with $p^\gamma\in
 CB^k(\Re^m)$ and
$$
\left\|{\prt\over \prt x_{\alpha_1}}\dots {\prt\over \prt
x_{\alpha_k}}p^\gamma\right\|_{L_\infty}\leq \cs_m, \quad
\alpha_1,\dots,\alpha_k\in\{1,\dots,m\}.
$$
Due to (\ref{45}), the measures $P^\gamma(X(t)\in \cdot)$ weakly
converge to $P(X(t)\in \cdot)$, $\gamma\to 0+$. This implies the
needed statement. The lemma is proved.

Thus, our further goal is to construct the grids $\Gf^\gamma$ in
the special way in order to provide (\ref{e_mal}) to hold true.
Let us mention that $\Xi^\gamma$ is invariant w.r.t.
$T_{rh_i^\gamma}^{\sbg_i^\gamma}$, and $\1_{\Xi^\gamma}\in
W_p^1(\Gf^\gamma,\Re)$ with $D_p^{\Gf^\gamma} \1_{\Xi^\gamma}=0$,
$p\in(1,+\infty)$. This means that the "censoring" operation
$P\mapsto P^\gamma$ described above is adjusted with the
differential structure. On the other hand, the following
proposition shows that $P^\gamma$ is some kind of a mixture of the
Bernoulli and uniform distributions. Such a measure appears to be
more convenient for us to deal with, than the initial Poisson one.
Below we omit the superscript  $\gamma$ in the notation for
$\Xi^\gamma$ and $\bg_i^\gamma$ (but not for $P^\gamma$).

\begin{prop}\label{p41} Denote
$$
\Xi^0_i=\{\{\tau\in
\Df\cap[0,1)|\bp(\tau)\in\bg_i\}=\emptyset\},\quad
\Xi_i^1=\Xi\backslash\Xi^0_i=\{\{\tau\in
\Df\cap[0,1)|\bp(\tau)\in\bg_i\}=\{\tau_i\}\}.
$$
Then

a) $P^\gamma(\Xi^0_i)={1\over 1+\lambda_i},
P^\gamma(\Xi_i^1)={\lambda_i\over 1+\lambda_i}$;

b) the distribution of $\tau_i$ w.r.t.
$P(\cdot|\Xi_i^1)=P^\gamma(\cdot| \Xi_i^1)$ coincides with the
uniform distribution on $[0,t]$ (below we denote this distribution
by $\lambda_t^1$);

c) the distribution of $\bp(\tau_i)$  w.r.t. $P(\cdot|\Xi_i^1)$ is
equal to $\mu_i(\cdot)={\sbpi(\cdot\cap \sbg_i)\over \sbpi(\sbg_i)}$;

d) for any  $i_1,\dots, i_k\in \NN, i_j\not=i_l, j\not=l$,
$a_1,\dots, a_k\in \{0,1\}$ the sets
$\Xi_{i_1}^{a_1},\dots,\Xi_{i_k}^{a_k}$ are jointly independent
w.r.t. $P^\gamma$;

e) for any  $i_1,\dots, i_k\in \NN, i_j\not=i_l, j\not=l$, the
variables
$\tau_{i_1},\dots,\tau_{i_k},\bp(\tau_{i_1}),\dots,\bp(\tau_{i_k})
$ are jointly independent w.r.t. $P^\gamma(\cdot| \bigcap_{j=1}^k
\Xi_{i_j}^1)$.
\end{prop}

{\it Proof.} Denote by $\bnu_i, i\in \NN$ the point measures,
defined on $\ax\times \Re^{m+1}$ by
$$
\bnu_i([0, s]\times \Delta)=\bnu([0,s]\times (\bg_i\cap \Delta)),
\quad s\in \ax, \Delta\in \Bf(\Re^{m+1}), i\in \NN.
$$
The following facts (valid for any disjoint family of the sets
$\bg_i, i\in \NN$ with $\bpi(\bg_i)<+\infty)$  are well known in
the theory of the L\'evy processes:

(i) the measures $\{\bnu_i, i\in \NN\}$ are jointly independent;

(ii) for every $i\in \NN$, the domain of the point process
$\bp_i$, correspondent to $\bnu_i$, is a.s. locally finite;

(iii) for every $i\in \NN$ the sequences $\{\tau_1^i, \tau_2^i,
\dots\}$ and $\{\xi_1^i,\xi_2^i,\dots\}$ of the points of the
domain of $\bp_i$ (enumerated increasingly) and correspondent
values of $\bp_i$ are independent;

(iv) the process $N^i_s\equiv \#\{k|\tau_k\leq s\}$ is a Poisson
process with the intensity $\bpi(\bg_i)$;

(v) $\{\xi_k^i, k\geq 1\}$ are i.i.d. random vectors in
$\Re^{m+1}$ with their common distribution equal to
${\sbpi(\cdot\cap \sbg_i)\over \sbpi(\sbg_i)}$.

For any $i\in \NN$ the sets $\Xi_i^a,a=1,2$ belong to
$\sigma(N^i_\cdot)$, and $\Xi\equiv \Xi^\gamma=\bigcap_{i\in
\NN}[\Xi_i^0\cup \Xi_i^1]$. Using this, one can easily verify that
(i) -- (v) imply statements c),d),e). For a Poisson process $N$
with the intensity $\lambda$, we have that
$$
P(N_t=0)=e^{-t\lambda }, \quad P(N_t=1)=(t\lambda) e^{-t\lambda
}\Longrightarrow P(N_t=0| N_t\leq 1)={1\over 1+t\lambda },\quad
P(N_t=1| N_t\leq 1)={t\lambda \over 1+t\lambda}.
$$
This provides the statement a). At last, for the moment $\tau$ of
the first jump of the process $N$, the following relation holds:
$$
P(\tau\leq s|N_t=1)=\left[(t\lambda) e^{-t\lambda}\right]^{-1}
P(N_s=1,N_t=1)=\left[(t\lambda) e^{-t\lambda}\right]^{-1}
\left\{(s\lambda)e^{-s\lambda}\cdot
e^{-(t-s)\lambda}\right\}={s\over t}, \quad s\in [0,t].
$$
This provides the statement b). The proposition is proved.

 Consider the space
$\varOmega=\prod_{i\in\NN}(\{0,1\}\times [0,t]\times \Re^{m+1})$
with the measure $M=\prod_{i\in\NN}\Bigl(\hbox{Be}
({\lambda_i\over 1+\lambda_i})\times \lambda_t^1\times
\mu_i\Bigr)$, here Be$(p)$ denotes the Bernoulli distribution with
$P(1)=p$. For every $\varpi=(\theta_i,s_i,\bu_i, i\in\NN)\in
\varOmega$, we define the configuration $\omega=\omega(\varpi)$ in
the following way: it consists of the points $\{(s_i,\bu_i)\in
[0,t]\times\Re^{m+1}, i\in I^1\}$, where $I^1=\{i|\theta_i=1\}$.
Let the function $f\in L_0(\Omega,\Ff,P)$ depend only on the
values of the point measure on $[0,t]\times\Re^{m+1}$, define
$\tilde f(\varpi)=f(\omega(\varpi))$. Since $P^\gamma\ll P$,
Proposition \ref{p41} implies that the map $f\to \tilde f$ is well
defined, i.e. taking a $P$-modification of $f$ we obtain the
function that is $M$-a.s. equal to $\tilde f$. Further we omit the
sign $\tilde{\phantom{f}}$ and denote by $f$ both the function
defined on  $\Omega$ and its image defined on $\varOmega$.

Denote $\vO_i^j=\{\theta_i=j\}, j=0,1$ and $M_i^j(\cdot)=M(\cdot|
\vO_i^j),j=0,1$. Denote
$$
\se f=\int_{\vO}f(\varpi)\,M(d\varpi),\quad \se_i^j
f=\int_{\vO}f(\varpi)\,M_i^j(d\varpi),\quad j=0,1.
$$
Define the transformation $\eps^{s,\bu}_i: \vO\to \vO_i^1,
(s,\bu)\in [0,t]\times\bg_i$ in the following way: it does not
change all coordinates with indices not equal to $i$ and replaces
$(\theta_i,s_i,\bu_i)$ by $(1,s,\bu)$. The restriction of this
operator on $\vO_i^0$ is just an appropriate version of the
operator $\eps^+_{(s,\bu)}$ adding the point $(s,\bu)$ to
the configuration (see \cite{picard}). Denote, by the same symbol
$\eps^{s,\bu}_i$, the transformation
$$
L_0(\Omega,P)\ni f(\cdot)\mapsto f(\eps^{s,\bu}_i\cdot )\in
L_0(\Omega_i^0, M_i^0).
$$
Recall  (see the discussion in \cite{picard}, Section 1) that, for
two different modifications $f_1,f_2$ of $f\in L_0(\Omega,P)$, the
functions $\eps^{s,\bu}_i f_1, \eps^{s,\bu}_i f_2$ may be not
equal to $M_i^0$ a.s. for the given $(s,\bu)$. But the set
$\{(s,\bu):\eps^{s,\bu}_i f_1\not=\eps^{s,\bu}_i f_2\}$ has zero
$\lambda_t^1\times \mu_i$-measure. This means that the family of
the transformations $\{\eps^{s,\bu}_i,(s,\bu)\in[0,t]\times
\bg_i\}$ is well defined in the $L_0([0,t]\times
\bg_i,\lambda_t^1\times \mu_i)$ sense.

The following formula is a simple corollary of Proposition \ref{p41}
and is, in fact, the main purpose of the  construction given
above.

\begin{prop}\label{p42} For any $f\in L_1(\Omega,P),i\in \NN$,
\be\label{46} \se_i^1 f={1\over
t}\int_0^t\int_{\sbg_i}\bigl[\se_i^0 \eps^{s,\bu}_i
f\bigr]\,\mu_i(d\bu)ds.\ee
\end{prop}

Now we are going to  proceed with the proof of Theorem \ref{t17}. We
will do this in two steps.

{\it Proof of Theorem \ref{t17}: the case $m=1$.}

Consider the functionals $f=X(t)\1_{\Xi}$ (we omit the initial
value $x$ in the notation for $X(x,t)$) and $g_i=D_{h_i}^{\bg_i}
f, i\in \NN$. The latter derivative exists since $D_{h_i}^{\sbg_i}
\1_{\Xi}=0$. Due to Theorem \ref{t31}, one has
$$
g_i=Jh_i(\tau_i)\Ef_{\tau_i}^t\Bigl[a\Bigl(X(\tau_i-)+p(\tau_i)\Bigr)-a\Bigl(X(\tau_i-)\Bigr)\Bigr]
\1_{\Xi_i^1},\quad i\in\NN.
$$
  Since $\nabla a$ is
bounded, $|\Ef_{\tau_i}^t|\leq \Cd$ and
$\Bigl|a\Bigl(X(\tau_i-)+p(\tau_i)\Bigr)-a\Bigl(X(\tau_i-)\Bigr)\Bigr|\leq
\Cd |p(\tau_i)|$. We recall that $Jh_i=(\eps_{n(i)}^{-1}\wedge 1)Jh$
and $\|Jh\|_\infty<+\infty$, thus
\be\label{47}
\sum_{i\in\NN} g_i^2\leq \Cd\sum_{i\in\NN}
p^2_1(\tau_i)(\eps_{n(i)}^{-2}\wedge 1)\1_{\Xi_i^1}\leq
\Cd\sum_{i\in\NN}(1\wedge \eps_{n(i)}^{2}) \1_{\Xi_i^1}.
\ee
We have
$$
E\sum_{i\in\NN} (1\wedge
\eps_{n(i)}^{2})\1_{\Xi_i^1}=P(\Xi)\sum_{i\in\NN} (1\wedge
\eps_{n(i)}^{2}){\lambda_i\over 1+\lambda_i}<\sum_{i\in\NN}
(1\wedge \eps_{n(i)}^{2})\lambda_i=
$$
$$=\sum_{n\in \ZZ}\sum_{k=1}^{K_n}(1\wedge
\eps_{n}^{2}){t\Pi(I_{n})\over K_{n}}=t \sum_{n\in \ZZ}(1\wedge
\eps_{n}^{2})\Pi(I_{n})\leq t2^2\int_{\Re} (u^2\wedge
1)\Pi(du)<+\infty.
$$
Here we used that $\eps_n\leq 2\eps_{n+1}\leq |u|$ for $u\in I_n$.
Thus the series on the right-hand side of (\ref{47}) converges in the $L_1$
sense, $g=(g_i)\in L_2(\Omega,P,\ell_2)$ and $f\in
W_2^1(\Gf^\gamma)$ with $D^\Gf f=g$.

We put $Z=\|g\|_{\ell_2}^2\geq 0.$ For any function $F\in
C_b^\infty$, one has
$$
(Z+c)^{-1}(D^{\Gf^\gamma}F(f), g)_{\ell_2}={\sum_{i\in\NN}
F'(f)g_i^2\over \sum_{i\in\NN} g_i^2+c}\to F'(f)
\1_{\{Z>0\}},\quad c\to 0+
$$
almost surely and in every $L_p$.  We will show below  that
$\{Z>0\}=\Xi$\, almost surely.  Thus, in order to estimate
$E^\gamma F'(f)=EF'(f)\1_{\Xi^\gamma}$, it is enough to estimate
$EF'(f){g_i^2\over Z+c}$ in such a  way that is uniform in $c$ and
allows the summation over $i$.  The key point here is   the
following moment estimate. For a given $k\in \NN, i_1,\dots,i_k\in
\NN, \bu_1,\dots\bu_k\in \Re^2, s_1,\dots, s_k\in [0,t]$, we
denote \be\label{4_notation} \mathsf{E}_{i_1,\dots,
i_k}^0[\cdot]=\mathsf{E}[\cdot|\theta_{i_1}=\dots\theta_{i_k}=0],\quad
Z_{i_1,\dots,
i_k}^{\bu_1,\dots,\bu_k}(s_1,\dots,s_k)=\eps^{s_1,\bu_1}_{i_1}\dots
\eps^{s_k,\bu_k}_{i_k}\left[\sum_{i\not=i_1,\dots,i_k}g_i^2\right].
\ee
\begin{lem}\label{l43}
 Let $a\in \Kb_r$ and ${t\over 2r}{e-1\over
e}\bro_{2r}> \alpha$ for some $r\in \NN,\alpha\in[0,+\infty)$.
Then, for every $k\in \NN$, there exists $\delta>0$ such that, under
an appropriate choice of the constants $B,\beta$ in the
construction of the grids $\Gf^\gamma$,
\be\label{48}
\sup_\gamma\sup_{l\leq k}\sup_{i_1,\dots,i_l\in \NN}
\sup_{\bu_1\in \sbg_{i_1},\dots,\bu_l\in
\sbg_{i_l}}\sup_{s_1,\dots,s_l\in[0,t]}\biggl[\se_{i_1,\dots,
i_l}^0 [Z_{i_1,\dots,
i_l}^{\bu_1,\dots,\bu_l}(s_1,\dots,s_l)]^{-\alpha-\delta}\biggr]<+\infty.
\ee
\end{lem}

{\it Proof.} In order to shorten the notation, we consider only the
case $k=1$, the general case is completely analogous
 (namely, the
only change in the proof will be that the term
$B-1$ in (\ref{414}) should be replaced by $B-k$). Everywhere in the proof of the
lemma, we omit the subscript near $i, \bu, s$.

 We use the arguments that are not the simplest possible here, but appear to be appropriate
both for the case $m=1$, and for the general case considered in
Lemma \ref{l412}  below.   We return from the "censored" probability
space $(\vO_i^0, M_i^0)$ to the initial one $(\Omega, P)$ and
  provide (\ref{48}) by the arguments analogous to those used
in the proof of Theorem \ref{t11}.

We have $P(\Xi_i^0)=P(\Xi){1\over 1+\lambda_i}\geq \Cd>0$,
and thus $\se_i^0[\cdot]=[P(\Xi_i^0)]^{-1}E[\cdot\cap \Xi_i^0]\leq
\Cd E[\cdot].$
 Let us denote
$$
Z_i=\sum_{\tau_k\in \Df:\bp(\tau_k)\not\in
\bg_i}[Jh(\tau_k)(|p(\tau_k)|^{-1}\wedge 1)]^2
\Bigl(a\Bigl(X(\tau_k-)+p(\tau_k)\Bigr)-a\Bigl(X(\tau_k-)\Bigr),
[(\Ef_0^{\tau_k})^*]^{-1}v\Bigr)^2_{\Re^m}
$$
and estimate  $E [\eps^{s,\bu} Z_i]^{-\alpha-\delta}$, where
$\eps^{s,\bu}$ denotes the operator adding the point
$(s,\bu)$ to the configuration.

For $D\equiv[D(a,r)\wedge 1]$ ($D(a,r)$ is given in
Definition \ref{d15}), we have
\be\label{49} \eps^{s,\bu} Z_i\geq
\Cd\!\! \sum_{\tau_k\in \Df:\bp(\tau_k)\not\in \sbg_i}
[p(\tau_k)]^{2r}\1_{|p(\tau_k)|\leq D}\1_{\tau_k\in
[\beta,t-\beta]}, \ee here we used that $\Ef_0^\cdot$ is separated
both from $0$ and from $+\infty$ by some non-random constants.

Denote
$$A_i(\varkappa)\equiv \Bigl\{\{\tau_k\in \Df\cap
[\beta,t-\beta]:\bp(\tau_k)\not\in \bg_i, |p(\tau_k)|\leq
D,|p(\tau_k)|> \varkappa\}=\emptyset\Bigr\},\quad  \varkappa>0.$$
Due to the Chebyshev inequality, we have
$$
P(\eps^{s,\bu}Z_i< \Cd\varkappa^{2r})\leq
$$
$$ \leq\Cd E\exp\left\{-\varkappa^{-2r}\sum_{\tau_k\in \Df\cap
[\beta,t-\beta]:\bp(\tau_k)\not\in \bg_i}
[p(\tau_k)]^{2r}\1_{|p(\tau_k)|\leq D}\1_{|p(\tau_k)|\leq
\varkappa}\right\}\1_{A_i(\varkappa)}=
$$
\be\label{410} =\Cd E\prod_{\tau_k\in \Df\cap
[\beta,t-\beta]:\bp(\tau_k)\not\in \sbg_i, |p(\tau_k)|\leq
D}\Psi(\varkappa,\tau_k),\ee where
$$
\Psi(\varkappa,\tau_k)=\begin{cases}
\exp\left\{-\varkappa^{-2r}[p(\tau_k)]^{2r}\right\},&
|p(\tau_k)|\leq
\varkappa,\\
0, &|p(\tau_k)| > \varkappa.
\end{cases}
$$
Denote
$$
\phi(\varkappa)=E\prod_{\tau_k\in \Df\cap
[\beta,(t-\beta)]:\bp(\tau_k)\not\in \sbg_i, |p(\tau_k)|\leq
D}\Psi(\varkappa,\tau_k),
$$
$$\phi^n(\varkappa)=E\prod_{\tau_k\in
\Df^n\cap [\beta,(t-\beta)]:\bp(\tau_k)\not\in \sbg_i,
|p(\tau_k)|\leq D}\Psi(\varkappa,\tau_k),
$$
we have $\phi^n\to \phi, n\to+\infty$. We may assume that the
(locally finite) set $\{\tau_k\}=\Df^n$ is ordered in the natural
monotonous way. Denote, by $\Pi_i$, the projection on the first
coordinate of the measure $\bpi_i(\cdot)=\bpi(\cdot\backslash
\bg_i)$. For every $k$ ($\tau_k\in\Df^n$) the value of the jump
$p(\tau_k)$ is independent of
$\Ff_k\equiv\Ff_{\tau_k-}\vee\sigma(\tau_k)$, and the distribution
of the jump is equal to $[\Pi_i(\{|u|\geq {1\over
n}\})]^{-1}\cdot\Pi_i(\cdot\cap \{|u|\geq {1\over n}\})$.
 Take $\varkappa<D, n>{1\over \varkappa}$ and denote $\gamma_i^n=\Pi_i(\{|u|\geq {1\over
n}\})$. Then
$$
E[\Psi(\varkappa,\tau_k)\1_{\bp(\tau_k)\not\in \bg_i,
|p(\tau_k)|\leq D}|\Ff_k]\leq[\gamma_i^n]^{-1}\int_{{1\over n}\leq
|u|\leq
\varkappa}\exp\{-\varkappa^{-2r}u^{2r}\}\Pi_i(du)=
$$
\be\label{411} =1-[\gamma_i^n]^{-1}\left\{\Pi_i(|u|\geq
\varkappa)+\int_{{1\over n}\leq |u|\leq
\varkappa}[1-\exp\{-\varkappa^{-2r}u^{2r}\}]\Pi_i(du)\right\}. \ee
It follows from (\ref{411}) that \be\label{412} \phi^n\leq
E\left[1-[\gamma_i^n]^{-1}\left\{\Pi_i(|u|\geq
\varkappa)+\int_{{1\over n}\leq |u|\leq
\varkappa}[1-\exp\{-\varkappa^{-2r}u^{2r}\}]\Pi_i(du)\right\}\right]^{N(n,i,D,\beta)},
\ee where $N(n,i,D,\beta)=\#\{k|\tau_k\in[\beta,t-\beta], {1\over
n}\leq |p(\tau_k)|\leq D\}$ is the Poissonian random variable with
its intensity equal to
$\gamma(n,i,D,\beta)\equiv(t-2\beta)\Pi_i({1\over n}\leq |u|\leq
D).$ We have ${\gamma(n,i,D,\beta)\over \gamma_n}\to
(t-2\beta)$, and thus (\ref{412}) implies that \be\label{413}
\phi(\varkappa)\leq \lim\sup_{n\to +\infty} \phi^n(\varkappa)\leq
\exp\left\{-(t-2\beta)\left[\Pi_i(|u|\geq \varkappa)+\int_{|u|\leq
\varkappa}[1-\exp\{-\varkappa^{-2r}u^{2r}\}]\Pi_i(du)\right]\right\}.\ee
It follows from the construction of the grid that \be\label{414}
\Pi_i(\cdot)\geq {B-1\over B}\Pi(\cdot),\ee because while one cell
$\bg_i$ is removed,  the "row" with the number $n(i)$  still
contains $K_{n(i)}-1$ "copies" of this cell. Then, using
(\ref{413}) and the elementary inequality $1-\exp(-x)\geq
{e-1\over e}x, x\in[0,1],$ we obtain that
$$
\phi(\varkappa)\leq\exp\left\{-(t-2\beta){e-1\over e}{B-1\over
B}\left[\Pi(|u|>\varkappa)+\varkappa^{-2r}\int_{|u|\leq
\varkappa}u^{2r}\Pi(du)\right]\right\}=
$$
$$
=\exp\left\{-(t-2\beta){e-1\over e}{B-1\over B}\ln\Bigl[{1\over
\varkappa}\Bigr]\rho_{2r}(\varkappa)\right\}=\varkappa^{(t-2\beta){e-1\over
e}{B-1\over B}\rho_{2r}(\varkappa)},
$$
and consequently, for $\kappa=\varkappa^{2r}$,
\be\label{415}
P(\eps^{s,\bu}Z_i<\Cd\kappa)\leq \Cd \kappa^{{t-2\beta\over
2r}{e-1\over e}{B-1\over B}\rho_{2r}(\kappa^{1\over
2r})}.\ee
Now we put $\delta={1\over 2}[t{e-1\over 2er}\bro_{2r}-\alpha]$
and choose $\beta$ and  $B$ in such a way that ${t-2\beta
(B-1)\over 2rB}{e-1\over e}\bro_{2r}>\alpha+{4\delta\over 3}$.
Then (\ref{415}) implies that
$$
\lim_{\kap\to 0+}\sup_{\gamma, i,s,\bu} \kap^{-\alpha-\delta}
P(\eps^{s,\bu}Z_i<\kappa)<+\infty,
$$
that proves the needed statement. The lemma is proved.

Let $i$ be fixed. We can write
$$
EF'(f){g_i^2\over Z+c}=P(\Xi){\lambda_i\over 1+\lambda_i}
\mathsf{E}_i^1 F'(f){g_i^2\over Z+c},
$$
since $g_i=0$ on $\Omega\backslash \Xi_i^1$. Using (\ref{46}), we write
$$
\mathsf{E}_i^1 F'(f){g_i^2\over Z+c}={1\over t}
\int_{\sbg_i}\se_i^0\left[\int_0^tF'(f_{\bu}(s)){g_{i,\bu}^2(s)\over
Z_\bu(s)+c}ds\right]\mu_i(d\bu)=
$$
\be\label{416} ={1\over t} \int_{\sbg_i}\se_i^0\left[\int_0^t
F'(f_{\bu}(s)) G_{i,\bu}(s) Y_{i,\bu,c}(s) ds\right]\mu_i(d\bu),
\ee where the following notation is used:
$f_{\bu}(s)=\eps^{s,\bu}_i f,g_{i,\bu}(s)=\eps^{s,\bu}_{i}
g_i,G_{i,\bu}(s)=[Jh_i(s)]^{-1}g_{i,\bu}(s),Z_{\bu}(s)=\eps^{s,\bu}_i
Z,
 Y_{i,\bu,c}(s)=[Jh_i(s)]^2\cdot  {G_{i,\bu}(s)\over Z_{\bu}(s)+c}.$

We are going to write the integration-by-parts formula for the
integral w.r.t. $ds$ in (\ref{416}). In order to do this, we need some
notation and preliminary results.

\begin{dfn}\label{d44} The function $f:\ax\mapsto \Re$ is called to
belong to the class ACPD (absolutely continuous + purely
discontinuous) if $f\in BV_{loc}(\ax)$ and there exists the
function $g\in L_{1,loc}(\ax)$ such that
$$
f(r-)-f(0+)=\int_0^r g(s)\,ds+ \sum_{s\in(0,r)}[f(s+)-f(s-)],
\quad r\in\ax.
$$
The function $g$ $\lambda^1$-a.s. coincides with the derivative of
$f$. Therefore we denote $g=f'={\prt\over \prt s} f$.

If $f$ belongs to ACPD and is continuous, then it is absolutely
continuous. In this case, we say that it belongs to the class
$AC$.
\end{dfn}

The following statement is quite standard, and therefore we just
outline its proof.

\begin{prop}\label{p45} Let $f_1,\dots, f_m$ belong to the class
ACPD. Then, for every $F\in C^1(\Re^m)$, the function
$F(f_1,\dots, f_m)$ belongs to the same class with
$$
[F(f_1,\dots, f_m)]'=\sum_{k=1}^m F'_k(f_1,\dots, f_m) f_k',
$$
$$
[F(f_1,\dots, f_m)](s+)-[F(f_1,\dots, f_m)](s-)=[F(f_1(s+),\dots,
f_m(s+))]-[F(f_1(s-),\dots, f_m(s-))]
$$
(the first equality should be understood in the  $\lambda^1$-a.s.
sense).
\end{prop}

{\it Sketch of the proof.} The statement of the proposition is
trivial when $f_1,\dots, f_m$ have only finite family
$\{s_1<\dots<s_m\}$  of the points of discontinuity, and belong to
the class $C^1$ on every interval $[s_k, s_{k+1}], k=1,\dots,
m-1$. If the functions  $f_1,\dots, f_m$ belong to the class AC on
every interval $[s_k, s_{k+1}]$, then one can prove the needed
statement for them, approximating them, together with their
derivatives, in $L_1$ sense on these intervals by smooth
functions,  and then passing to the limit. In the general case,
one should first approximate every function $f_j$ by the functions
$f^\eps_j, \eps>0$, defined by the relations
$$
f_j^\eps(r-)-f^\eps_j(0+)=\int_0^r f_j'(s)\,ds+\sum_{s\in
(0,r)}[f(s+)-f(s-)]\1_{|f(s+)-f(s-)|>\eps},
$$
and then again pass to the limit as $\eps\to 0+$.

\begin{prop}\label{p46} There exist the modifications of the
processes $X(\cdot), \Ef_0^\cdot$ such that, for any $\bu\in\bg_i$,

1) for every $r\in[0,t]$, the function $s\mapsto \eps_i^{s,\bu}
X(r)$ belongs to $AC$ with its derivative equal to
$$
{\prt\over \prt s}\eps_i^{s,\bu} X(r)= (\eps_i^{s,\bu}\Ef_{s}^r)
\Bigl[a\Bigl(X(s-)+u\Bigr)-a\Bigl(X(s-)\Bigr)\Bigr]\1_{[0,r]}(s),\quad
s\in[0,t];
$$

2) for every $r\in [0,t]$, the function $s\mapsto \eps_i^{s,\bu}
\Ef_{0}^r$ belongs to $AC$ with
$$
{\prt\over \prt s} \eps_i^{s,\bu} \Ef_{0}^r=
(\eps_i^{s,\bu}\Ef_0^r)
\Bigl[a\Bigl(X(s-)+u\Bigr)-a\Bigl(X(s-)\Bigr)\Bigr] \int_s^r
a''(\eps_i^{s,\bu}X(z))\,dz\cdot \1_{[0,r]}(s),\quad s\in[0,t];
$$

3) the function $s\mapsto \Ef_0^s$ belongs to $AC$ with
${\prt\over\prt s}\Ef_0^s=a'(X(s-))\Ef_0^s$;

4) the function $s\mapsto X(s-)$ belongs to ACPD with
${\prt\over\prt s} X(s-)=\tilde a (X(s-)),\quad \tilde a(x)\equiv
a(x)-\int_{|u|\leq 1}u\,\Pi(du). $ The set of jumps of this
function coincides with $\{s_j| \theta_j=1\}$, and the value of
the jump at the point $s_j$ is equal to $u_j$.
\end{prop}

\demo Statements 3),4)  follow straightforwardly from the
construction of $X(\cdot), \Ef_0^\cdot$. Statement 1) is just
the statement of Theorem \ref{t31} reformulated to the other form.
Statement 2) follows from the considerations completely
analogous to those given in the proof of Theorem \ref{t31}. The
proposition is proved.

 As a corollary, we obtain the following statement.

\begin{prop}\label{p47}
 There exist the modifications of the
functions $f, g_i$ such that, everywhere on $\varOmega_i^0$ for
every $\bu\in\bg_i$, the function $Y_{i,\bu,c}(\cdot)$ belongs to
the class ACPD, and
 the following integration-by-parts formula holds:
 \be\label{417}
\int_0^t F'(f_{\bu}(s))G_{i,\bu}(s) Y_{i,\bu,c}(s) ds=-\int_0^t
F(f_{\bu}(s))[Y_{i,\bu,c}]'(s)\,ds-\sum_{s\in[0,t]}F(f_{\bu}(s))
\bigl[Y_{i,\bu,c}(s+)-Y_{i,\bu,c}(s-)\bigr]. \ee
\end{prop}

{\it Proof.} It follows from Proposition \ref{46} that
$[f_\bu]'=G_{i,\bu}$ belongs to ACPD with \be\label{418}
|G_{i,\bu}(s)|\leq C_2(a)|u|,\quad |[G_{i,\bu}]'(s)|\leq
C_2(a)|u|(1+|X(s-)|),\quad |G_{i,\bu}(s_j+)-G_{i,\bu}(s_j-)|\leq
C_2(a)|u||u_j|,\quad j\not=i,\ee where the constant $C_2(a)$
depends only on $\|a'\|_\infty,\|a''\|_\infty$. Analogously, for
$j\not=i$, the function
$$
s\mapsto G_{i,j,\bu}(s)=(\eps_i^{s,\bu}\Ef_{s_j}^t)
\Bigl[a\Bigl(\eps_i^{s,\bu}X(s_j-)+u_j\Bigr)-a\Bigl(\eps_i^{s,\bu}X(s_j-)\Bigr)\Bigr]
\1_{\{\theta_j=1\}}
$$
belongs to AC with
\be\label{419}
|G_{i,j,\bu}(s)|\leq C_2(a) |u_j|,\quad |[G_{i,j,\bu}]'(s)|\leq
C_2(a) |u||u_j|. \ee
 Then the  function
$\sum_{j\not=i}[Jh_j(s_j)]^2G_{i,j,\bu}^2(\cdot) $ belongs to AC
with its  derivative dominated by $|u|
(C_2(a)\cdot\|Jh\|_\infty)^2\xi,$ where \be\label{420}
\xi=2\sum_{j\in\NN} (u_j^2\wedge 1)\1_{\{\theta_j=1\}}\in
\bigcap_{p>1}L_p(\vO, M).\ee Therefore the function
$$
Z_{\bu}(\cdot)=G_{i,\bu}^2(\cdot)
[Jh_i(\cdot)]^2+\sum_{j\not=i}[Jh_j(s_j)]^2G_{i,j,\bu}^2(\cdot)
$$ belongs to the class ACPD. At last,  $Z_{\bu}(s)+c\geq c>0,$ and,
applying Proposition \ref{p45} with $F\in C^1(\Re^{2})$  such that
$F(x,y)={x\over y}$ for $x\in\Re, y>c$, we obtain  that
$Y_{i,\bu,c}$ belongs to ACPD. Applying once again Proposition
\ref{p45}, we obtain (\ref{417}) (we use here that
$Jh_i(0)=Jh_i(t)=0$, and thus
$Y_{i,\bu,c}(0+)=Y_{i,\bu,c}(t-)=0$). Proposition is proved.

Estimates (\ref{418}),(\ref{419}) straightforwardly imply the following
estimates for $[Y_{i,\bu,c}]'$ and
$\bigl[Y_{i,\bu,c}(s)-Y_{i,\bu,c}(s-)\bigr]$ that do not involve
$c$.

\begin{prop}\label{p48}
 1) For every $s\in[0,t]$,
$$|[Y_{i,\bu,c}]'(s)|\leq  2(|u|\wedge 1)
(C_2(a)\cdot\|Jh\|_\infty)^2(\xi+1+|X(s-)|)[\Sigma_{i}^\bu(s)]^{-2}.$$

2) For every $j\not=i$,
$$
\bigl|Y_{i,\bu,c}(s_j)-Y_{i,\bu,c}(s_j-)\bigr|\leq
C_2(a)(|u|\wedge 1)(|u_j|\wedge 1)[\Sigma_{i}^\bu(s_j)]^{-1}
\1_{\{\theta_j=1\}}.
$$

The constant  $C_2(a)$ depends only on
$\|a'\|_\infty,\|a''\|_\infty$.
\end{prop}

Now we can write down the integration-by-parts formula for  the
functionals of $f=X(t)$  on $(\Xi^\gamma, P^\gamma)$. Denote, by
$E^\gamma$, the expectation w.r.t. $P^\gamma$ and put $Y_{i,\bu}\equiv
Y_{i,\bu,0}$.

\begin{lem}\label{l49} Let $a\in \Kb_r$ and ${t\over 2r}{e-1\over
e}\bro_{2r}> 2$ for some $r\in \NN$. Suppose that the constants $\beta,B$
in the construction of the grids $\Gf^\gamma$ are given by Lemma \ref{l43}
 with $\alpha=2, k=2$.
Then \be\label{421} E^\gamma F'(f)=-{1\over t}\sum_{i\in
\NN}{\lambda_i\over
\lambda_i+1}\int_{\sbg_i}\se_i^0\left[\int_0^1F(f_\bu(s))Y_{i,\bu}'(s)\,
ds+\sum_{j\not=i}
F(f_\bu(s_j))[Y_{i,\bu}(s_j+)-Y_{i,\bu}(s_j-)]\right]\mu_i(d\bu)
  \ee
for every $F\in C_b^1(\Re)$, and \be\label{422} \sup_\gamma
\sum_{i\in \NN}{\lambda_i\over
\lambda_i+1}\int_{\sbg_i}\se_i^0\left[\int_0^1|Y_{i,\bu}'(s)|\,
ds+\sum_{j\not=i}
|Y_{i,\bu}(s_j+)-Y_{i,\bu}(s_j-)|\right]\mu_i(d\bu)<+\infty. \ee
\end{lem}

{\it Remark.} Two terms on the right-hand side of (\ref{421}) can
be naturally interpreted as the integrals of $F(f)$ w.r.t. some
signed measures. Estimate (\ref{422}) shows that these measures
have finite total variation. The essential point here is that the
 second term in the integral w.r.t. the measure that is, in fact, singular
w.r.t. the initial probability. This motivates us to call (\ref{421}) the
\emph{singular type} integration-by-parts formula.

{\it Proof.}  We have
 $\sup_{s\in[0,t]}E|\xi+1+X(s-)|^p<+\infty$ for
every $p<+\infty$, thus statement 1) of Proposition \ref{p48} and
Lemma \ref{l43} provide that
$$
\int_{\sbg_i}\se_i^0\int_0^1|Y_{i,\bu}'(s)|\, ds\mu_i(d\bu) \leq
\Cd(\eps_{n(i)}\wedge 1),\quad i\in \NN, c>0.
$$
Next, we  use statement 2) of Proposition \ref{p46} and Proposition
\ref{p42} to write
$$
\int_{\sbg_i}\se_i^0\sum_{j\not=i}
\bigl|Y_{i,\bu,c}(s_j)-Y_{i,\bu,c}(s_j-)\bigr| ds\mu_i(d\bu)\leq
$$
$$
\leq{C_2(a)\over t}(\eps_{n(i)}\wedge
1)\int_{\sbg_i}\sum_{j\not=i}\int_0^t\int_{\sbg_j}\se_{i,j}^0{\lambda_j\over
1+\lambda_j}(\eps_{n(j)}\wedge 1) [\Sigma_{i,j}^{\bu,\tilde
\bu}(s,\tilde s)]^{-1}\mu_j(d\tilde\bu) d\tilde s\mu_i(d\bu)\leq
$$
$$
\leq{C_2(a)}(\eps_{n(i)}\wedge 1)\left[\sum_j \lambda_j
(\eps_{n(j)}\wedge 1)\right]\left[\se_{i,j}^0\sup_{i,j,\bu,\tilde
\bu, s,\tilde s} [\Sigma_{i,j}^{\bu,\tilde \bu}(s,\tilde
s)]^{-1}\right]\leq \Cd (\eps_{n(i)}\wedge 1), \quad i\in \NN, c>0
$$
(see (\ref{4_notation} for the notation $Z^{\bu, \tilde
bu}_{i,j}$). In the last inequality, we used Lemma \ref{l43} and
the fact that, due to condition (\ref{13}),
$$
\sum_j \lambda_j (\eps_{n(j)}\wedge 1)\leq 2\int_{\Re}(|u|\wedge
1)\Pi(du)<+\infty.
$$
Once again, we use $\sum_i \lambda_i (\eps_{n(i)}\wedge
1)<+\infty$ and deduce (\ref{421}) and (\ref{422}).  The lemma is proved.

{\it Remark.} The explicit estimates given above show that there
exists a constant $\cs_1<+\infty$ such that, {\it for every} grid
$\Gf^\gamma$ constructed in the way given above {\it for any}
$\gamma>0$, the expression on the left-hand side of (\ref{422}) is
dominated by $\cs_1$.

The last thing we need to complete the proof of Theorem \ref{t17} is
to iterate (\ref{421}) in order to provide an estimate for
$EF^{(n)}(f)$ in the terms of $\sup_x|F(x)|$ ($F^{(n)}$ denotes
the $n$-th derivative of $F$). The essential point here is that
the measure $M_i^0$ is also the product measure and possesses
the constructions given before for the measure $M$.

Let us  rewrite (\ref{421}) to the form that is convenient to the
further iterative procedure. For a given $n$, we denote, by $\Theta(n)$,
the family of all partitions $\theta=(\theta_1,\dots,
\theta_r)$ of the set $\{1,\dots, n\}$ into non-overlapping parts
(for instance, $\Lambda(2)$ contains two partitions
$(\{1\},\{2\})$ and $(\{1,2\})$). Denote also, by $\NN_d^n$, the set
of all vectors $i_1,\dots, i_n$ with all coordinates not equal to one
another. For a given $\bar i\equiv (i_1,\dots, i_n)\in
\NN_d^n$, $\bar u\equiv (\bu_1,\dots, \bu_n), \bar s\equiv
(s^1,\dots, s^r)$, and a partition
$\theta=(\theta_1=\{\theta_1^1,\dots,
\theta_1^{l_1}\},\dots,\theta_r=\{\theta_r^1,\dots,
\theta_r^{l_r}\})\in\Lambda(n)$, we denote
$$
\eps_{\bar i,\theta}^{\bar s,\bar
\bu}=[\eps_{i_{\theta_1^1}}^{s_1,\bu_1}\circ
\eps_{i_{\theta_1^2}}^{s_1,\bu_2}\circ\dots\circ
\eps_{i_{\theta_1^{l_1}}}^{s_1,\bu_{l_1}}]\circ[
\eps_{i_{\theta_2^1}}^{s_2,\bu_{l_1+1}}\circ\dots\circ
\eps_{i_{\theta_2^{l_2}}}^{s_2,\bu_{l_1+l_2}}] \circ\dots \circ[
\eps_{i_{\theta_r^1}}^{s_r,\bu_{n-l_r+1}}\circ\dots\eps_{i_{\theta_r^{l_r}}}^{s_r,\bu_{n}}].
$$
Now, using the statement analogous to the one of  Proposition \ref{p42},
applied to $\se_i^0$ instead of $\se$, we can write (\ref{421}) in the
form
\be\label{423}
E^\gamma F'(f)=\sum_{\theta\in\Theta(2)}\sum_{\bar i\in
\NN_d^2}\int_{\sbg_{i_1}\times
\sbg_{i_2}}\int_{[0,t]^{r(\theta)}}\!\se_{i_1,i_2}^0F(\eps_{\bar
i,\theta}^{\bar s,\bar \bu}f)\, Y_{\bar i,\theta}^{\bar \bu}(\bar
s)\, d\bar s\, [\mu_{i_1}\times\mu_{i_2}](d\bar \bu),
\ee
where $r(\theta)$ is the number of the components in the partition
$\theta$, and the functions  $Y_{\bar i,\theta}^{\bar \bu}$ are
either a derivative or a jump of the function $Y_{i,\bu}$ (in the
notation of (\ref{421})) multiplied by
$-{\lambda_{i_1}\lambda_{i_2}\over
t^2(\lambda_{i_1}+1)(\lambda_{i_2}+1)}$ or
$-{\lambda_{i_1}\lambda_{i_2}\over
t(\lambda_{i_1}+1)(\lambda_{i_2}+1)}$, respectively.

Take $F\in C_b^2(\Re)$ and apply (\ref{423}) to $\tilde F=F'.$ Then the
terms of the type $\se_{i_1,i_2}^0F'(\eps_{\bar i,\theta}^{\bar
s,\bar \bu}f)\, Y_{\bar i,\theta}^{\bar \bu}(\bar s)$ occur on the
right-hand side of (\ref{423}). For every such a term, we write
$$
\se_{i_1,i_2}^0F'(\eps_{\bar i,\theta}^{\bar s,\bar \bu}f)\,
Y_{\bar i,\theta}^{\bar \bu}(\bar
s)=\sum_{i\not=i_1,i_2}\se_{i_1,i_2}^0  F'(\eps_{\bar
i,\theta}^{\bar s,\bar \bu}f)\, Y_{\bar i,\theta}^{\bar \bu}(\bar
s)\cdot  \eps_{\bar i,\theta}^{\bar s,\bar \bu} \left[{g_i^2\over
\Sigma_{i_1,i_2}}\right]=
$$
\be\label{424}
=\sum_{i\not=i_1,i_2}{\lambda_i\over
t(\lambda_i+1)}\int_{\sbg_i}\int_0^t \eps_i^{s,\bu}[F'(\eps_{\bar
i,\theta}^{\bar s,\bar \bu}f)\, Y_{\bar i,\theta}^{\bar \bu}(\bar
s)]\cdot  \eps_{\bar i,\theta}^{\bar s,\bar \bu} \left[{g_i^2\over
\Sigma_{i_1,i_2}}\right]\,ds\mu_i(d\bu),
\ee
where $\Sigma_{i_1,i_2}=\sum_{i\not=i_1,i_2} g^2_i$. From
Proposition \ref{p46}, we get that the function $s\mapsto
\eps_i^{s,\bu}\eps_{\bar i,\theta}^{\bar s,\bar \bu}f$ belongs to
AC with
$$
{\prt\over \prt s} \eps_i^{s,\bu}\eps_{\bar i,\theta}^{\bar s,\bar
\bu}f =[Jh_i(s)]^{-1} \eps_i^{s,\bu}  \eps_{\bar i,\theta}^{\bar
s,\bar \bu} {g_i}.
$$
The function $s\mapsto [Jh_i(s)] \eps_i^{s,\bu}  \eps_{\bar
i,\theta}^{\bar s,\bar \bu}\left[{g_i\over \Sigma_{i_1,i_2}
}\right]$ belongs to ACPD with its derivative and jumps satisfying
the estimates analogous to those given in Proposition \ref{p48},
but with $\Sigma_i^\bu$ replaced by $\eps_i^{s,\bu}  \eps_{\bar
i,\theta}^{\bar s,\bar \bu} \Sigma_{i,i_1,i_2, i}$, where
$\Sigma_{i,i_1,i_2}=\sum_{j\not=i,i_1,i_2} g^2_j$. At last, using
Proposition \ref{p46} and the explicit form of $Y_{\bar
i,\theta}^{\bar \bu}(\bar s)$, one can verify that the function
$s\mapsto \eps_i^{s,\bu} Y_{\bar i,\theta}^{\bar \bu}(\bar s)$
also belongs to ACPD  with its derivative and jumps dominated by
$$
\Cd\xi\cdot\lambda_{i}\lambda_{i_1}\lambda_{i_2}(\eps_{n(i)}\wedge
1)(\eps_{n(i_1)}\wedge 1)(\eps_{n(i_2)}\wedge
1)[\Sigma_{i,i_1,i_2}]^{-3},
$$
where the constant $\Cd$ depends only on the coefficient $a$, and
the variable $\xi$ belongs to $\cap_p L_p$. This means that, under
an appropriate moment condition imposed on
$[\Sigma_{i,i_1,i_2}]^{-3}$, we can write the integration-by-parts
formula on the right-hand side of (\ref{424}) and obtain the analog
of (\ref{423}) with $E^\gamma F''(f)$ on the left-hand side. Let us
formulate this statement  for the derivative of an
arbitrary order. For a given $\bar i\in \NN_d^n$, we denote $\bg_{\bar
i}=\bg_{i_1}\times\dots\times \bg_{i_n}, \mu_{\bar i}=
\mu_{i_1}\times\dots\times\mu_{i_n}, \se_{\bar
i}^0=\se_{i_1,\dots, i_n},\Sigma_{\bar i}=\sum_{i\not\in\bar i}g_i^2$.

\begin{lem}\label{l410} Let $n\in \NN$ be fixed,
 $a\in \Kb_r$ and ${t\over 2r}{e-1\over
e}\bro_{2r}> 2n$ for some $r\in \NN$. Suppose that the constants
$\beta,B$ in the construction of the grids $\Gf^\gamma$ are given
by Lemma \ref{l43} with
 $\alpha=2n, k=n$.

 Then there exists a set of the functions
$\{Y_{\bar i, \theta}^{\bar \bu}:[0,t]^{r(\theta)}\to \Re, \bar
i\in \NN_d^{2n}, \theta\in\Theta(2n),\bar\bu\in[\Re^2]^{2n} \}$
such that \be\label{425} E^\gamma
F^{(n)}(f)=\sum_{\theta\in\Theta(2n)}\sum_{\bar i\in
\NN_d^{2n}}\int_{\sbg_{\bar i}}\int_{[0,t]^{r(\theta)}}\!\se_{\bar
i}^0F(\eps_{\bar i,\theta}^{\bar s,\bar \bu}f)\, Y_{\bar
i,\theta}^{\bar \bu}(\bar s)\, d\bar s\, \mu_{\bar i}(d\bar \bu),
\ee and \be\label{426} \int_{\sbg_{\bar
i}}\int_{[0,t]^{r(\theta)}}\!\se_{\bar i}^0|Y_{\bar
i,\theta}^{\bar \bu}(\bar s)|\, d\bar s\, \mu_{\bar i}(d\bar
\bu)\leq C(n,\delta)\lambda_{i_1}\dots \lambda_{i_{2n}}
\eps_{n(i_1)}\dots
\eps_{n(i_{2n})}(1+\eps_{n(i_1)})^{M(n,\delta)}\dots(1+\eps_{n(i_{2n})})^{M(n,\delta)},
\ee $\bar i\in \NN_d^{2n},
 \theta\in\Theta(2n),$
 where $C(n,\delta), M(n,\delta)$ are some constants depending only on $n$ and the number
 $\delta$ given by Lemma \ref{l43}.
\end{lem}
{\it Proof.} The iterative procedure described before shows how
one can deduce formula (\ref{425}) for a given $n$ from the
same formula for $n-1$: one should take one term in (\ref{425})
 and write
down the formula analogous to (\ref{424}) for it. This explains how the coefficients
$Y_{\bar i, \theta}^{\bar \bu}$ of the order $n$ (i.e., with $\bar i\in \NN_d^{2n}$)
are constructed: one should take all $\bar i\in \NN_{d}^{2n-2}$, $i\not\in \bar i$
 and calculate the derivative and
the jump part of the function $s\mapsto {\lambda_i[Jh_i(s)]^2\over
t(\lambda_i+1)}\eps_{i}^{s,\bu} [Y_{\bar i, \theta}^{\bar \bu}
{g_i\over \sum_{j\not\in\bar i}g_j^2}]$. All such functions are
exactly the new coefficients $Y_{\bar i, \theta}^{\bar \bu}$. Such
a description of the family $\{Y_{\bar i, \theta}^{\bar \bu}\}$
allows one to rewrite it to the form $\{Y_{\bar i, \theta}^{\bar
\bu}(\bar s)= {H_{\bar i, \theta}^{\bar \bu}(\bar s) [\eps_{\bar
i, \theta}^{\bar s,\bar \bu} \Sigma_{\bar i}]^{-2n}}\}$, where the
functions  $\{H_{\bar i, \theta}^{\bar \bu}\}$ are defined
 iteratively. The power $2n$ here appears, since the power of the
  denominator increases by $1$ twice on one step of the induction:  the first time when the term ${g_i\over
\sum_{j\not\in\bar i}g_j^2}$ is added, and the second one when
either a derivative or the jump part is calculated.

Using the explicit expressions for the derivatives and jumps of
the processes $X(\cdot),\Ef_0^\cdot$ (which the functions
$\{g_i\}$, and thus the functions  $\{H_{\bar i, \theta}^{\bar
\bu}\}$,  are expressed through)  one can deduce by induction on
$n$ that, for every index sets $\bar i\in\NN_d^n,\bar j$ with $\bar
l=\bar i\cup\bar j={l_1,\dots, l_N},$ for every ordered sets
$p=(p_1,\dots,p_k)\subset\bar l^k, o=(o_1,\dots,o_k)\in\{0,1\}^k$,
 \be\label{427a} |\prt_{p_1}^{o_1}\dots \prt_{p_k}^{o_k}
 \eps^{\bar s^1,\bar \bu^1}_{\bar j\backslash \bar i} H_{\bar i, \theta}^{\bar \bu}(\bar
s)|\leq \Cd\eps^{\bar s^1,\bar \bu^1}_{\bar j\backslash \bar
i}\eps_{\bar i,\theta}^{\bar s,\bar \bu} [1+\max_{s\leq
t}|X(s)|]^{M(N,k)} \lambda_{l_1}\dots \lambda_{l_N}
(\eps_{n(l_1)}\wedge 1)\dots (\eps_{n(l_N)}\wedge 1), \ee $\bar
i\in \NN_d^{2n},
 \theta\in\Theta(2n),\bar s\in[0,t]^{r(\theta)},\bar u\in\bg_{\bar
 i},\bar s^1\in [0,t]^{N-n},\bar u\in\bg_{\bar j\backslash \bar
 i},
 $
 where $\prt_p^0$ denotes the derivative w.r.t. the variable
with the number $p$, $\prt_p^1$ denotes the jump w.r.t. the
same variable, $\eps^{\bar s^1,\bar \bu^1}_{\bar j\backslash \bar
i}\equiv \eps^{\bar s^1,\bar \bu^1}_{\bar j\backslash \bar
i,\theta_*}$ with $\theta_*=(\{1\},\dots,\{N-n\})$. We do not need
estimate (\ref{427a}) in its full generality, we only need the
partial case $\bar j=\bar i, k=0$. In this case, we have the
estimate
 \be\label{427} |H_{\bar i, \theta}^{\bar \bu}(\bar
s)|\leq \eps_{\bar i,\theta}^{\bar s,\bar \bu}  [1+\max_{s\leq
t}|X(s)|]^{M(n)} \lambda_{i_1}\dots \lambda_{i_{2n}}
(\eps_{n(i_1)}\wedge 1)\dots (\eps_{n(i_{2n})}\wedge 1),\ee $\bar
i\in \NN_d^{2n},
 \theta\in\Theta(2n),\bar s\in[0,t]^{r(\theta)},\bar u\in\bg_{\bar i},
 $
 where $M(n)$ is some constant. Note that estimate (\ref{427}) is
 not well designed to be proved by induction on $n$, while
 (\ref{427a}) is; this was the only reason for us to write firstly
 estimate (\ref{427a}).
  Now, using Lemma \ref{l43}, we obtain
$$
\int_{\sbg_{\bar i}}\int_{[0,t]^{r(\theta)}}\!\se_{\bar
i}^0|Y_{\bar i,\theta}^{\bar \bu}(\bar s)|\, d\bar s\, \mu_{\bar
i}(d\bar \bu)\leq\left[
 \int_{\sbg_{\bar i}}\int_{[0,t]^{r(\theta)}}\!\se_{\bar i}^0
 \eps_{\bar i,\theta}^{\bar s,\bar
 \bu}[1+\max_{s\leq
t}|X(s)|]^{M(n)(2n+\delta)\over \delta}d\bar s\mu_{\bar
 i}(d\bu)\right]^{\delta\over 2n+\delta} \times
$$
$$
\times \lambda_{i_1}\dots \lambda_{i_{2n}} (\eps_{n(i_1)}\wedge
1)\dots (\eps_{n(i_{2n})}\wedge 1).
 $$
 Since $\nabla a$ is bounded and
 $\int_{\{|u|>1\}}|u|^p\Pi(du)<+\infty$ for every $p$, there exists
 such a constant $\tilde C(n)$ that
$$
\se_{\bar i}^0
 \eps_{\bar i,\theta}^{\bar s,\bar
 \bu}[1+\max_{s\leq
t}|X(s)|]^{M(n)(2n+\delta)\over \delta}\leq \left[\tilde
C(n)(1+\|u_1\|^{M(n)})\dots(1+\|u_{2n}\|^{M(n)})\right]^{2n+\delta\over
\delta}.
$$
This provides (\ref{426}). The lemma is proved.

Now we can complete the proof of Theorem \ref{t17} in the case
$m=1$. We apply Lemma \ref{l410} for $n\leq k+1$. Equality
(\ref{425}) and estimate (\ref{426}) immediately imply that
(\ref{e_mal}) holds true.  Thus the needed statement holds true
due to Lemma \ref{l_mal}. The proof is complete.

{\it Proof of Theorem \ref{t17}: the case $m>1$.} All the
technique, that is necessary for the proof of Theorem
\ref{t17} in the general case, was already introduced in the proof
of the case $m=1$. Our aim now is to adapt this technique to the
multidimensional situation.

Again, denote $f=X(t)\1_\Xi$,  $g_i=D_{h_i}^{\sbg_i}f$,
now $f, g_i$ are the random vectors in $\Re^m$. Considerations analogous to those given
after estimate (\ref{47}) show that $g=(g_i)\in L_2(\Omega,P,\Re^m\otimes \ell_2)$ and $f\in
W_2^1(\Gf^\gamma,\Re^m)$ with $D^\Gf f=g$. We put
$$
Z=\sum_{i}g_i\otimes g_i,\quad Z_{\bar i}=\sum_{i\not\in \bar i}g_i\otimes g_i,\quad \bar i\in
\NN_d^n,n\geq 1,
$$
$Z$ is the Malliavin matrix for the vector $f$. We can write down
the estimate analogous to (\ref{47}) for $\|Z\|_{\Re^{m^2}}$
 and then prove (for instance, calculating the Fourier transform
of the right-hand side and then estimating its
derivatives of all the orders) that $\|Z\|_{\Re^{m^2}}\in \cap_pL_p$.

We use the notation $\alpha\in \{1,\dots, m\}$ and
$\ba=(\alpha_1,\dots,\alpha_n)\in \{1,\dots, m\}^n$ for the
indices and multiindices, $\pa\equiv{\prt\over \prt x_\alpha},
\pba\equiv {\prt\over \prt x_{\alpha_1}}\dots {\prt\over \prt
x_{\alpha_n}}.$ Let us write down the analogs of (4.21) and
(4.25). First, we do this formally, without taking care of
the terms involved in the corresponding integration-by-parts
formula to belong to $L_1$. The necessary moment estimates will be
given later on, in the second part of the proof.

Denote $Y_{i,\bu}(s)=Jh_i(s)[\eps_{i}^{s,\bu}Z]^{-1}g_i,
s\in[0,t], \bu\in\bg_i,i\in \NN,$ and let $Y_{i,\bu}^\alpha$
denote the $\alpha$-th component of the vector $Y_{i,\bu}.$ Using
Proposition \ref{p45} and an appropriate analog of
Proposition \ref{p47}, one can obtain the following analog of
the integration-by-parts formula (\ref{421}): \be\label{428}
E^\gamma [\pa F](f)=-{1\over t}\sum_{i\in \NN}{\lambda_i\over
\lambda_i+1}\int_{\Gamma_i}\se_i^0\left[\int_0^1F(\eps_{i}^{s,\bu}
f)[Y_{i,\bu}^\alpha]'(s)\, ds+\sum_{j\not=i} F(\eps_{i}^{s,\bu}
f)[Y_{i,\bu}^\alpha(s_j+)-Y_{i,\bu}^\alpha
(s_j-)]\right]\mu_i(d\bu),
  \ee
for every $F\in C_b^1(\Re^m)$ and $\alpha\in \{1,\dots, m\}$. One
can rewrite (\ref{428}) to the form analogous to (\ref{423}) and then
iterate this formula in the way described before the formulation
of Lemma \ref{l410}. The inverse matrix $Z^{-1}$ can be expressed in the
form $[\det Z]^{-1} Q$, where the elements of the matrix $Q$ (the
\emph{cofactor matrix} for $Z$) are certain polynomials of the
elements of $Z$.  At last, for every $\bar i_1\subset \bar i_2$,
$\det Z_{\bar i_2}\leq \det Z_{\bar i_1}$. Summarizing all these
considerations, we can formulate the following statement.

\begin{prop}\label{p411} For every $F\in
C_b^n(\Re^m),n\geq 1, \ba\in\{1,\dots,m\}^n$, \be\label{429}
E^\gamma [\pba F](f)=\sum_{\theta\in\Theta(2n)}\sum_{\bar i\in
\NN_d^{2n}}\int_{\sbg_{\bar i}}\int_{[0,t]^{r(\theta)}}\!\se_{\bar
i}^0F(\eps_{\bar i,\theta}^{\bar s,\bar \bu}f)\, Y_{\bar
i,\theta}^{\bar \bu,\sba}(\bar s)\, d\bar s\, \mu_{\bar i}(d\bar
\bu).\ee Here the family $\{ Y_{\bar i, \theta}^{\bar \bu,
\sba}(\bar s)\}$ possesses the point-wise representation
$\{Y_{\bar i, \theta}^{\bar \bu, \sba}(\bar s) = {H_{\bar i,
\theta}^{\bar \bu,\sba}(\bar s) [\eps_{\bar i, \theta}^{\bar
s,\bar \bu}\det Z_{\bar i}]^{-2n}}\}$ with the functions
$\{H_{\bar i, \theta}^{\bar \bu,\sba}\}$ estimated by
\be\label{430} |H_{\bar i, \theta}^{\bar \bu}(\bar s)|\leq
\eps_{\bar i,\theta}^{\bar s,\bar \bu}  [1+\max_{s\leq
t}\|X(s)\|]^{M(n)} \lambda_{i_1}\dots \lambda_{i_{2n}}
(\eps_{n(i_1)}\wedge 1)\dots (\eps_{n(i_{2n})}\wedge 1),\ee $ \bar
i\in \NN_d^{2n},
 \theta\in\Theta(2n),\bar s\in[0,t]^{r(\theta)},\bar u\in\bg_{\bar i},
 $
 where the constant $M(n)$ depends only on $n$.
\end{prop}

Equality (\ref{429}) is now nothing more than the formal
expression, since the variables $Y_{\bar i,\theta}^{\bar
\bu,\sba}$ may not belong to $L_1$. However, estimate (\ref{430})
allows one to separate the case where this equality becomes
meaningful and rigorous.

\begin{cor} Suppose that the grids $\Gf^\gamma$ were constructed in such a way
 that, for some $n\in \NN,\delta>0$,
 \be\label{431}
 \sup_\gamma \sup_{l\leq 2n}\sup_{\bar i\in \NN_d^{l}}\sup_{\bar \bu\in\sbg_{\bar i}}\sup_{\bar
 s\in [0,t]^{l}}\se_{\bar i}^0[\eps_{i_1}^{s_1,\bu_1}\dots \eps_{i_{l}}^{s_{l},
 \bu_{l}}\det Z_{\bar i}]^{-2n-\delta}.
 \ee
 Then (\ref{429}) holds true with
 \be\label{432}
\sup_\gamma \sum_{\theta\in\Theta(2n)}\sum_{\bar i\in
\NN_d^{2n}}\int_{\sbg_{\bar i}}\int_{[0,t]^{r(\theta)}}\!\se_{\bar
i}^0| Y_{\bar i,\theta}^{\bar \bu,\sba}(\bar s)|\, d\bar s\,
\mu_{\bar i}(d\bar \bu)=\cs_n<+\infty. \ee
\end{cor}
Thus, the only essential fact, that it is left to prove, is the
following multidimensional analog of Lemma \ref{l43}.

\begin{lem}\label{l412} Let $a\in \Kb_r$ and ${t\over 2r}{e-1\over
e}\bro_{2r}> (\alpha+4)m-4$ for some $r\in
\NN,\alpha\in[0,+\infty)$. Then, for every $k\in \NN$, there exists
$\delta>0$ such that, under an appropriate choice of the constants
$B,\beta$ in the construction of the grids $\Gf^\gamma$,
$$
\sup_\gamma\sup_{l\leq k}\sup_{i_1,\dots,i_l\in \NN}
\sup_{\bu_1\in \sbg_{i_1},\dots,\bu_l\in
\sbg_{i_l}}\sup_{s_1,\dots,s_l\in[0,t]}\biggl[\se_{i_1,\dots,
i_l}^0 [\eps_{i_1}^{s_1,\bu_1}\dots \eps_{i_{l}}^{s_{l},
 \bu_{l}}\det Z_{(i_1,\dots,
i_l)}]^{-\alpha-\delta}\biggr]<+\infty.
$$
\end{lem}
\demo We consider only the case $k=1$, the general case is
completely analogous.  We have  $g_i=\Ef_0^t q_i$, $
q_i\equiv Jh_i(\tau_i)[\Ef_0^{\tau_i}]^{-1}
\Bigl[a\Bigl(X(\tau_i-)+p(\tau_i)\Bigr)-a\Bigl(X(\tau_i-)\Bigr)\Bigr]\1_{\Xi_i^\gamma}$.
Define
$$
Q=\sum_{i}q_i\otimes q_i,\quad Q_{i}=\sum_{j\not=i}q_j\otimes q_j,
$$
then $Z=\Ef_0^t\cdot Q\cdot [\Ef_0^t]^*, Z_i=\Ef_0^t\cdot Q_i\cdot
[\Ef_0^t]^*$. Since $\nabla a$ is bounded, $|\det \Ef_0^t|$ is
separated from $0$  by some non-random constant (see Proposition
\ref{p52} below for the explicit estimate). Thus, in order to prove the
statement of the lemma for $k=1$, it is enough to prove that
\be\label{433}
\sup_\gamma\sup_{i\in \NN} \sup_{\bu\in
\sbg_{i}}\sup_{s\in[0,t]}\biggl[\se_{i}^0 [\eps_{i}^{s,\bu} \det
Q_i]^{-\alpha-\delta}\biggr]<+\infty.
\ee
The calculations given in the proof of Lemma 1
\cite{kom_takeuchi_simpl} provide that, in order to verify  (\ref{433}),
it is enough to prove that
\be\label{434}
\sup_\gamma\sup_{i\in \NN} \sup_{\bu\in
\sbg_{i}}\sup_{s\in[0,t]}\sup_{v:\|v\|=1}\biggl[\se_{i}^0
[(\eps_{i}^{s,\bu}
Q_iv,v)_{\Re^m}]^{4-m(\alpha+4)-\delta}\biggr]<+\infty.
\ee
We do this analogously to the proof of Lemma \ref{l43}.  Let us
return from the "censored" probability space $(\vO_i^0, M_i^0)$ to
the initial one $(\Omega, P)$ and estimate $E [(\eps^{s,\bu}
Q_iv,v)]^{4-m(\alpha+4)-\delta}$, where $\eps^{s,\bu}$ denotes the
operator adding the point $(s,\bu)$ to the configuration.
We have
$$
(Q_iv,v)_{\Re^m}=\!\!\!\!\sum_{\tau_k\in \Df:\bp(\tau_k)\not\in
\bg_i}\![Jh(\tau_k)]^2 \|(\Ef_0^{\tau_k})^*]^{-1}v\|^2
\Bigl(a\Bigl(X(\tau_k-)+p(\tau_k)\Bigr)-a\Bigl(X(\tau_k-)\Bigr),
{[(\Ef_0^{\tau_k})^*]^{-1}v\over
\|(\Ef_0^{\tau_k})^*]^{-1}v\|}\Bigr)^2_{\Re^m}[\|p(\tau_k)\|\wedge
1]^2.
$$
Since $\nabla a$ is bounded,
$$
\mathrm{essinf}\,\inf_{\|v\|=1} \|(\Ef_0^{\tau_k})^*]^{-1}v\|\geq
\Cd>0
$$
for every $k$ (see Proposition \ref{p52} below). Thus, we deduce
that, for every $\varrho\in(0,1)$, the following inequality holds
true for $D\equiv[D(a,r,\varrho)\wedge 1]$: \be\label{435}
(Q_iv,v)\geq \Cd\!\! \sum_{\tau_k\in \Df:\bp(\tau_k)\not\in \bg_i}
\Bigl(p(\tau_k), w(\tau_k) \Bigr)^{2r}_{\Re^m}\1_{p(\tau_k)\in
V(w(\tau_k),\varrho)}\1_{|p(\tau_k)|\leq D}\1_{\tau_k\in
[\beta,t-\beta]},\ee where  we denoted $w(\tau)\equiv
w(X(\tau-),{[(\Ef_0^{\tau_k})^*]^{-1}v\over
\|(\Ef_0^{\tau_k})^*]^{-1}v\|})$ (see Definition \ref{d15} for the
notation $w(\cdot,\cdot)$). Denote
$$A_i(\varkappa)\equiv \{\tau_k\in \Df\cap
[\beta,t-\beta]:\bp(\tau_k)\not\in \bg_i,p(\tau_k)\in
V(w(\tau_k),\varrho), |p(\tau_k)|\leq D,|(p(\tau_k),\eps^{s,\bu}
w(\tau_k))|> \varkappa\}=\emptyset\}, \varkappa>0.$$
 Due to the Chebyshev inequality, we have
$$
P((\eps^{s,\bu}Q_iv,v)<\Cd \varkappa^{2r})\leq
$$
$$ \leq\Cd E\exp\left\{-\varkappa^{-2r}\sum_{\tau_k\in \Df\cap
[\beta,t-\beta]:\bp(\tau_k)\not\in \bg_i} \Bigl(p(\tau_k),
\eps^{s,\bu} w(\tau_k) \Bigr)^{2r}_{\Re^m}\1_{p(\tau_k)\in
V(w(\tau_k),\varrho)}\1_{|p(\tau_k)|\leq
D}\1_{|(p(\tau_k),\eps^{s,\bu} w(\tau_k))|\leq
\varkappa}\right\}\times $$
$$
\times \1_{A_i(\varkappa)}= \Cd E\prod_{\tau_k\in \Df\cap
[\beta,t-\beta]:\bp(\tau_k)\not\in \bg_i,p(\tau_k)\in
V(w(\tau_k),\varrho), |p(\tau_k)|\leq D}\Psi(\varkappa,\tau_k),
$$
where
$$
\Psi(\varkappa,\tau_k)=\begin{cases}
\exp\left\{-\varkappa^{-2r}\Bigl(p(\tau_k), \eps^{s,\bu} w(\tau_k)
\Bigr)^{2r}_{\Re^m}\right\},& |(p(\tau_k),\eps^{s,\bu}
w(\tau_k))|\leq
\varkappa,\\
0, &|(p(\tau_k),\eps^{s,\bu} w(\tau_k))|> \varkappa.
\end{cases}
$$
One has that $p(\tau_k)$ is independent of $\Ff_k\equiv
\Ff_{\tau_k-}\vee \sigma(\tau_k)$, and $w(\tau_k)$ is
$\Ff_k$-measurable.  Thus, repeating the arguments given in the
proof of Lemma \ref{l43}, one can obtain analogously to (\ref{411} --
\ref{415}) that
$$
P((\eps^{s,\bu}Q_iv,v)_{\Re^m}< \Cd \kappa)\leq \Cd
\kappa^{{t-2\beta\over 2r}{e-1\over e}{B-1\over
B}\rho_{2r}(\kappa^{1\over 2r},\varrho)},\quad \kap\in(0,D),
$$
and, under an appropriate choice of $\varrho, \beta, B$,
$$
\lim_{\kap\to 0+}\sup_{\gamma, i,s,\bu,\|v\|=1}
\kap^{4-m(\alpha+4)-\delta}P((\eps^{s,\bu}Q_iv,v)_{\Re^m}<\kap)=0
$$
for $\delta={1\over 2}\Bigl[{t\over 2r}{e-1\over e}\bro_{2r}-
(\alpha+4)m+4\Bigr]$. The lemma is proved.

\begin{cor} Under condition of Theorem \ref{t17}, the grids
$\Gf^\gamma$ can be constructed in such a way that the integration-by-parts formula (\ref{429}) together with the moment estimate
(\ref{432}) hold true for $n\leq k+m$.
\end{cor}

This corollary immediately implies that estimates
(\ref{e_mal}) hold true for $n\leq m+k$.  Now the statement of
Theorem \ref{t17} follows from Lemma \ref{l_mal}. The
theorem is proved.

Let us make a conclusive remark. The first and second terms in
the integration-by-parts formula (\ref{421}) can be interpreted as the
"volume integral" and "surface integral", respectively, since the measure in
the second term is supported, in fact, by the
countable union of the sets $I_{i,j}\equiv \{s_i=s_j\},
i,j\in\NN$, and each of these sets can be interpreted as a "level
 set" (or  "codimension 1 set"). This is the main reason for the
calculus of variations, developed in this section, to be
substantially different from the classical (Malliavin's) form of
the stochastic calculus of variations, since, in the latter one,
the new measure is absolutely continuous w.r.t. the initial one, i.e.
in the integration-by-parts formula  only the "volume integral" is
present.

It should be mentioned that the differential structure in our
case is not like the one for the manifold with a (smooth)
boundary.
  The "surface measure"  again admits the  similar
regular structure, and the integration-by-parts formula for such a
measure generates the  "codimension 1" and "codimension 2" terms,
and so on. Thus one can informally say that the phase space of the
Poisson random measure, considered with the differential structure
generated by the time-stretching transformations, looks like  the
"infinite-dimensional complex". The crucial point in our
construction is that, on every "side of codimension $k$" of such a
complex, there still remains an infinite family of admissible
directions.

\section{Smoothness of the density of the invariant
distribution}

In this section, we consider the stationary process $\{X(s),
s\in\Re\}$ satisfying the equation \be\label{51}
X(t)-X(s)=\int_s^ta(X(r))\,dr+ U_t-U_s, \quad -\infty<s\leq
t<+\infty, \ee with the L\'evy process $U$  defined on $\Re$
by the standard construction
$$
U_t=\begin{cases}U^1_t,&t\geq 0\\
-U^2_{(-t)-}, &t<0\end{cases},
$$
where $U^1,U^2$ are two independent copies of the L\'evy process
defined on $\ax$. The coefficient $a$ is supposed to satisfy the
conditions formulated in subsection 1.3.

In order to prove the regularity of the distribution of $X(t)$
(i.e., the statement of Theorem \ref{t114}), we need to modify
slightly the constructions from the Sections 3 and 5. The reason
is that now one cannot suppose the probability space
 $(\Omega,{\Ff},P)$ to satisfy the condition
${\Ff}=\sigma(U)$. Such a supposition is, in fact, the claim to
(\ref{51}) to possess a strong solution on $\Re$ and is, in general,
a non-trivial restriction. In order to avoid such a restriction, we
make the following modifications of the constructions given above.

Denote  $H=L_2(\Re)$. Let $H_0\subset L_\infty(\Re)$ be the set of
functions with a bounded support. For $h\in H_0$ denote
$Jh(\cdot)=\int_{-\infty}^\cdot h(s)\,ds,b(h)=\sup\{r|h(v)=0,
v\leq r\}.$ For a fixed $h\in H_0$, we define the family
$\{T_h^t,t\in\Re\}$ of transformations of the axis $\Re$ by
putting $T^t_hx, x\in\Re$ equal to the value at the point $s=t$ to
the solution of the Cauchy problem (\ref{21}).

For every $h\in H_0, \Gamma\in\pf$, the transformation $T_h^\Gamma$
of the random measure $\nu$ associated with $U$ is well defined.
Since $T_h^tx\equiv x, x\leq b(h)$, the transformation
$T_h^\Gamma$ does not change the values of $\nu$ on every subset
of $(-\infty,b(h)]\times \Re^m$. Equation (\ref{51}) considered
as the Cauchy problem with $s$ fixed possesses the strong solution.
Thus, one can define the transformation $T_h^\Gamma$ of the process
$X$ in such a way that $T_h^\Gamma X(t)=X(t), t\leq b(h)$,
\be\label{52} T_h^\Gamma X(t)=X(b(h))+\int_{b(h)}^ta(T_h^\Gamma
X(r))\,dr+ T_h^\Gamma(U_t-U_{b(h)}), \quad t\geq b(h). \ee
 Like
in the proof of Theorem \ref{t17}, we enlarge the probability
space and suppose that the random measure $\nu$ associated with
the process $U$ is the projection on the first $m$  coordinates
of the random measure $\bnu$ defined on $\Re\times \Re^{m+1}$,
with its intensity measure being equal to $\lambda^1\times\bpi, \bpi\equiv
\Pi\times \left[\lambda^1|_{[0,1]}\right]$. One possible formal
way to do this is to define $(\Omega,{\Ff}, P)$ as the product of
two probability spaces $(\Omega^1,{\Ff}^1,P^1)$,
$(\Omega^2,{\Ff}^2,P^2)$, where ${\Ff}^1=\sigma(X)$, and
$\Omega^2=[0,1]^\infty, P^2=\prod_{l\in
\NN}\left[\lambda^1|_{[0,1]}\right].$ We enumerate jumps of the
process $X$ in some measurable way and put
$$
\mathbf{X}(t)=\begin{cases}(X(t),0),& X(t)=X(t-)\\
(X(t), \xi_{l(t)}),& X(t)\not=X(t-)\end{cases},
$$
where $\{\xi_l\}$ is the sequence of coordinate functionals on
$\Omega^2$ (i.e., every  $\xi_l$ has uniform distribution on
$[0,1]$), and $l(t)$ denotes the number of the jump that happens
at the moment $t$. Then $\sigma(\mathbf{X})=\Ff$, and the random
measure $\bnu$ and the corresponding point process $\bp(\cdot)$
can be constructed from $\mathbf{X}$ in the obvious way. For every
$h\in H_0, \bg\in\bpi_{fin}$, the transformation $T_h^{\sbg}$ of
the process $\mathbf{X}$ is well defined (the first coordinate $X$
is  transformed accordingly to (\ref{52}), and the transformation
of the last coordinate $\xi_{l(t)}$ is defined by the condition
$T_h^{\sbg}[l(t)]=l(T_{-h}t)$).

Further we suppose that ${\Ff}=\sigma(\bnu)$. Under this
condition, one can easily verify that an analog of Lemma \ref{l21}
holds true, and $T_h^{\sbg}$ is, in fact, the admissible
transformation of $(\Omega,{\Ff},P)$ (the explicit formula for
$p_h$ differs slightly from the one given in subsection 3.1). The
notions of the stochastic and a.s. derivatives associated with
such admissible transformations can be introduced, and then the
statement of Theorem \ref{t31} holds true for every given $h\in
H_0$ with the trivial replacements: $0$ should be replaced by
$b(h)$ and $x$ should be replaced by $X(b(h))$.

We introduce the notion of a differential grid in the same way
with Definition \ref{d26}, with $\ax$ replaced by $\Re$ and $a_i$
claimed to belong to $\Re$ (i.e., $a_i$ should not be equal to
$-\infty$) for every $i$. For every such a grid, the Sobolev
classes associated with the grid are defined in the same way with
Definition \ref{d28}.

Now let us proceed with the proof of Theorem \ref{t114}. Since $X$
is a stationary process, it is enough to study the distribution of
$X(t)$ at one fixed point $t$, say, $t=0$. For every given
$\gamma\in(0,{1\over 2})$, we construct the grid $\Gf^\gamma$ in
the way analogous to one given at the beginning of subsection 5.2.
We take the same sequence $\{\eps_n\}$ and consider all sets of
the type \be\label{53} [-N,-N+1)\times I_n \times [{k-1\over
K_{n,N}}, {k\over K_{n,N}}),\quad k=1,\dots, K_{n,N}, n\in\ZZ,
N\in \NN,  \ee recall that
$I_n=\{u|\|u\|\in[\eps_{n+1},\eps_n)\}$. We enumerate sets
(\ref{53}) by $i\in \NN$ in an arbitrary way and denote, by $n(i),
N(i)$ and $k(i)$, such numbers that the corresponding components
in the set with the number $i$ are equal to $[-N(i),-N(i)+1)$,
$I_{n(i)}$, and $[{k(i)-1\over K_{n(i),N(i)}}, {k(i)\over
K_{n(i),N(i)}})$. The numbers $K_{n,N}$ are defined for every
given $B>0, \gamma$ by
$$
K_{n,N}=\left[\max\left(B, 2t\Pi(I_n), {3\over
\gamma}\cdot2^{|n|-N-1}t^2\Pi(I_n)\right)\right]+2,
$$
 and therefore

1) $K_{n,N}\geq B$ (the constant $B$  will be determined below);

2) ${1\over K_{n,N}} \Pi(I_n)<{1\over 2}$;

3) ${1\over K_{n,N} }\Pi^2(I_n)<{2\gamma\over 3}2^{-|n|-N}$.

We define the grids $\Gf^\gamma$ by the equalities $[a_i^\gamma,
b_i^\gamma)=[-N(i),-N(i)+1)$, $\bg_i=I_{n(i)}\times[{k(i)-1\over
K_{n(i),N(i)}}, {k(i)\over K_{n(i),N(i)}}),$
$$
h_i^\gamma(s)=A^{-N(i)}(\eps_{n(i)}^{-1}\wedge 1)h(s+N(i)), \quad
s\in \Re,
$$
where $A>1$ will be determined later on, and  $h\in C^\infty$ is some
given function such that $Jh=0$ outside $[0,1]$, $Jh>0$ on $(0,1)$,
and $Jh=1$  on $[{1\over 3}, {2\over 3}]$.

The construction of the grids $\Gf^\gamma$ provides that the estimate
analogous to (\ref{45}) holds true. Next, for the function
$f=X(0)\1_{\Xi^\gamma}$, the estimate analogous to (\ref{47}) can be
written, and one can prove that $f\in \cap_pW_p^1(\Gf^\gamma,
\Re^m)$ with
$$
g_i\equiv D_{h_i}^{\sbg_i}f=
Jh_i(\tau_i)\Ef_{\tau_i}^0
\Bigl[a\Bigl(X(\tau_i-)+p(\tau_i)\Bigr)-a\Bigl(X(\tau_i-)\Bigr)\Bigr]\1_{\Xi_i^\gamma}
$$
(here and below, we use the notation from subsection 5.2).
 Repeating step-by-step the considerations given in
subsection 5.2, we obtain  the following analog of Proposition
\ref{p411}. Denote, by $S(\bar i, \theta)$ for $\bar i\in
\NN_{d}^{2n},\theta \in \Theta(2n)$, the set of $(s_1,\dots,
s_{2n})\in (-\infty,0)$ such that

1) $s_{i_k}\in [-N({i_k}), -N(i_k)+1)$, $k=1,\ldots,2n$;

2) $s_{i_k}=s_{i_j}$ for every $k,j$ such that $i_k, i_j$ belong
to the same set w.r.t. the partition $\theta$.

Denote by $\lambda_{\bar i,\theta}$ the uniform distribution on
$S(\bar i, \theta)$, i.e the surface measure on $S(\bar i,
\theta)$ considered as a subset of $\Re^{2n}$ with the Lebesgue
measure $\lambda$.

\begin{prop}\label{p51} For every $F\in
C_b^n(\Re^m),\ba\in\{1,\dots,m\}^n$ \be\label{54} E^\gamma [\pba
F](f)=\sum_{\theta\in\Theta(2n)}\sum_{\bar i\in
\NN_d^{2n}}\int_{\sbg_{\bar i}}\int_{S(\bar i, \theta)}\!\se_{\bar
i}^0F(\eps_{\bar i,\theta}^{\bar s,\bar \bu}f)\, Y_{\bar
i,\theta}^{\bar \bu,\sba}(\bar s)\, \lambda_{\bar i,\theta}(d\bar
s)\, \mu_{\bar i}(d\bar \bu).\ee The family $\{Y_{\bar i,
\theta}^{\bar \bu, \sba}(\bar s)\}$ possesses the point-wise
representation $\{Y_{\bar i, \theta}^{\bar \bu, \sba}(\bar s) =
{H_{\bar i, \theta}^{\bar \bu,\sba}(\bar s) [\eps_{\bar i,
\theta}^{\bar s,\bar \bu}\det Z_{\bar i}]^{-2n}}\}$ with the
functions $\{H_{\bar i, \theta}^{\bar \bu,\sba}\}$ estimated by
  \be\label{55} |H_{\bar i,
\theta}^{\bar \bu,\sba}(\bar s)|\leq \Cd \eps_{\bar
i,\theta}^{\bar s,\bar \bu}
[1+\max_{s\in[\min(s_1,\dots,s_n),0]}|X(s)|]^{M(n)}
\lambda_{i_1}\dots \lambda_{i_{2n}} (\eps_{n(i_1)}\wedge 1)\dots
(\eps_{n(i_{2n})}\wedge 1)A^{-\max (N(i_1),\dots,N(i_{2n}))}, \ee
$\bar i\in \NN_d^{2n},
 \theta\in\Theta(2n),\bar s\in S(\bar i,
\theta),\bar u\in\bg_{\bar
 i}$.
\end{prop}

{\it Remark.} In the estimates for $|H_{\bar i, \theta}^{\bar
\bu,\sba}(\bar s)|$ we dominate all terms of the type
$Jh_i(\tau_i), [Jh_i]'(\tau_i),\cdots$, by some constant $\Cd$
except the terms of the same type with $N(i)=\max
(N(i_1),\dots,N(i_{2n}))$. These terms are dominated by $\Cd
A^{-\max (N(i_1),\dots,N(i_{2n}))}$, that provides the term
$A^{-\max (N(i_1),\dots,N(i_{2n}))}$ in (\ref{55}).

Let us repeat the cautions made after Proposition \ref{p411}:
equality (\ref{54}) is just a formal one; in order to make it
 rigorous the proof that $Y_{\bar i,
\theta}^{\bar \bu, \sba}(\bar s)$ are integrable  w.r.t.
$\se_{\bar i}^0$ is needed. Such a proof should contain two parts:
the estimate of the moment of $\eps_{\bar i,\theta}^{\bar s,\bar
\bu} [1+\max_{s\in[\min(s_1,\dots,s_n),0]}|X(s)|]^{M(n)}$, and the
estimate of the moment of $[\eps_{\bar i, \theta}^{\bar s,\bar
\bu}\det Z_{\bar i}]^{-2n}$.

The first part of the proof is more or less standard. The variable
$\eps_{\bar i,\theta}^{\bar s,\bar \bu} X(\min(s_1,\dots,s_n)-)$
is, in fact, equal to $X(\min(s_1,\dots,s_n)-)$, and thus its
distribution w.r.t. $\se_{\bar i}^0$ is equal to the initial
invariant distribution $P^*$. This distribution was supposed in
the formulation of Theorem \ref{t114} to have all the moments.
Moreover, the gradient $\nabla a$ is globally bounded, and thus we
can deduce from the standard martingale inequalities and the
Gronwall lemma that there exists a constant
$C(a)\equiv\sup_x\|\nabla a(x)\|]$ such that, for every $p>1$,
\be\label{56} \Bigl[\se_{\bar i}^0 \eps_{\bar i,\theta}^{\bar
s,\bar \bu}
[1+\max_{s\in[\min(s_1,\dots,s_n),0]}|X(s)|]^{p}\Bigr]^{1\over
p}\leq \Cd(1+\|\eps_{n(i_1)}\|)\dots
(1+\|\eps_{n(i_{2n})}\|)e^{-C(a) \min(s_1,\dots,s_n)}
 \ee
 with the constant $\Cd$ depending on $p,n$ and the moments of $P^*.$

 The second part of the proof contains the estimate for $[\eps_{\bar i, \theta}^{\bar s,\bar
\bu}\det Z_{\bar i}]^{-2n}$, and is  yet another version of
Lemma \ref{l43}.

\begin{lem}\label{l52} Let $\Pi$ possess the wide cone condition and
$a\in \Kb_\infty$. Let $k\in \NN$  be fixed, and let the constant $B$
in the construction of the grids $\Gf^\gamma$ be taken greater
than $k$. Then, for every $\alpha>0$ under an arbitrary choice of the
constant $A$ in the construction of the grids $\Gf^\gamma$,
$$
\sup_\gamma\sup_{l\leq k}\sup_{i_1,\dots,i_l\in \NN}
\sup_{\bu_1\in \sbg_{i_1},\dots,\bu_l\in
\sbg_{i_l}}\sup_{s_1,\dots,s_l\in[0,t]}\biggl[\se_{i_1,\dots,
i_l}^0 [\eps_{i_1}^{s_1,\bu_1}\dots \eps_{i_{l}}^{s_{l},
 \bu_{l}}\det Z_{(i_1,\dots,
i_l)}]^{-\alpha}\biggr]<+\infty.
$$
\end{lem}

\demo Again we consider only the case $k=1$, the general case is
completely analogous. Like in the proofs of Lemmae
\ref{l43},\ref{l412}, we return from the "censored" probability
space $(\vO_i^0, M_i^0)$ to the initial one $(\Omega, P)$ and
estimate $E [\eps^{s,\bu} \det Z_i]^{-\alpha}$, where
$$
Z_i=\sum_{\tau_k\in \Df, \bp(\tau_k)\not\in\sbg_i}
g(\tau_k)\otimes g(\tau_k),
$$
$$
g(\tau_k)\equiv
Jh(\tau_k-[\tau_k])A^{[\tau_k]}(\|p(\tau_k)\|^{-1}\wedge
1)\Ef_{\tau_k}^0
\Bigl[a\Bigl(X(\tau_k-)+p(\tau_k)\Bigr)-a\Bigl(X(\tau_k-)\Bigr)\Bigr]
$$
Let $N$ be fixed. We denote by $\Df(i,N)$ the set of
$\tau_k\in\Df$ such that $\bp(\tau_k)\not\in\bg_i$ and $\tau_k\in
\bigcup_{r=1}^N[-r-{1\over 3},-r-{2\over
 3}]$. Let us estimate the variable
$$
\det Z_{i,N}, \quad Z_{i,N}\equiv\sum_{\tau_k\in\Df(i,N)}
g(\tau_k)\otimes g(\tau_k).
$$
It is clear that $\eps^{s,\bu}\det Z_{i,N}\leq \eps^{s,\bu} \det
Z_i$, thus the lower estimate for $\eps^{s,\bu}\det Z_{i,N}$
provides also the lower estimate for $\eps^{s,\bu} \det Z_i$.

We write the decomposition  $g(\tau_k)=\Ef_{-N}^0 q(\tau_k)$,
$$
q(\tau_k)\equiv
Jh(\tau_k-[\tau_k])A^{[\tau_k]}(\|p(\tau_k)\|^{-1}\wedge
1)[\Ef_{-N}^{\tau_k}]^{-1}
\Bigl[a\Bigl(X(\tau_k-)+p(\tau_k)\Bigr)-a\Bigl(X(\tau_k-)\Bigr)\Bigr],
$$
and define
$$
Q_{i,N}=\sum_{\tau_k\in\Df(i,N)} q(\tau_k)\otimes q(\tau_k).
$$
Then $Z_{i,N}=\Ef_{-N}^0\cdot Q_{i,N}\cdot [\Ef_{-N}^0]^*$ and
$\det Z_{i,N}=\det Q_{i,N}\cdot \Bigl[\det \Ef_{-N}^0 \Bigr]^2$.
It would be convenient for us to formulate all the estimates
concerned to $\Ef_{-N}^\cdot$ in a separate statement.

\begin{prop}\label{p52} The following estimates hold true almost surely for every $T>0$:

1) $\det \Ef_{-T}^0\geq \exp\{-mTC(a)\};$

2) $\inf_{\|v\|=1} \|(\Ef_{-T}^{s})^*]^{-1}v\|\geq
\exp\{-TC(a)\}$, $s\in[-T,0]$;

3)$\|\Ef_{-T}^0\|,\|[\Ef_{-T}^{s}]^{-1}\|\in[\exp\{-TC(a)\},\exp\{TC(a)\}]$.
\end{prop}

\demo The first estimate is implied by the representation
\be\label{57}
 \det
\Ef_{-T}^0=\exp\{\int_{-T}^0\mathrm{trace}(\nabla a(X(s)))\, ds\}.
\ee
 This representation follows from the same one for ODE's, that
is a classical fact in theory of ODE's. In order to deduce
(\ref{57}) in the framework of the equations with the L\'evy noise
one should first prove (\ref{57}) for a compound Poisson process $U$
by just applying (\ref{57}) for ODE's piecewisely and then use
an approximation procedure.

In order to deduce the second estimate, we use the equality
$$
[(\Ef_{-T}^{s})^*]^{-1}v=v-\int_{-T}^s [\nabla
a(X(r))]^*[(\Ef_{-T}^{r})^*]^{-1}v\,dr
$$
that implies that, for every $v\in S_m$, the function
$V(s)=\|[(\Ef_{-T}^{s})^*]^{-1}v\|$ satisfies the inequality
\be\label{58} V(s)\geq v-\int_{-T}^{s}C(a) V(r)\, dr,\quad s\geq
{-T}. \ee Inequality (\ref{58}) can be written in the form of the
equation
 $$
 V(s)=1+\Delta(s)-\int_{-T}^{s}C(a) V(r)\, dr
 $$
 with the condition $\Delta(s)\geq 0$, and the solution to this
 equation can be given in the form
 $$
 V(s)=\exp\{-(s+T)C(a)\}+\int_{-T}^s
 \exp\{-(s-r)C(a)\}\Delta(r)\,dr\geq
 \exp\{-(s+T)C(a)\}.
 $$
The last estimate follows from the Gronwall lemma, on the one
hand, and from the arguments given in the proof of the second
estimate, on the other hand. The proposition is proved.

One can see that the same estimates with those given made in
Proposition \ref{p52} hold true for $\eps^{s,\bu} \Ef_{-T}^\cdot$.
Due to statement 1), $\eps^{s,\bu} \det Z_{i,N}\geq \det
\eps^{s,\bu}Q_{i,N}\cdot\exp\{-2mNC(a)\}$. Let us estimate
$\eps^{s,\bu}\det Q_{i,N}$. In order to do this, we will
appropriately modify the arguments given in the proof of Lemma
1 \cite{kom_takeuchi_simpl}.

Due to the condition on $\Pi$, there exists $\varrho\in(0,1)$ such
that $\Pi(V(w,\varrho))=+\infty$ for every cone $V(w,\varrho),
w\in S_m.$ Let $a\in\Kb_r, r\in\NN$, further we denote
$D=[D(a,r,\varrho)\wedge 1]$. For every given $\Lambda>0$, there
exists $\delta=\delta(\Lambda,\varrho)$ such that
$$
\Pi\Bigl(u|u\in V(w,\varrho), \|u\|<D,
|(u,w)_{\Re^m}|>\delta\Bigr)\geq \Lambda
$$
(one can prove this using the Dini theorem analogously to
Lemma \ref{l34}).

 We take an arbitrary $v\in S_m,$ and denote by $ \Df(i,N,v,\delta)$
 the subset of  $\Df$ containing all the points $\tau_k$
 such that $\tau_k\in \bigcup_{r=1}^N[-r-{1\over 3},-r-{2\over
 3}]$,
   $\bp(\tau_k)$ does not belong to the
cell $\bg_i$, and
$$\eps^{s,\bu}\left|\left(
\Bigl[a\Bigl(X(\tau_k-)+p(\tau_k)\Bigr)-a\Bigl(X(\tau_k-)\Bigr)\Bigr],
[(\Ef_{-N}^{\tau_k})^*]^{-1}v\right)_{\Re^m}\right|\geq
D\delta^r\cdot\eps^{s,\bu}\left\|[(\Ef_{-N}^{\tau_k})^*]^{-1}v\right\|_{\Re^m}.
$$
The same arguments with those used in the proofs of Theorem
\ref{t11} and Lemmae \ref{l43}, \ref{l412} provide that $$
P(\Df(i,N,v,\delta)=\emptyset) \leq \exp\{-{B-1\over 3B}\Lambda
N\}. $$

For every $v,\tilde v\in S_m$, due to statement 3) of Proposition
\ref{p52} we get for every $\tau_k\in \Df(i,N,v,\delta)$
$$
\eps^{s,\bu}\left|
\Bigl(a\Bigl(X(\tau_k-)+p(\tau_k)\Bigr)-a\Bigl(X(\tau_k-)\Bigr),
[(\Ef_{-N}^{\tau_k})^*]^{-1}v\Bigr)_{\Re^m}-\right.
$$
$$
\left.-\Bigl(a\Bigl(X(\tau_k-)+p(\tau_k)\Bigr)-a\Bigl(X(\tau_k-)\Bigr),
[(\Ef_{-N}^{\tau_k})^*]^{-1}\tilde v\Bigr)_{\Re^m}\right|\leq
2C(a)\|[(\Ef_{-N}^{\tau_k})^*]^{-1}\|\|v-\tilde v\|\leq
2C(a)e^{C(a)N}\|v-\tilde v\|.
$$
Let us choose the vectors $v_1,\dots,v_{\exp[2C(a)N]}$ on the
sphere $S_m$ in such a way that, for every $v\in S_m$, $\inf_{l\leq
\exp[2C(a)N]} \|v-v_l\|\leq \exp[-(2C(a)+{1\over 2})N]$ (one can
do this for $N$ large enough). Consider the event
$$
\Omega(i,N,\delta)=\bigcap_{l=1}^{\exp[2C(a)N]}\{\Df(i,N,v_l,\delta)\not=\emptyset\},\quad
P(\Omega(i,N,\delta))\geq 1-\exp\left\{\Bigl[2C(a)-{B-1\over
3B}\Lambda \Bigr]N\right\}.
$$
Take $v\in S_m$ and $l\leq \exp[2C(a)N]$ such that $\|v-v_l\|\leq
\exp[-(2C(a)+{1\over 2})N]$. Then, for every $\omega\in
\Omega(i,N,\delta)$, there exists $\tau_k\in
\Df(i,N,v_l,\delta)\subset\Df(i,N)$. For such $\tau_k$, we have
$$
\eps^{s,\bu}\left|\Bigl(a\Bigl(X(\tau_k-)+p(\tau_k)\Bigr)-a\Bigl(X(\tau_k-)\Bigr),
[(\Ef_{-N}^{\tau_k})^*]^{-1}v\Bigr)_{\Re^m}\right|\geq
$$
$$
\geq\eps^{s,\bu}\left|\Bigl(a\Bigl(X(\tau_k-)+p(\tau_k)\Bigr)-a\Bigl(X(\tau_k-)\Bigr),
[(\Ef_{-N}^{\tau_k})^*]^{-1}v_l\Bigr)_{\Re^m}\right|-2C(a)\exp[-(C(a)+{1\over
2})N]\geq
$$
\be\label{59} \geq
D\delta^r\cdot\eps^{s,\bu}\left\|[(\Ef_{-N}^{\tau_k})^*]^{-1}v_l\right\|_{\Re^m}-2C(a)\exp[-(C(a)+{1\over
2})N]\geq \left(D\delta^r-2C(a)\exp\Bigl[-{N\over
2}\Bigr]\right)\exp[-C(a)N] \ee (the last inequality in (\ref{59})
holds true due to statement 2) of Proposition \ref{p52}).
Take $N$ large enough for $D\delta^r-2C(a)\exp\Bigl[-{N\over
2}\Bigr]\geq {1\over 2}D\delta^r$. Then, due to the construction of the
grid $\Gf^\gamma$ and inequality (\ref{59}) for every $v\in S_m$,
we have the estimate
$$
\eps^{s,\bu}(Q_{i,N}v,v)_{\Re^m}\geq \!\sum_{\tau_k\in \Df(i,N)}
A^{-2N}\eps^{s,\bu}\Bigl(a\Bigl(X(\tau_k-)+p(\tau_k)\Bigr)-a\Bigl(X(\tau_k-)\Bigr),
[(\Ef_{-N}^{\tau_k})^*]^{-1}v\Bigr)^2_{\Re^m}\geq
$$
$$\geq {1\over 4} D^2\delta^{2r} A^{-2N}\cdot
e^{-2C(a)N}\1_{\Omega(i,N,\delta)}.$$

Thus, for every $\omega\in\Omega(i,N,\delta)$, we have the estimate
\be\label{510} \eps^{s,\bu}\det Z_{i,N}\geq
[\eps^{s,\bu}\inf_{v\in S_m}(Q_{i,N}v,v)]^{-m}e^{-2mC(a)N}\geq
{1\over 4^m} D^{2m}\delta^{2rm}\exp[-C(a,A,m)N],\ee
$C(a,A,m)\equiv 2m[\ln A+2C(a)]$.

At last, take $\Lambda$ large enough for ${B-1\over 3B}\Lambda-
2C(a)>(\alpha+1) C(a,A,m)$ and consider the sequence $t_N={1\over
4^m} D^{2m}\delta^{2rm}\exp[-C(a,A,m)N], N\geq 1$ (recall that
$\delta$ is defined by $\Lambda$).  Then (\ref{510}) provides that,
for $N$ large enough,
$$
P(\eps^{s,\bu}\det Z_{i}\leq t_N)\leq P(\eps^{s,\bu}\det
Z_{i,N}\leq t_N)\leq 1-P(\Omega(i,N,\delta))\leq \Cd
[t_N]^{\alpha+1}.
$$
Since $t_N\to 0+$ with $\lim\sup_{N}{t_N\over t_{N+1}}<+\infty$,
this completes the proof of the lemma. The lemma is proved.
\begin{cor} Let $N$ be fixed. Then
under conditions of Theorem \ref{t114} one can construct the grids
$\Gf^\gamma$
 in such a way that, for every
$\ba\in\{1,\dots,m\}^n, n\leq N$, the integration-by-parts formula
(\ref{54}) holds true with \be\label{511} \sup_\gamma
\sum_{\theta\in\Theta(2n)}\sum_{\bar i\in
\NN_d^{2n}}\int_{\sbg_{\bar i}}\int_{S(\bar i, \theta)}\!\se_{\bar
i}^0|Y_{\bar i,\theta}^{\bar \bu,\sba}(\bar s)|\, \lambda_{\bar
i,\theta}(d\bar s)\, \mu_{\bar i}(d\bar \bu)=\cs_n<+\infty.\ee
\end{cor}

\demo In estimate (\ref{56}), the term $\min(s_1,\dots,
s_{2n})$ can be replaced by $\max(N(i_1,\dots,N(i_{2n}))$. Now,
let us take the constant $A$ in the construction of the grid to be equal
to $2e^{C(a)}$. Then Lemma \ref{l52} and estimates
(\ref{55}),(\ref{56}), together with the H\"older inequality,
provide that
$$
 \int_{\sbg_{\bar i}}\int_{S(\bar i, \theta)}\!\se_{\bar
i}^0|Y_{\bar i,\theta}^{\bar \bu,\sba}(\bar s)|\, \lambda_{\bar
i,\theta}(d\bar s)\, \mu_{\bar i}(d\bar \bu)\leq
$$
$$
\leq\Cd \lambda_{i_1}\dots \lambda_{i_{2n}} (\eps_{n(i_1)}\wedge
1)\dots (\eps_{n(i_{2n})}\wedge 1)\Cd(1+\|\eps_{n(i_1)}\|)\dots
(1+\|\eps_{n(i_{2n})}\|)2^{-\max (N(i_1),\dots,N(i_{2n}))},
$$
$\bar i\in \NN_d^{2n}, \theta\in\Theta(2n).$ Taking the sum over
$\bar i,\theta$ we obtain (\ref{511}).

{\it End of the proof of Theorem \ref{t114}.} The corollary
given above implies that, for every given $k\in\NN$, one can
construct the grids $\Gf^\gamma$ in such a way that estimates
(\ref{e_mal}) hold true for every $n\leq k+m$. Thus, due to
Lemma \ref{l_mal}, $P^*(dy)=p^*(y)dy$ with $p^*\in \bigcap_k
CB^k(\Re^m)=C_b^\infty(\Re^m)$. The theorem is proved.

\subsection*{Acknowledgement.} The technique exposed in Section 6 and the result
of Theorem \ref{t114} were motivated by the question of A.A.
Dorogovtsev, who asked the author whether the stochastic calculus
of variations, described in the Section 3, can provide the
smoothness of the invariant density. It was a quite new idea for
the author that the differential properties of the invariant
distribution can be essentially different from those of the
distribution of the solution to the Cauchy problem. The author
would like to express his gratitude to A.A. Dorogovtsev for the
question that opened for the author a new and fruitful research
field.

\end{document}